\documentclass[onefignum,onetabnum, final]{siamart171218}
\usepackage{amsfonts}
\usepackage{amssymb}

\usepackage{subcaption}
\usepackage{bbm}
\usepackage{bm}
\usepackage{esvect}
\usepackage{booktabs}
\usepackage{multirow}
\usepackage{color}
\usepackage{cite}
\usepackage{algorithmic}
\Crefname{ALC@unique}{Line}{Lines} 

\usepackage{lipsum}
\usepackage{amsfonts}
\usepackage{graphicx}
\usepackage{epstopdf}
\usepackage{algorithmic}
\ifpdf
  \DeclareGraphicsExtensions{.eps,.pdf,.png,.jpg}
\else
  \DeclareGraphicsExtensions{.eps}
\fi

\usepackage{enumitem}
\setlist[enumerate]{leftmargin=.5in}
\setlist[itemize]{leftmargin=.5in}

\newsiamremark{remark}{Remark}
\newsiamremark{hypothesis}{Hypothesis}
\crefname{hypothesis}{Hypothesis}{Hypotheses}
\newsiamthm{claim}{Claim}

\headers{Statistical Learning for Probability-Constrained Stochastic Optimal Control}{A. Balata, M. Ludkovski, A. Maheshwari and J. Palczewski}

\title{Statistical Learning for Probability-Constrained Stochastic Optimal Control}

\author{Alessandro Balata\thanks{School of Mathematics, University of Leeds, Woodhouse Lane, Leeds LS2 9JT, United Kingdom
  (\email{alessandro.balata@live.com}; \email{J.Palczewski@leeds.ac.uk}).}
\and Michael Ludkovski\thanks{Department of Statistics and Applied Probability, University of California, Santa Barbara, CA 93106 USA
  (\email{ludkovski@pstat.ucsb.edu}; \email{maditya0310@gmail.com}).} \funding{Michael Ludkovski and Aditya Maheshwari were supported in part by NSF grant DMS-1736439. Alessandro Balata was supported by the Natural Environment Research Council Doctoral training partnership Leeds-York}
  \corresponding{Aditya Maheshwari} 
\and Aditya Maheshwari$^\dag$ 
\and Jan Palczewski$^*$}

\usepackage{amsopn}

\ifpdf
\hypersetup{
  pdftitle={},
  pdfauthor={A. Balata, M. Ludkovski, A. Maheshwari and J. Palczewski}
}
\fi

\newcommand{\cA}{\mathcal A}
\newcommand{\cC}{\mathcal C}
\newcommand{\cL}{\mathcal L}

\newcommand{\cU}{\mathcal U}
\newcommand{\CD}{\mathcal D}
\newcommand{\cH}{\mathcal H}
\newcommand{\cF}{\mathcal F}

\newcommand{\cX}{\mathcal X}

\newcommand{\cG}{\mathcal G}
\newcommand{\cK}{\mathcal K}
\newcommand{\cW}{\mathcal W}
\newcommand{\cN}{\mathcal N}
\newcommand{\cT}{\mathcal T}

\newcommand{\bI}{\mathbf I}
\newcommand{\bK}{\mathbf K}

\newcommand{\E}{\mathbb E}
\newcommand{\bX}{\mathbf X}

\newcommand{\bx}{\mathbf{x}}
\newcommand{\by}{\mathbf{y}}

\newcommand{\vb}{\bm{\beta}}

\newcommand{\va}{\bm{\alpha}}

\newcommand{\blu}[1]{\textcolor{black}{#1}}

\DeclareMathOperator\argmin{arg\, min}
\DeclareMathOperator\argmax{arg\, max}

\begin{document}

\maketitle

\begin{abstract}
We investigate Monte Carlo based algorithms for solving stochastic control problems with {local} probabilistic constraints. Our motivation comes from microgrid management, where the controller tries to optimally dispatch a diesel generator while maintaining low probability of blackouts {at each step}. The key question we investigate are empirical simulation procedures for learning the {state-dependent} admissible control set that is specified implicitly through a probability constraint on the system state. We propose a variety of relevant statistical tools including logistic regression, Gaussian process regression, quantile regression and support vector machines, which we then incorporate into an overall Regression Monte Carlo (RMC) framework for approximate dynamic programming. Our results indicate that using logistic or Gaussian process regression to estimate the admissibility probability outperforms the other options. Our algorithms offer an efficient and reliable extension of RMC to probability-constrained control. We illustrate our findings with two case studies for the microgrid problem.
\end{abstract}



\section{Introduction}
Stochastic control with probabilistic constraints is a natural relaxation of deterministic restrictions which tend to generate high costs forcing the avoidance of extreme events no matter their likelihood of occurrence. In contrast, with probabilistic constraints, constraint violation is tolerated up to a certain level offering a better trade-off between admissibility and cost. We refer to \cite{Geletu2013} for an overview of probability-constrained problems and list below some of our motivating settings and references:

\begin{enumerate}
\item Microgrid management: An electric power microgrid is a collection of intermittent renewable generator units, a conventional dispatchable diesel generator (or grid interconnection), and a battery energy storage system. The microgrid supplies electricity to a community in islanded mode, balancing fluctuating demand and supply. The operator achieves this by optimizing the use of the battery storage and the back-up dispatchable generator. Since perfect balancing is very expensive, it is common to allow for a small frequency of blackouts, i.e.~occurrences where demand outstrips supply. Mixed-integer linear programming approaches to this problem through approximating with more conservative convex constraints appear in \cite{Liu2017,adhithya18}.

\item Hydropower optimization: control of a hydropower dam with probabilistic constraints was discussed in \cite{delara17}. Within this setup, the controller observes random inflows from precipitation, as well as fluctuating electricity prices. \blu{His objective} is to control the downstream outflow from the dam to maximize profit from power sales,  while ensuring a minimum dam capacity with high probability. Other related works are  \cite{vanAckooij2014, Andrieu2010, Prekopa1978}.

\item Motion planning: finding the minimum-cost path for a robot from one location to another while avoiding colliding with objects that obstruct its path. Stochasticity in the environment implies that the robot motion is only partially controlled. Robust optimization that guarantees obstacle avoidance might be infeasible, making probabilistic constraints a viable alternative. Dynamic programming methods for unmanned aerial vehicles were introduced in \cite{Ludkovski2016} and the probabilistic-constrained motion of a robot was solved in \cite{Janson2018}.
\end{enumerate}

\textbf{Contribution.} {We consider state-dependent probabilistic constraints that are expressed through an expectation constraint at each system state. Our setup involves a continuous-state, continuous-time model with a discrete-step control, where the constraints are imposed step-by-step. Therefore, at each step and at each state, the controller must estimate which controls are admissible and then optimize over the latter. While this setting is simpler than global probabilistic constraints, it is much harder than unconstrained control, because a secondary numerical procedure is needed as part of dynamic programming to repeatedly compute the admissibility sets.} The canonical setup involves finite-horizon control of a stochastic process described through a stochastic differential equation of It\^{o} type. The solution paradigm involves the Bellman or Dynamic Programming equation (DPE), which works with discretized time-steps, but with a smooth spatial variable. In this context, we develop algorithms to solve stochastic optimal control problems with probabilistic constraints using regression Monte Carlo (RMC). To make this highly nontrivial extension to RMC, we  investigate tools from  machine learning (including support vector machines (SVM), Gaussian process (GP) regression, parametric density estimation, logistic regression and quantile regression) to {statistically estimate the admissible set as a function of the system state}. Our algorithm handles the two parts of the problem---the constraint estimation and the approximation of the conditional expectation---in parallel and with significantly lower simulation budget compared to a \textit{naive} implementation.

{After benchmarking the proposed approaches on two practical case-studies from energy battery management}, our main finding is to recommend logistic regression and GP-smoothed probability estimation as the best procedures. These methods are stable, relatively fast and allow for a variety of further adjustments and speed-ups. In contrast, despite theoretical appeal, quantile regression and SVM are not well-suited for this task. On a higher level, our main take-away is that DPE-based stochastic control with probabilistic constraints (SCPC) is well within reach of cutting-edge RMC methods. Thus, it is now computationally feasible to tackle such problems, opening the door for new SCPC models and applications.

\textbf{Related Models for Probabilistic Constraints.} There is an extensive literature on one-period optimization with chance constraints and on global multi-period probabilistic constraints.
For the one-period formulations, the most popular approach is to transform the problem into  linear or non-linear programs over a set of scenarios~\cite{Calafiore2006, Campi2009, Nemirovski2006, Alejandra2019}. In particular, Monte Carlo scenarios as employed below are very common, but the typical setup involves a single optimization problem, while we face an infinite family of them indexed by the system state $\bx$ and the time-step $n$. Global probabilistic constraints in multi-period settings are tackled from multiple perspectives. The dynamic programming method~\cite{Ono2015} incorporates the constraint into the objective function via  a Lagrange multiplier. The solution  is then obtained  by iteratively solving for the optimal control and the Lagrange multiplier. However, the solution is sub-optimal due to the duality gap. Mixed-integer linear programming~\cite{ShapiroReview} works by linearizing the constraints and requires to discretize the state space. An alternative~\cite{Blackmore2010, Blackmore2011} is to transform into a static problem which is however computationally feasible only under strict assumptions on the system dynamics and noise distribution. Yet another option is the stochastic viability approach~\cite{delara17,delara2010} that focuses on maximizing the probability of being admissible, which is defined both in terms of profit targets and satisfying constraints at every time step.

Compared to above, our model of a multi-period optimization with one-step probabilistic constraints applied at each time step is new in the literature. To our knowledge, the closest setup is studied recently in \cite{Tankov2017} to compute the hedging price of a portfolio whose risk is defined in terms of its future value  with respect to a set of stochastic benchmarks.  Besides a local probabilistic constraint, the authors also provide dynamic programming equations for multi-period constraints. However, their solution is driven by very specific loss functions and state processes. In contrast, we develop general purpose numerical schemes using statistical learning methods.

\section{Problem formulation}

We study numerical resolution of stochastic control problems on finite horizon $[0,T]$ with local implicit constraints, specifically we work with constraints defined through probabilistic conditions on the controlled state. Let $(\bX(t))_{t \geq 0}\in\cX \subset\mathbb{R}^d$ be a continuous time controlled Markov process adapted to a given filtration $(\cF_t)$. {The control is an $(\cF_t)$-adapted process $(u(t))_{t \geq 0}$, taking values  in  $\mathcal{W} \subset\mathbb{R}$. We further assume that control decisions are made at discrete epochs  $\{t_0,t_1,\ldots,t_N=T \}$; between time-steps the value of $u(t)$ remains constant. Thus, the control process is piecewise-wise constant and \blu{c\`adl\`ag} (right-continuous with left limits), and will be alternatively represented as $u(t) = \sum_{n=0}^{N-1} u_n \mathbbm{1}_{ [t_n, t_{n+1}) }(t)$.}

{Our \blu{setup} is essentially discrete-time as far as control is concerned, but we introduce the continuous-time system state because admissibility of actions depends on the trajectory of $\bX(t)$  between control points. In our motivating examples, the dynamics of the system is described by a stochastic differential equation:}
\[
\mathrm{d}\bX(t) = b(t,\bX(t),u(t))\mathrm{d} t + \sigma(t,\bX(t),u(t))\mathrm{d} B(t),
\]
where $(B(t))$ is a $m-$dimensional Brownian motion generating the filtration $(\cF_t)$ and $b: \mathbb{R}_+\times \cX \times \cW \rightarrow \mathbb{R}^d$ and $\sigma: \mathbb{R}_+\times \cX \times \cW \rightarrow \mathbb{R}^{d \times m}$ are measurable functions, such that a {unique} (weak) solution exists for admissible controls defined below and takes values in $\cX$. {For the convenience of notation, we will write $\bX_n$ for $\bX(t_n)$, $n=0, \ldots, N$, and will not indicate explicitly the dependence on the control $(u(t))$ if the latter is clear from the context.}

Admissibility is defined in feedback form via
\begin{equation}
\label{eq:generalConstraint}
\cU_{n:N} (\bX_n) =\Big\{ (u_{k})_{k=n}^N: P_k(\bX_k,u_k) \in \mathcal{A}_k\ \forall k \in \{n,\,\dots,\,N-1\} \Big\},
\end{equation}
for {deterministic functions} $P_k: \cX \times \cW \to \mathbb{R}$ and {given subsets} $\mathcal{A}_k \subset \mathbb{R}$. For the rest of the article we assume $P_n(\bX_n,u_n)$ and $\cA_n$ to be of the form
\begin{equation}
 P_n(\bX_{n},u_n ) \equiv p_n(\bX_{n},u_n ):= \mathbb{P}\Big( \blu{ \cG_n(  (\bX(s))_{s\in[t_n,t_{n+1})}  )} >0 \big|\, \bX_{n},u_n \Big) \ \text{ and } \ \mathcal{A}_n := [0,p),
 \label{eq:def1_FA}
\end{equation}
{where $\cG_n$ is a functional defined on trajectories of $(\bX(s))$ over the time interval ${s\in[t_n,t_{n+1})}$.}
In other words, we target  the set of controls such that the conditional probability of the ``failure" functional $\cG_n(\cdot)$ of $\bX$ being greater than zero is bounded by a threshold $p$, i.e.
\begin{equation}
\label{eq:specificConstraint}
\cU_{n}(\bX_{n}) :=  \Big\{ u \in \cW:\ p_n(\bX_{n},u_n )  <p  \Big\}.
\end{equation}
The threshold $p$ in equation~\eqref{eq:specificConstraint} is interpreted as relaxing the strong constraint $\cG_n(\cdot) \le 0$ which may not be appropriate in a stochastic environment. Typical values of $p$ would generally be small \blu{($p \approx 0.05$)}. {We assume in the following that at least one admissible control exists at any state and at any time, hence $\mathcal{U}_{n:N}$ is determined by $\mathcal{U}_k(\bx) = \mathcal{U}_{k:k}(\bx)$, $k=n, \ldots, N-1$, the non-empty sets of admissible controls satisfying the constraints at decision epochs $t_k$ conditional on $\bX_k = \bx$.}

The controller must optimize as well as evaluate the feasibility of proposed actions, with the performance criterion of the form
\begin{equation}
\label{eq:continuous_time_problem}
V_n(\bX_{n}) = \inf_{(u_{s})_{s=n}^N\in\mathcal{U}_{n:N}(\bX_n)} \Big\{\E \Big[ \sum_{k=n}^{N-1} \int_{t_k}^{t_{k+1}} \pi_s(\bX(s),u_k) \mathrm{d} s + W(\bX(t_N)) \Big| \,\bX_n \Big]  \Big\},
\end{equation}
where $W(\cdot)$ represents the terminal penalty and $\pi_t(\cdot, \cdot)$ the running cost.  We re-write~\eqref{eq:continuous_time_problem} in terms of the corresponding dynamic programming equation at step  $n$:
\begin{equation}
\label{eq:DPE}
\begin{split}
V_n(\bX_{n}) = &\inf_{u \in \cU_{n}(\bX_{n})} \Big\{ \mathcal{C}_{n}(\bX_{n},u) \Big\}, \\
 \text{where }\quad  \mathcal{C}_{n}(\bX_{n},u) & =\mathbb{E}  \left[ \int_{t_ n}^{t_{n+1}} \pi_s(\bX(s),u) \mathrm{d} s +V_{n+1}(\bX(t_{n+1})) \Big| \,\bX_n,u \right].
  \end{split}
\end{equation}
Above $\mathcal{C}_{n}(\bX_{n},u)$ is the continuation value, i.e.~reward-to-go plus expectation of future rewards, from using the control $u$ over $[t_n, t_{n+1})$. Moreover, given the state $\bX_n$, we say that $u^* \in \cU_{n}(\bX_{n})$ is an optimal control if $V_n(\bX_n) =  \mathcal{C}_{n}(\bX_{n},u^*) $. Since the admissible set $\cU_{n}(\bX_{n})$ is both time and state dependent, we need to estimate the continuation value $\mathcal{C}_{n}(\cdot, \cdot)$ and the admissible control set $\cU_{n}(\cdot)$ at every time~step. This is the major distinction from the standard scenario approach \cite{Nemirovski2006} in chance-constrained optimization where there is only a single problem to optimize over a fixed $\cU$, but no further indexing by $x$ and by $n$. The latter require a combination interpolation and optimization as part of the solution.

\textbf{Alternative Formulation of Admissibility.}
 We denote by $G_n(\bX_{n},u_n)$ as the regular conditional distribution  \cite{Karatzas1998} of the functional $\cG_n(\cdot)$  given $(\bX_{n},u_n)$:
\begin{equation}
G_n(\bX_{n},u_n) := \mathcal{L} \Big( \cG_n((\bX(s))_{s\in[t_n,t_{n+1})})\Big|\,\bX_{n},u_n \Big),
\label{eq:generic_G}
\end{equation}
{where $\mathcal{L}(\cdot | X_n, u_n)$ stands for a conditional law.}
When writing $\mathbb{P} \big(G_n(\bX_n, u_n) > z\big)$ or $\mathbb{E}\big[g\big(G_n(\bX_n, u_n)\big)\big]$ we mean the probability or the expectation with respect to this conditional distribution.

We may rewrite equation~\eqref{eq:def1_FA} through the corresponding $(1-p)^{th}$ quantile $q(\bX_n,u_n)$ of $G_n(\bX_{n},u_n)$:
\begin{equation}
 q_n(\bX_n,u_n) :\,(\bX_n,u_n) \mapsto \text{arg}\inf_z \Big\{ \mathbb{P}\Big( G_{n}(\bX_n,u_n)>z\Big)\le p\Big\}.
  \label{eq:def2_FA}
\end{equation}
Then using
\begin{equation} \label{eq:specificConstraint1}
\cU_{n}(\bX_{n}) :=  \left\{ u: p_n(\bX_{n},u )  <p  \right\} = \left\{ u: q_n(\bX_{n},u )  \le 0  \right\},
\end{equation}
we can set $P'_n := q_n$ and $\tilde{\mathcal{A}} = (-\infty,0]$ in \eqref{eq:generalConstraint}.
We will exploit this equivalence to propose quantile-based methods (Section~\ref{sec:admissibleSetEstimation}) for the admissible set.

\begin{remark}
\label{def:umin}
Assuming a one dimensional control $u_n \in \cW \subset \mathbb{R}$, and the probability $p_n(\bX_n,u_n )$ monotonically decreasing in $u_n$, estimating the admissible set $\cU_{n}(\bX_{n})$ is equivalent to estimating the \emph{minimum} admissible control
\[
u^{\min}_n(\bX_{n}) := \inf_{u \in \cW} \Big\{u: p_n(\bX_n,u ) < p \Big\}.
\]
The corresponding admissible set is $\cU_{n}(\bX_{n}) = \{u \in \cW :  u \geq u^{\min}_n(\bX_{n}) \}$.
\end{remark}

\begin{remark}
A more general version are implicit constraints of the form
$$
 \Big\{ u\in \mathcal{W}: \mathbb{E}\Big[ g\Big(G_n(\bX_{n},u)\Big) \Big] \le p \Big\},
$$
for a function $g: \mathbb{R} \to \mathbb{R}$. Even more abstractly, we can think of a generic implicit map $P_k(\cdot, \cdot)$ in \eqref{eq:generalConstraint} that defines $\cU_{n:N}(\cX_n)$, with the idea that inverting this map is numerically nontrivial (i.e.~$P_k$ is expensive to evaluate) and hence \emph{a priori} it is not clear which controls satisfy constraints and which do not.
\end{remark}

\begin{remark}
Equation~\eqref{eq:specificConstraint1} describes admissible controls $u$ for a given state $\bx$. The ``dual'' perspective is to consider the set of states $\cX^a_{n}(u) \subset \cX$ for which a given control $u$ is admissible:
\begin{equation}
\label{eq:specificConstraintThroughPartition}
\cX^a_{n}(u) :=  \Big\{ \bx \in \cX: p_n(\bx, u) <p  \Big\}.
\end{equation}
Often the cardinality of $\cX$ is infinite, while the control space $\mathcal{W}$ is finite, so that enumerating  \eqref{eq:specificConstraintThroughPartition} over $u \in \mathcal{W}$ is considerably easier than enumerating the uncountable family of sets $ \bx \mapsto \cU_n(\bx)$ in equation \eqref{eq:specificConstraint}. Furthermore, if $u \mapsto p_n(\bx, u)$ is decreasing for all $\bx \in \cX$, then we obtain an ordering $\cX^a_{n}(u_1) \subseteq \cX^a_{n}(u_2)$ for $u_1 \le u_2$.
The latter nesting feature corresponds to ranking the controls in terms of their ``riskiness'' with respect to $G_n$, so that the safest control will have a very large $\cX^a_{n}(u) $ (possibly all of $\cX$), while the riskiest control will have a very small admissibility domain.
\end{remark}

\subsection{Regression Monte Carlo}
\label{sec:rmcintro}
In this article we focus on simulation-based techniques to solve \eqref{eq:continuous_time_problem}. The overall framework is based on solving  equation~\eqref{eq:DPE} through backward induction on $n = N-1,N-2,\ldots$, replacing the true $V_n(\bx)$ with an estimate $\hat{V}_n(\bx)$. Since neither the conditional expectation, nor the admissibility constraint are generally available \blu{explicitly, those} terms must also be replaced with their estimated counterparts. As a result, we work with the approximate Dynamic Programming recursion
\begin{equation}
\label{eq:DPE-approx}
\begin{split}
\hat{V}_n(\bX_{n}) = &\inf_{u_n \in \hat{\cU}_{n}(\bX_{n})} \Big\{ \hat{\mathcal{C}}_{n}(\bX_{n},u_n) \Big\}, \\
 \text{where }\quad  \hat{\mathcal{C}}_{n}(\bX_{n},u_n) & :=\hat{\mathbb{E}} \left[\int_{t_n}^{t_{n+1}} \pi_s(\bX(s),u_n) \mathrm{d} s +  \hat{V}_{n+1}(\bX(t_{n+1})) \Big|\bX_n,u_n\right].
  \end{split}
\end{equation}
Above, $\hat{\mathbb{E}}$  is the approximate projection operator and the set of admissible controls $\hat{\cU}_n$ is also approximated via either
 $\hat{p}_n(\cdot,\cdot)$, i.e., $\hat{\cU}_{n}(\bX_{n}) :=  \big\{ u:\hat{p}_n(\bX_{n},u) < p \big\}$,  or $\hat{q}_n(\cdot,\cdot)$, i.e., $\hat{\cU}_{n}(\bX_{n}) =\big\{ u:\hat{q}_n(\bX_{n},u) \le 0 \big\}$, see \eqref{eq:specificConstraint1}. The estimated optimal control $\hat{u}_n \in \hat{\cU}_n(\bX_n) $ satisfies $\hat{V}_n(\bX_{n}) = \hat{\mathcal{C}}_{n}(\bX_{n},\hat{u}_n). $

The key idea underlying our algorithm and defining the Regression Monte Carlo paradigm is that $\hat{\mathbb{E}}$ and $\hat{\cU}$ are implemented through empirical regressions based on Monte Carlo simulations. In other words, we construct \emph{random}, probabilistically defined approximations based on realized paths of $\bX$. This philosophy allows to simultaneously handle the numerical integration (against the stochastic shocks in $\bX$) and the numerical interpolation (defining $\hat{V}_n(x)$ for arbitrary $x$) necessary to solve \eqref{eq:DPE-approx}.

To understand RMC, recall that specifying $\hat{\mathbb{E}}$ is equivalent to approximating the  conditional expectation map $(\bx,u) \mapsto \mathbb{E}[ \psi\big((\bX(s))_{s\in[t_n,t_{n+1}]}\big) | \bX_n=\bx, u_n = u] =: f(\bx,u)$ where we specifically substitute
\[
\psi\big((\bX(s))_{s\in[t_n,t_{n+1}]}\big) = \int_{t_n}^{t_{n+1}} \pi_s(\bX(s),u_n) \mathrm{d} s +  \hat{V}_{n+1}(\bX(t_{n+1}) ).
\]
To do so, we consider a dataset consisting of inputs $(\bx^{1}_n,u^1_n), \ldots, (\bx^{M_c}_n,u^{M_c}_n)$ and the corresponding pathwise realizations $y^1, \ldots, y^{M_c}$ with $y^j = \psi\big((\bx(s))^{j}_{s\in[t_n,t_{n+1}]}\big)$, where $(\bx(s))^{j}_{s\in[t_n,t_{n+1}]}$ is an independent draw from the distribution of the process $(\bX(s))_{s\in[t_n,t_{n+1}]}|(\bx^{j}_n,u^{j}_n)$. Then we use the training set $\{\bx^j_n, u^j_n, y^j\}_{j=1}^{M_c}$ to {compute} $\hat{f}$, an estimator of $f$, via regression.

Similarly, estimating $\cU_n$ is equivalent to learning the conditional probability map $p_n(\bx,u)$ (or the conditional quantile map ${q}_n(\bx,u)$ in \eqref{eq:specificConstraint}) and then comparing to the threshold value $p$ (zero, respectively). This statistical task, whose marriage with RMC is our central contribution, is discussed in Section \ref{sec:admissibleSetEstimation}.

The technique of using regressions for the approximation of the continuation value was developed in the celebrated works by \cite{ls2001} and  \cite{tvr} in the context of American option pricing and further enhanced in \cite{warin12,ludkovski15}. This was extended for storage problems and controlled state process in \cite{ludkovski10, boogert08, boogert12, balata17, aditya2018}. Among the approaches for approximating $f$ we mention \cite{warin} and \cite{aditya2018} who exploit the structure of the problem to reduce the dimensionality of the regressions, \cite{balata17,balata18} who harness the distribution of process to reduce the variance of $\hat{f}$ and \cite{warin,langrene15, aditya2018} who use non-parametric regression methods for $\hat{f}$. \blu{Regression based approach has also been discussed in the context of stochastic dual dynamic programming for solving high dimensional storage problems in \cite{warin17}.}

In contrast to the above well-developed literature, {very little exists about estimating the set of admissible controls $\hat{\cU}_{n}(\bX_{n})$, which requires approximating $p(\bx,u)$ (or $q(\bx,u)$) in equation~\eqref{eq:specificConstraint}. There are results about learning a single global admissibility set $\hat{\cU}$ in a one-period setup, but those approaches do not transfer to our context of state- and time-dependent admissibility constraints.}
A {naive} approach is to estimate $\hat{\cU}_{n}(\bX_{n})$  for every state realized during the backward induction through nested Monte Carlo. Namely for each pair $(\bx, u)$ encountered, we may estimate the probability of violating the constraint by simulating $M_b$ samples from the conditional distribution $G_{n}(\bx,u)$ as $\{g_n^b(\bx,u) \}_{b=1}^{M_b}$. We then set $u \in \hat{\cU}_{n}(\bx)$ if $\bar{p}_n(\bx,u) < p$, where
\begin{align}\label{eq:barp}
\bar{p}_n(\bx,u):= \sum_{b=1}^{M_b} \frac{\mathbbm{1}_{g_n^b(\bx,u)>0}}{M_b}
\end{align}
 is the empirical probability. Although simple to implement, this Nested Monte Carlo (NMC) method is computationally intractable for even the easiest problems. As an example, a typical RMC scheme employs $M_c \approx 100,000$ and assuming $M_b=1000$ for inner simulations, which is necessary for good estimates of small probabilities $p \le 0.1$, would require $10^8$ simulation budget at every time-step to implement NMC. Note furthermore that NMC returns only the local estimates $\bar{p}(\bx, u)$; no \emph{functional} estimate of $\cU_n(\bx)$ or $\cX^a_{n}(u) $ is provided for an arbitrary $\bx$ or $u$, respectively. As a result, any out-of-sample evaluation (i.e.~on a future sample path of $\bX$) requires further inner simulations, making this  implementation even more computationally prohibitive.

An important challenge in using $\hat{\cU}$ is \emph{verifying} admissibility. Since we are employing Monte Carlo samples to decide whether $u$ is admissible at $\bx$, this is a probabilistic statement and admissibility can never be guaranteed 100\%. We may use statistical theories to quantify the accuracy of estimators of $\cU$, for example, by applying Central Limit Theorem tools for the estimator $\bar{p}(\bx,u)$ of the true $p(\bx,u)$. In particular,  we develop tools based on confidence intervals in order to make statements (with asymptotic guarantee) such as ``$u \in \cU$ with 95\% confidence'' (equivalent to $p(\bx,u) < p$ with 95\% probability conditional on the data collected). Achieving reasonable confidence levels calls for ``conservative'' estimators of $\hat{\cU}$. As we show, not doing so can make learning $\cU$ highly unreliable, frequently causing decisions that are inadmissible with respect the imposed probability constraint. Thus, the related construction of $\hat{\cU}^{(\rho)}$ with specified confidence level $\rho$ is a running theme in Section~\ref{sec:admissibleSetEstimation}.

\subsection{Motivation: controlling blackout probability in a microgrid}
\label{sec:motivationMicrogrid}
\begin{figure*}[!ht]
        \centering
        \begin{subfigure}[b]{0.4\textwidth}
            \centering
            \includegraphics[width=\textwidth]{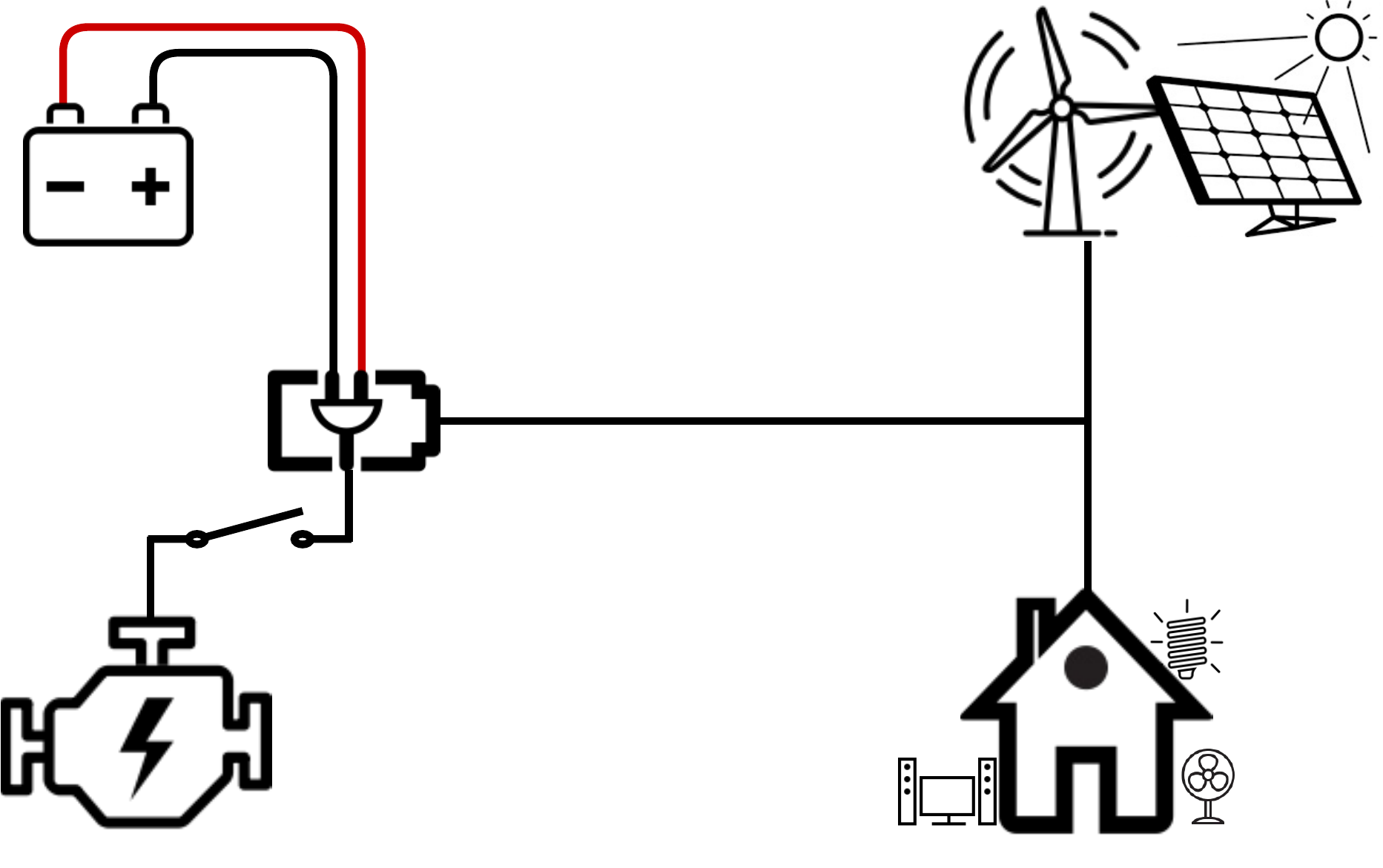}
            \label{fig:microgridTopology}
        \end{subfigure}
        \quad
        \begin{subfigure}[b]{0.4\textwidth}
            \centering
            \includegraphics[width=0.9\textwidth]{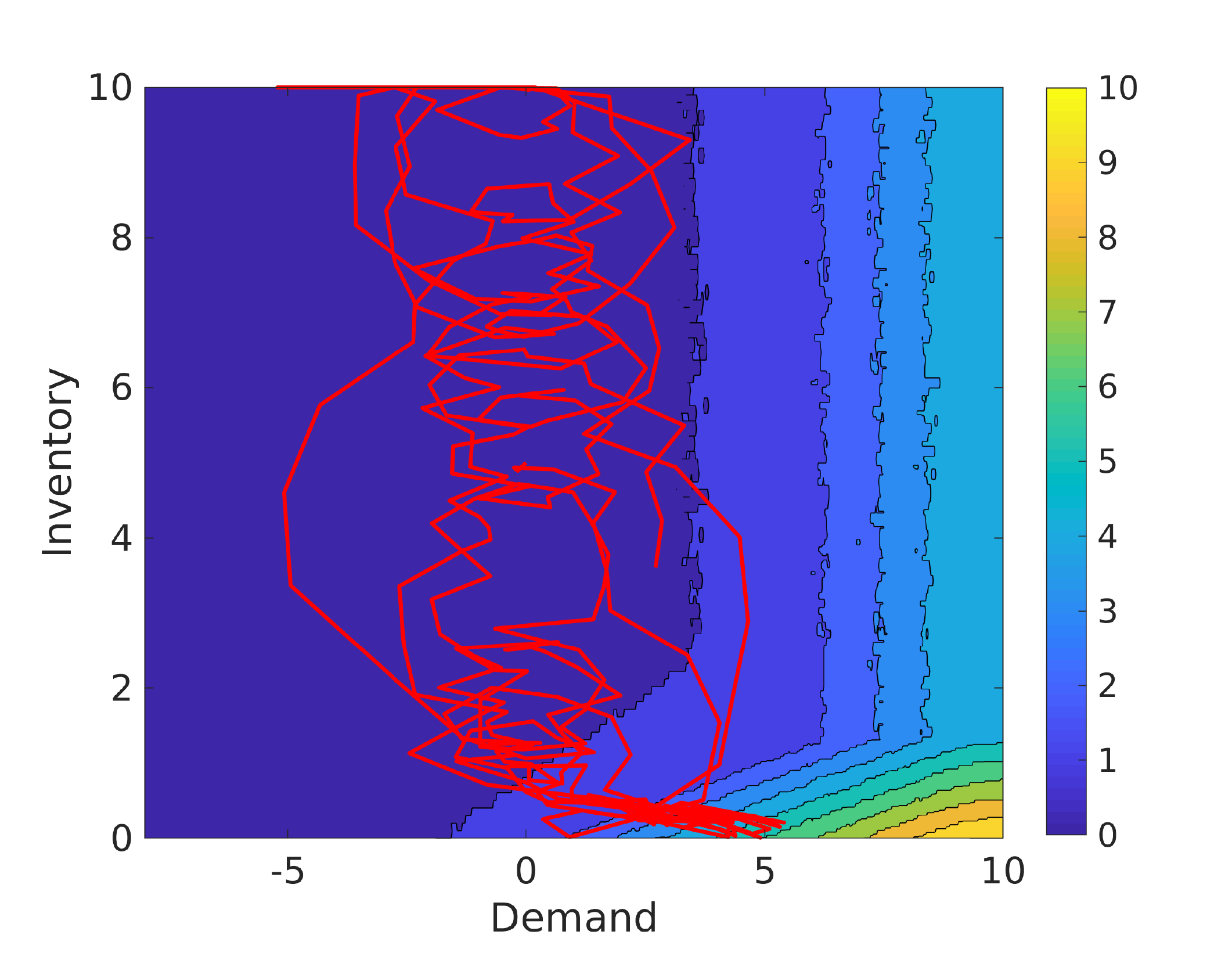}
            \label{fig:contour plot with pathwise strategy}
        \end{subfigure}
        \caption[ ]{\footnotesize Left panel: Microgrid topology: the load, the diesel generator, the battery and the renewables. Right: Contour plot for minimum admissible diesel output $(L,I,C) \mapsto u^{\min}_n(L,I,C)$ (see Remark~\ref{def:umin}).
         For $L < 0$, the constraint is not binding and $u^{\min}_n(L,I,C) = 0$. As demand increases, the constraint becomes more stringent, i.e.~$u^{\min}_n(L,I,C)$ increases in $L$. Red curve represents a path of the controlled demand-inventory pair $(L^{u^*}_n, I^{u^*}_n, C^{u^*}_n)$ following a myopic strategy choosing the minimum admissible control $u_n(L_n,I_n,C_n) = u^{\min}_n(L_n,I_n,C_n)$. \blu{The regime $C$ can be visualised by observing when the red line crosses on the R.H.S. of the first contour line, indicating the the diesel generator should be turned on.}}
        \label{fig:topology_pathwise}
    \end{figure*}

To make our presentation concrete, we illustrate the framework of \eqref{eq:DPE} by formalizing
the motivating application from microgrid management. A microgrid comprises renewable and traditional generation sources, along with a medium of storage, designed and managed to provide \blu{electrical power} to a community in a decentralized way. We consider a system composed of a dispatchable diesel generator, a renewable energy source and an electric battery storage.  The microgrid  topology is illustrated in the left panel of Figure~\ref{fig:topology_pathwise} and is same as the example discussed in \cite{aditya2018}.

In this context, the state variables are $\bX(t) = (L(t),I(t), C(t))$, where $L(t)$ is the net demand (demand net of renewable generation),  $I(t) \in [0, I_{\max}]$ is the state of charge of the battery, referred to as ``the inventory'', and $C(t) \in \{0,1 \}$  is the state of the diesel generator. $C(t)=0$ refers to diesel being OFF and $C(t)=1$ implies ON. The controller is in charge of the diesel through the control $u(t)$, which indicates the power output of the unit. We assume, for clarity of exposition, that the net demand $L(t)$ is an exogenous process, while $I(t)$ is controlled. We reiterate that the control decisions are made at discrete epochs $\{t_0,t_1,\ldots,t_{N-1}\}$, however these decisions affect the state of the system continuously. The choice of $u(t_n) \equiv u_{n}$ at time $t_n$ is based on minimizing the cost of running the microgrid, as well as controlling the probability of a blackout (i.e.~failing to match the net demand) during $[t_n,t_{n+1})$. The blackout is described through the imbalance process $S(s) := L(s) - u_{n} - B(s)$, $\forall s \in [t_n,t_{n+1})$, representing the difference between the demand and supply,  while the diesel output is held constant over the time step. The power output from the battery is a deterministic function of net demand, inventory and the control, $B(s) = \varphi(L(s),I(s),u_n) $ constrained by the physical limitations of the battery. $B(s) > 0$ implies supply of power from the battery and $B(s) < 0$ implies battery charging. The set of admissible controls is thus:
\begin{equation}
\cU_{n}(L_n,I_n,C_n) :=  \left\{ u: \mathbb{P} \Big(\sup_{s \in [t_n,t_{n+1})} S(s)>0 \Big|(L_n,I_n,C_n,u)\Big)  < p \right\}.
\label{eq:microgridConstraint}
\end{equation}
Thus in the context of microgrid, the conditional distribution $G_n$ of equation~\eqref{eq:generic_G} and the corresponding $p_n(L_n,I_n,C_n)$ are:
\begin{equation}
\begin{aligned}
G_{n}(L_n,I_n,C_n,u_n) &=  \mathcal{L} \Big(\sup_{s \in [t_n,t_{n+1})} S(s) \Big|(L_n,I_n,C_n,u_n) \Big),\\
p_n(L_n,I_n,C_n,u_n) &= \mathbb{P}(G_{n}(L_n,I_n,C_n,u_n) > 0).
\end{aligned}
\label{eq:microgrid_G}
\end{equation}
Because $p_n$ is not (in general) available analytically, the admissibility condition $p_n(L_n,I_n,C_n)<p$ is implicit.
Recall that we denote by $\cW = 0 \cup [\underline{u}, \bar{u}]$ the unconstrained control set. We assume that $u(t) = 0$ means that the diesel is OFF, while $u(t) > 0$ means that it is ON, and at output level $u(t)$. Thus, we define $C(s) = \mathbbm{1}_{ \{u_n > 0\} }$ $\forall s \in (t_n,t_{n+1}]$ with the time interval left-open in order to allow for identification of switching on and off of the diesel generator at times $t_n$. Notice also that the process $C(t)$ does not satisfy the controlled diffusive dynamics, but this slight extension of the framework does not impact on the methods and results presented.
We then look at the following formulation of the general problem:

\begin{equation}
\label{eq:microgrid_value_fun}
\begin{aligned}
&V_n(L_n,I_n,C_n)=\underset{\{u_{k}\}_{k=n}^{N-1}}{\min}\left\{\mathbb{E}\left[ \sum_{k=n}^{N-1} \Big(\mathbbm{1}_{\{{C}_{k}=0,u_{k}>0\}}\cK +\rho(u_{k})\Delta t_{k}  \Big)+W(L_N,I_N,C_N) \Big| (L_n,I_n,C_n) \right] \right\}, \\[5pt]
& \text{subject to } \qquad \mathbb{P}\Big(\underset{s \in [t_k,t_{k+1})}{\sup} S(s)>0\Big|(L_k,I_k,C_k, u_k) \Big) < p, \qquad \,k=n, \ldots, N-1,
 \end{aligned}
\end{equation}
where $\Delta t_{k} = t_{k+1}-t_k$, $\rho(u_k)$ is the instantaneous cost of running the diesel generator with power output $u_k$ and $\cK$ is the cost of switching it ON. We assume zero cost to turn the generator off.
The DPE corresponding to \eqref{eq:microgrid_value_fun} is the same as in \eqref{eq:DPE} with the integral running cost $\int_{t_n}^{t_{n+1}} \pi_s(\bX(s),u_n) \mathrm{d} s$ replaced by
\begin{equation*}
\mathbbm{1}_{\{{C}_{n}=0,u_{n}>0\}}\cK +\rho(u_{n})\Delta t_{n}.
\end{equation*}

\begin{remark}{The admissible set $\cU \subseteq \cW$ for this problem has the special structure: if $u \in \cU(\bx)$, then $\forall \ \cW \ni \tilde{u}>u, \tilde{u} \in  \cU(\bx) $. Hence, we may represent $\cU(\bx) = [u^{\min}_n(\bx), \bar{u}] \cap \cW$ in terms of the minimal admissible diesel output $u^{\min}_n(\bx)$.  Conversely, the admissibility domains for a fixed $u \in \cW$ are nested: if $u_1 \le u_2$ then $\cX^a_{n}(u_1) \subseteq \cX^a_{n}(u_2)$. This suggests to compute $\cX^a_{n}(u)$ sequentially as $u$ is increased and then invert to get $\cU(\bx)$.} \label{re:microgridAdmissibleSet} \end{remark}

To visualize the minimum admissible control $u^{\min}_n(\bx)$, the right panel of Figure~\ref{fig:topology_pathwise}  presents the map $\bx \rightarrow u^{\min}_n(\bx)$ under a constraint of  $p=0.01$ probability of blackout. We also present a path for $(L(t),I(t),C(t))_{t\geq 0}$ using a \emph{myopic} strategy where the controller employs the minimum admissible control at each point, $u_n := {u}^{\min}_n( L_{n}, I_{n},C_n) \;\forall n$. Notice how for the most part, $u^{\min}_n(\cdot)=0$ so that $\cU_n(\cdot) = \cW$ and the blackout constraint is not binding. This is not surprising, as blackouts are only possible when $L(t) \gg 0$ is strongly positive and the battery is close to empty, $I(t) \simeq 0$. Thus, except for the lower-right corner, any control is admissible.
As a result, only a small subset of the domain $\cX$ actually requires additional effort to estimate the admissible set $\cU(\bx)$. In our experience this structure, where the constraint is not necessarily binding and where we mostly perform unconstrained optimization, is quite common.

\section{Dynamic emulation algorithm}

In this section we present our Dynamic emulation algorithm which provides approximation for the admissible set $\hat{\cU}_n(\cdot)$ and the continuation value function $\hat{\mathcal{C}}_{n}(\cdot, \cdot)$. The crux of the algorithm are the following two steps, implemented in parallel at every time-step:
\begin{align}
\begin{split}
    \text{Generate design} & \rightarrow \text{Generate 1-step paths \& admissibility statistic} \rightarrow \text{Estimate admissible set}\\
    \text{Generate design} & \rightarrow \text{Generate 1-step paths \& pathwise profits} \rightarrow \text{Estimate continuation function}
\end{split}
\label{eqn:core_loop}
\end{align}
To estimate $\hat{\mathcal{C}}_{n}(\cdot, \cdot)$'s and $\hat{\cU}_n(\cdot)$'s, we proceed iteratively backward in time starting with known terminal condition $W(\bX)$ and sequentially estimate $\hat{\cU}_n$ and $\hat{\cC}_n$ for $n=N-1, \ldots, 0$. Assuming we have estimated $\hat{\cU}_{n+1},\ldots, \hat{\cU}_{N-1}$ and $\hat{\cC}_{n+1},\ldots,\hat{\cC}_{N-1}$, we first explain the estimation procedure for $\hat{\cU}_n$ and $\hat{\cC}_n$. This corresponds to a \texttt{fit} task. In the subsequent backward recursion at step $n-1$ we also need the \texttt{predict} task to actually \emph{evaluate} $\hat{V}_n(\bX_{n})$ which requires evaluating $\hat{\cC}_n(\cdot)$ at new (``out-of-sample'') inputs $\bX_{n},u_n$ which of course do not coincide with the training inputs $(\bx^{1}_n,u^1_n),\,\dots,\, (\bx^{M_c}_n,u^{M_c}_n)$.

\subsection{Estimating the set of admissible controls}
To estimate the set of admissible controls $\hat{\cU}_n(\cdot)$ at time-step $n$, we choose design $\CD^a_n:=(\bx_n^i,u_n^i, i=1,\ldots,M_a)$ and simulate trajectories of the state process $(\bX(s))^i_{s\in[t_n,t_{n+1})}$ starting from $\bX^i(t_n) = \bx_n^i$ and driven by control $u_n^i$. To evaluate the functional $\cG \big((\bX(s))^i_{s\in[t_n,t_{n+1})}\big)$, we discretize the time interval $[t_n,t_{n+1})$ into $K$ finer sub-steps with $\Delta n_k := t_{n_{(k+1)}} - t_{n_k}$ and define the discrete trajectory $\bx_{n}^{i}=\bx_{n_0}^{i}, \bx_{n_1}^{i}, \ldots, \bx_{n_{(K-1)}}^{i}, \bx_{n_K}^{i}$. We then record
\begin{equation}
\label{eq:wni}
w_n^i := \mathbbm{1}\Big(\cG((\bx_{n_k}^{i})_{k\in \{0,\ldots,K-1\} })>0\Big), \qquad i = 1,\ldots, M_a,
\end{equation}
where, formally, we extend $(\bx_{n_k}^{i})_{k\in \{0,\ldots,K-1\}}$ to a piecewise constant trajectory on $[t_n, t_{n+1})$.

Analogous to standard RMC, we now select an approximation space $\cH^a_n$ to estimate the probability $\hat{p}_n$ or the quantile $\hat{q}_n$, using the loss function $\mathcal{L}^a_n$ and apply empirical projection:
\begin{equation}
\label{eq:l2project_prob}
\hat{p}_{n} := \arg \min_{f_n^a \in \cH^a_n} \sum_{i=1}^{M_a}  \cL_n^a(f_n^a, w_n^i;\bx_{n}^{i},u_{n}^{i}).
\end{equation}
See Section~\ref{sec:admissibleSetEstimation} for concrete examples of $\cH^a$ and $\cL^a$.
Note that the approximations $\hat{p}_n$ and $\hat{q}_n$ must be trained on joint state-control datasets $\{\bx_n^{i},u_n^{i},w^i_n \}_{i=1}^{M_a}$ with $w^i_n$ dependent on the method of choice and moreover yield random estimators ($\hat{p}_n$ is a random variable).

Using the distribution of $\hat{p}_n(\bx,u)$ we may obtain a more conservative estimator that provides better guarantees on the ultimate admissibility of $(\bx, u)$.
As a motivation, recall the NMC estimator $\bar{p}_n(\bx,u)$ from \eqref{eq:barp}; for reasonably large $M_b \gg 20$, the distribution of $\bar{p}_n(\bx,u)$ is approximately Gaussian with mean $p_n(\bx,u)$ and variance \blu{ $ \frac{p_n(\bx,u)(1-p_n(\bx,u))}{M_b}$}. Defining 
\begin{align}
\hat{p}_n^{(\rho)}(\bx,u) &:= \bar{p}_n(\bx,u) +\xi^{(\rho)}_n(\bx,u) \\
& := \bar{p}_n(\bx,u) + z_{\rho} \sqrt{\frac{\bar{p}_n(\bx,u)(1-\bar{p}_n(\bx,u))}{M_b}},
\end{align}
where $z_{\rho}$ is the standard normal quantile at level $\rho$ and $\xi^{(\rho)}_n(\bx,u)$ represents a ``safe'' margin of error for $\bar{p}_n$ at confidence level $\rho$. The corresponding \blu{approximate} admissible set with  confidence $\rho$ is
\begin{equation}
\label{eq:admset_adaptive}
 \hat{\cU}_n^{(\rho)}(\bx) :=  \hat{\cU}_n^{\xi^{(\rho)}}(\bx)  = \left\{ u : \hat{p}_n(\bx,u) +\xi_n^{(\rho)}(\bx,u) < p \right\}.
 \end{equation}
More generally, we set the admissible set for a site $\bx \in \cX$ to
\begin{equation}
\label{eq:admset}
 \hat{\cU}_n^{\xi}(\bx) = \left\{ u : \hat{p}_n(\bx,u) +\xi_n(\bx,u) < p \right\},
 \end{equation}
where $\xi_n(\bx,u)$ ensures ``stronger'' guarantee for the admissibility of $u$ at $\bx$. The margin of estimation error can also be fixed, $\xi_n(\bx,u) = c \; \forall (\bx,u) \in \cX\times\cW $, which can be applied when the sampling distribution of $\hat{p}_n(\bx,u)$ is unknown. The corresponding admissible set
\begin{equation}
\label{eq:admset_nonadaptive}
\hat{\cU}_n^{\xi = c}(\bx)  = \left\{ u : \hat{p}_n(\bx,u) +c < p \right\}.
 \end{equation}
is equivalent to estimating $\hat{\cU}_n^{\xi = 0}(\bx)$ with a shifted lower probability threshold $p-c$. To simplify notation,  we use  $\hat{\cU}_n(\bx)$ to denote the unadjusted admissible set,
$\hat{\cU}_n(\bx) :=  \hat{\cU}_n^{\xi = 0}(\bx)
$
in the context of NMC. Analogously, we can adjust equations~\eqref{eq:wni}-\eqref{eq:l2project_prob} based on learning the quantile $q_n(\bx,u)$ to add a margin of error, $\hat{\cU}_n^{\xi}(\bx) = \{ u : \hat{q}_n(\bx,u) +\xi_n(\bx,u) \le 0 \}$.

\begin{remark}
In the context of one-step static optimization, conservative estimates for the admissible set were explored in \cite{Campi2009}. The idea is to consider $u$ to be admissible iff $\overline{p}(\bx, u) = 0$, i.e.~there are no failures observed in all of the simulated samples. The choice of the sample size $M_b$ is then determined using Chebyshev inequality that provides a theoretical guarantee on $\mathbb{P}( p(\bx,u) > 0\, |\, \overline{p}(\bx,u; M_b) = 0)$. This approach is similar in our setup to a violation margin $c=p$ in~\eqref{eq:admset_nonadaptive}. However, the guarantee on $p(\bx,u)$ is only applicable locally at the design sites. We are not aware of any tools that would allow non-parametric guarantee for \emph{any} $\bx \in \cX$. The regression based approach offers a model-based guarantee for admissibility of $u$ at any $\bx$ by setting the parameter $\rho$ in~\eqref{eq:admset_adaptive} or $c$ in~\eqref{eq:admset_nonadaptive}.
\end{remark}

\subsection{Estimating the continuation value}
\label{sec:continuation_value}

 To estimate the continuation value $\mathcal{C}_{n}(\cdot,\cdot)$, we choose a simulation design $\CD^c_{n} := (\bx_n^j,u_n^j, j = 1 \ldots, M_c)$ (which could be independent or equivalent to $\CD^a_n$) and generate one-step paths for the state process $(\bX(s))^j_{s\in[t_n,t_{n+1})}$ starting from $\bX^j(t_n) = \bx_n^j$ and driven by control $u_n^j$, \blu{ comprising again finer sub-steps} $\bx_{n}^{j}=\bx_{n_0}^{j}, \bx_{n_1}^{j}, \ldots, \bx_{n_{(K-1)}}^{j}, \bx_{n_K}^{j}$ (in principle the sub-steps could differ from the time discretization for $\hat\cU_n$). Next, we compute the pathwise cost $y_{n}^j$:
\begin{equation}
y_{n}^j  = \sum_{k=0}^{K-1} \pi_{n_k}(\bx_{n_k}^j,u_n^j) \Delta n_k +  v_{n+1}^j, \quad \text{ where } v_{n+1}^j = \inf_{u \in \hat{\cU}_{n+1}(\bx_{n_K}^{j})} \hat{\mathcal{C}}_{n+1}(\bx_{n_K}^{j},u), \ \ j = 1 \ldots M_c, \label{eq:one-step-aheadprofit}
\end{equation}
and we replace the time integral in~\eqref{eq:DPE} with a discrete sum over $t_{n_k}$'s. At the key step, we project $\{ y_{n}^j \}_{j=1}^{M_c}$ onto an approximation space $\cH^c_n$ to evaluate the continuation value $\mathcal{C}_{n}(\cdot,\cdot)$:
\begin{equation}
\label{eq:l2project}
\hat{\mathcal{C}}_{n}(\cdot,\cdot) := \arg \min_{f_n^c \in \mathcal{H}_n^c} \sum_{n=1}^{M_c} |f_{n}^c(\bx_{n}^{j},u_{n}^{j}) - y_n^j |^2.
\end{equation}
The design sites $\{\bx_n^{j},u_n^{j} \}_{j=1}^{M_c}$ could be same or different from those used for learning the admissible sets in the previous subsection.
Two standard approximation spaces $\cH^c_n$ used in this context are: global polynomial approximation and piecewise continuous approximation.

\begin{remark}
In the microgrid example of Section~\ref{sec:motivationMicrogrid} the running cost over $[n,n+1)$ is known once the control $u_n$ is chosen. Thus it can be taken outside the conditional expectation and the data to be regressed is simply $y^j =  v_{n+1}^j$.
\end{remark}

\paragraph{Global polynomial approximation:}
\label{sec:globalpoly}
This is a classical regression framework with polynomial bases $\phi_k(\cdot, \cdot)$  and $\hat{\mathcal{C}}_n^{\va}(\bx,u) := \sum_k \alpha_k \phi_k(\bx,u)$. The coefficients $\bm{\alpha}$ are fitted via
\begin{equation}
\hat{\bm{\alpha}} := \arg \min_{\bm{\alpha}} \sum_{j=1}^{M_c} \Big|\sum_k \alpha_k \phi_k(\bx^j,u^j) - y^j \Big|^2.
\label{eq:global_poly}
\end{equation}

As an illustration, for the microgrid example of Section \ref{sec:motivationMicrogrid} we construct a quadratic polynomial approximation when diesel generator is ON, $u>0$, using $10$ bases $\{1,L,I,u,L^2,I^2,u^2,LI,Iu,LI \}$ and a separate quadratic approximation with the $6$ basis functions $\{1,L,I,L^2,I^2,LI,LI \}$ when diesel generator is OFF, $u=0$. Polynomial approximation is easy to implement but typically requires many degrees of freedom (lots of $\phi$'s) to properly capture the shape of $\mathcal{C}$ and can be empirically unstable, especially if there are sharp changes in the underlying function (see for example \cite{aditya2018,langrene15}).

\paragraph{Piecewise continuous approximation:}
\label{sec:piecewise}
This is a state-of-art tool in low dimensions, $d \le 3$. The main idea is to employ polynomial regression in a single dimension and extend to the other dimensions via linear interpolation. As an example, for the microgrid with diesel generator ON, we have three dimensions $(L,I,u)$. We discretize inventory $I$ as $\{ I^0, I^1, \ldots, I^{M_I} \}$ and control $u$ as $\{ u^1, u^2, \ldots, u^{M_u} \}$ and fit independent cubic polynomials in $L$ for each pair $(I^l,u^e)$ with $l \in \{0,1,\ldots,M_I \}$ and \blu{ $e \in \{1,\ldots,M_u \}$}, i.e., $f_n^{l,e}(L) = \sum_{k}\alpha_{k}^{l,e}\phi_k(L)$. For any $(I,u) \in [I^l,I^{l+1}] \times [u^e,u^{e+1}]$ we then provide the interpolated approximation $\hat{\mathcal{C}}_n(L,I,u)$ as
\begin{gather}
 \hat{\mathcal{C}}_n(L,I,u) = \frac{\begin{bmatrix} I^{l+1} - I & I - I^{l}\end{bmatrix}
  \begin{bmatrix}
   f_n^{l,e}(L)  &  f_n^{l,e+1}(L) \\
    f_n^{l+1,e}(L)  &  f_n^{l+1,e+1}(L)
   \end{bmatrix}
  \begin{bmatrix} u^{e+1}-u \\ u-u^e \end{bmatrix}
  }{(u^{e+1}- u^e)(I^{l+1} - I^l)}.
\end{gather}

\paragraph{Nonparametric approximation:}
Further alternatives for $\cH_n^c$ can be found in \cite{aditya2018} who used Gaussian process regression and \cite{langrene15, langrene17} who used local polynomial regression. For semi-parametric approximation, \cite{warin12} developed piecewise multivariate linear regression.

\medskip
There are several possibilities for choosing the designs $\CD^a_{\cdot}$ and $\CD^c_{\cdot}$, see~\cite{aditya2018} for a detailed discussion of different regression designs $\CD^c_{\cdot}$ and their impact on the quality of the final solution.

\subsection{Evaluation}
We analyze the quality of the solution by computing three quantities on the out-of-sample dataset:
\begin{itemize}
\item estimate of the value function $V_0(\bx_0)$ at $t = 0$ and state $\bx_0$;
\item empirical frequency of inadmissible decisions on the controlled trajectories $\bx^{\hat{u}}_\cdot$;
\item statistical test for the realized number of constraint violations (blackouts for the microgrid).
\end{itemize}
Good solutions should minimize costs and not apply inadmissible controls. However, since we employ empirical estimators, $\cU$ is never known with certainty and we must handle the possibility that constraints are violated with probability more than $p$. In turn this leads to the trade-off between complying with \eqref{eq:generalConstraint} and optimizing costs. Similar treatment of constraints in the context of sample average approximation of probabilistic constrained optimization problems have been discussed in \cite{Nemirovski2007, Luedtke}. Moreover, our framework implies that the whole algorithm is stochastic: multiple runs will lead to different results since both $\hat{p}_n$ and $\hat{\mathcal{C}}_n$ are impacted by the random samples $y^j_n$ and $w_n^i$.

\textbf{Estimate of the value function:}
We evaluate the value function $\hat{V}_0(\bx_0)$ at time $t_0 = 0 $ and state $\bx_0$ using $M'$ out-of-sample paths $(\bx_{0:N}^{\hat{u}, m'}), m'=1,\ldots, M'$. Each trajectory $(\bx_{0:N}^{\hat{u},m'})$ is generated by applying the estimated optimal control $\hat{u}_{0:N-1}$ based on the continuation value  and admissible sets $(\hat{\mathcal{C}}_n, \hat{\cU}_n)_{n = 0}^{N-1}$ leading to the realized pathwise cost
$$
v_0(\bx_{0:N}^{\hat{u},m'}) := \blu{\sum_{n=0}^{N-1} \sum_{k=0}^{K-1} \pi_{n_k}(\bx_{n_k}^{\hat{u},m'},\hat{u}_n^{m'}) \Delta n_k + W(\bx_{N}^{\hat{u}, m'})}.
$$
The resulting empirical Monte Carlo estimate is
\begin{align}
\hat{V}_0(\bx_0) \simeq \frac{1}{M'}\sum_{m'=1}^{M'} v_0(\bx_{0:N}^{\hat{u},m'})
\end{align}
 and represents an unbiased estimation of the value of the control policy and an \blu{asymptotic} upper bound estimation of the value function, provided all controls used are admissible.

\textbf{Empirical frequency of inadmissible decisions  on the controlled trajectories}: For the $M'$ out-of-sample paths, we compare the estimated optimal control $\{\hat{u}_n(\bx_n^{\hat{u},m'})\}_{n=m=1}^{N-1,M'}$ against the minimum admissible control  $\{u^{\min}_n(\bx_n^{\hat{u},m'})\}_{n=m=1}^{N-1,M'}$ assumed for a second to be known. Namely, for each path we compute the number of inadmissible decisions $w_0( \bx_{0:N}^{\hat{u},m'})$ and the respective empirical frequency $w_{freq}$ as:
\begin{equation}
w_0( \bx_{0:N}^{\hat{u},m'}) := \sum_{n}\mathbbm{1}_{ \hat{u}_n(\bx_n^{\hat{u},m'}) < u^{\min}_n(\bx_n^{\hat{u},m'})} \quad \text{and} \quad w_{freq} := \frac{1}{N \cdot M'}\sum_{m'=1}^{M'} w_0(\bx_{0:N}^{\hat{u},m'}),\label{eq:w_freq}
\end{equation} respectively.   We employ these metrics in Section~\ref{sec:numerics}, where a ``gold standard'' $\{u^{\min}_n(\bx_n^{\hat{u},m'})\}_{n=m=1}^{N-1,M'}$ is obtained by brute force, utilizing a simulation budget $10^5$ larger than  for the actual methods we are comparing. Empirical gold standard is a common technique when analytical benchmark is unavailable, see e.g.~\cite{gpdp}. A good estimation method should yield $w_{freq} \simeq 0$. \blu{However, if $w_{freq} > 0$, a controller can choose margin of error $\xi_n$ (Equation~\eqref{eq:admset}) to reduce $w_{freq}$ at the expense of higher cost $v_0(\bx_{0:N}^{\hat{u},m'})$.}
\blu{\begin{remark}
In Equation~\eqref{eq:w_freq} we exploit the structure of the admissible set for the microgrid control (cf. Remark~\ref{re:microgridAdmissibleSet}). Generally, $w_0( \bx_{0:N}^{\hat{u},m'}) := \sum_{n}\mathbbm{1}_{ \hat{u}_n(\bx_n^{\hat{u},m'}) \notin \cU_n(\bx_n)} $, where $\cU_n(\bx_n)$ is either known in closed form or computed via empirical gold standard. The general setting, without any assumptions on the structure of $\cU_n(\bx_n)$, is computationally very expensive and beyond the scope of this work.
\end{remark}}

\textbf{Statistical test:} Next we propose statistical tests using the controlled trajectories to validate different methods for admissible set estimation. Such a test is essential to affirm the use of a numerical scheme for $\cU_n$ in the absence of a benchmark. As an example, in the context of microgrid we want to test the null hypothesis $H_0$ that the realized probability of blackouts is bounded to the required level against the alternative $H_1$ that their probability is too high.  Let
\begin{equation}
    B_n^{m'} = \mathbbm{1}\Big( \cG(\bx_{s \in [t_n,t_{n+1})}^{\hat{u},m'})>0 \Big), \ \ n = 0,\ldots, N-1 \text{ and } m'= 1,\ldots, M'.
\end{equation}

Ignoring the correlation due to the temporal dependence in $\bx_n$, we assume that $B_n^{m'} \sim \text{Bernoulli}(\tilde{p})$ \blu{are i.i.d. and $\tilde{p}$ represents the true (unknown) probability of blackout.} We want to test:
\begin{equation}
    H_0: \tilde{p} \leq p \quad \text{ vs. } \quad H_0: \tilde{p} > p.
\end{equation}
A common approach to such composite null hypothesis is to replace $H_0$ with a more conservative hypothesis $\tilde{p} = p$ leading to the test statistic
\begin{equation}
    \cT := \frac{\sum_{m', n} (B_n^{m'} - p) }{\sqrt{M'\cdot N \cdot p \cdot (1-p)}} \sim \cN(0,1).
    \label{eq:test1}
\end{equation}
Hence, $H_0$ is rejected at a confidence level $\alpha$ if $\cT > z_{\alpha}$ with $z_{\alpha} = \Phi^{-1}(\alpha)$, e.g.~$z_\alpha = 1.65$ for $\alpha = 95\%$.

\begin{remark}
The above test assumes independence and identical distribution of $B_n^{m'}$'s. In the context of the microgrid example, neither of the two assumptions are valid; $B_n^{m'}$ \blu{has a} different distribution because the state of the system affects the probability of a blackout, thus $\tilde{p}$ varies with $n, m'$. Furthermore, $B_n^{m'}$ are not independent as they are derived from a single, sequentially controlled trajectory. 
\end{remark}

\begin{remark}
In the microgrid setup, the blackout constraint is frequently not binding (the net demand is negative half of the time). Therefore, $\cT$ as defined in equation~\eqref{eq:test1} is most likely negative leading to accept the $H_0$ even when the method fails to choose the admissible control when the constraint is binding. We fix this by evaluating the sum only when the constraint is binding, i.e.
\begin{equation}
    \tilde{\cT} := \frac{\sum_{m', n} (B_n^{m'} - p)\mathbbm{1}_{  \blu{\hat{u}^{\min}_n(\bx_n^{\hat{u},m'}) } > 0 }  }{\sqrt{ p \cdot (1-p) \cdot  M'\cdot N \cdot w_{bind} }} \quad \text{where }\quad w_{bind} = \frac{ \sum_{m', n} \mathbbm{1}_{  \blu{ u^{\min}_n(\bx_n^{\hat{u},m'}) } > 0 }}{M'\cdot N},
    \label{eq:test2}
\end{equation}
\blu{ where $\hat{u}^{\min}_n(\bx_n^{\hat{u},m'})$ represents model estimate of the true minimum admissible control  $u^{\min}_n(\bx_n^{\hat{u},m'})$. }
\end{remark}

To wrap up this section, Algorithm~\ref{algo_Generalized} (dubbed Dynamic Emulation due to similarities with a related algorithm for unconstrained stochastic control from \cite{aditya2018}) summarizes the overall sequence of steps. Lines 1-2 contain the parameters to the algorithm, Lines 3-8 (and 14-19) yield the stochastic simulator which generates designs and corresponding one-step paths. Line 10 (and again Line 20) computes pathwise one-step costs. Line 12 is the admissible set estimation. Line 13 is the estimation of the continuation value.

Algorithm~\ref{algo_Generalized} carries several advantages. First and foremost it is very general, and does not make any restrictions on the distribution $G_n(\bX_{n},u)$ defining $\cU_n$ or the form of the payoffs $\pi(\bx,u)$. Hence it can be generically applied across a wide spectrum of SCPC problems. Second, the same template (in particular based on having two independent sub-modules) accommodates a slew of techniques for learning $\cC$ and $\cU$ bringing plug-and-play functionality, such as straightforward switching from probability to quantile estimation. Third, it allows for computational savings through parallelizing the estimation of $\cU$ and $\cC$, or by re-using the same design and simulations $\CD^a_n \equiv \CD^c_n$ for the computation of the two sub-modules.

\begin{remark}
The challenge of RMC methods is that the errors recursively propagate backward. As a result, poor estimation at one step can affect the overall quality of the solution. In our algorithm, the errors at every step occur due to:

\begin{itemize}
\item Approximation architecture $\cH^a_n$  for  $\hat{\cU}_n$ $\Rightarrow$ Projection error in admissible control set estimation;
\item Approximation architecture $\cH^c_n$ for $\hat{\cC}_n$ $\Rightarrow$ Projection error in estimating continuation value;
\item Designs $\CD^a_{n}$ and $\CD^c_{n} \Rightarrow$ Finite-sample Monte Carlo errors (difference between empirical estimates and theoretical projection-based ones)
\item Discretization of the time interval $[t_n,t_{n+1})$ using $\Delta n_k$ $\Rightarrow$ Integration error in approximating the integral $\int_{t_ n}^{t_{n+1}} \pi_s(\bX(s),u) \mathrm{d} s$ and the admissible set $\cU_n$.
\item Numerical approximation of the solution of the controlled dynamics of $\bX(t)$.
\item Optimization errors in maximizing for $\hat{u}$ over $\hat{\cU}$, especially when the control set $\mathcal{W}$ is continuous.
\end{itemize}
\end{remark}

\begin{algorithm}
\caption{Dynamic Emulation Algorithm}
\label{algo_Generalized}
\begin{algorithmic}[1]
\STATE {$N$ (time steps), $M_c$ (simulation budget for conditional expectation),}
\STATE{$M_a$ (simulation budget for admissible set estimation)}
\STATE{Generate designs:}
\STATE{\hspace{1cm}  $\CD^a_{N-1} := (\bx_{N-1}^{\CD^a_{N-1}},u_{N-1}^{\CD^a_{N-1}})$ of size $M_a$ for estimating $\hat{\cU}$. }
\STATE{\hspace{1cm}   $\CD^c_{N-1} := (\bx_{N-1}^{\CD^c_{N-1}},u_{N-1}^{\CD^c_{N-1}})$ of size $M_c$ for estimating $\hat{\mathcal{C}}$.}
\STATE{Generate one-step paths: }
\STATE{\hspace{1cm}  $\bx_{N-1}^{i,\CD^a_{N-1}} \mapsto \bx_{N}^{i,\CD^a_{N-1}}$ using $u_{N-1}^{\CD^a_{N-1}}$  for $i=1,\ldots,M_{a}$}
\STATE{\hspace{1cm}  $\bx_{N-1}^{j,\CD^c_{N-1}} \mapsto \bx_{N}^{j,\CD^c_{N-1}}$ using $u_{N-1}^{\CD^c_{N-1}}$  for $j=1,\ldots,M_{c}$}
\STATE{Terminal condition: }
\STATE{\hspace{1cm} $y_{N-1}^j \leftarrow \sum_{k=0}^{K-1} \pi_{(N-1)_k}(\bx_{(N-1)_k}^{j,\CD^c_{N-1}},u_{(N-1)_k}^{j,\CD^c_{N-1}}) \blu{\Delta n_k}  +  W(\bx_{N}^{j,\CD^c_{N-1}})$ for $j=1,\ldots, M_c$ }
\FOR{$ n=N - 1, \ldots, 1$}
\STATE{Estimate $\hat{\cU}_n(\cdot)$ using methods in Section \ref{sec:admissibleSetEstimation} and paths $\bx_{n}^{i,\CD^a_{n}} \mapsto \bx_{n+1}^{i,\CD^a_{n}}$ \\}
\STATE{$\hat{\cC}_n(\cdot,\cdot) \leftarrow \underset{{f_n \in \mathcal{H}_n^c}}{\arg \min}\sum_{j=1}^{M_c} |f_{n}(\bx_{n}^{j,\CD^c_{n}},u_{n}^{j,\CD^c_{n}}) - y_n^j |^2$  \\}
\STATE{Generate designs: }
\STATE{\hspace{1cm}  $\CD^a_{n-1} := (\bx_{n-1}^{\CD^a_{n-1}},u_{n-1}^{\CD^a_{n-1}})$ of size $M_a$ for estimating $\hat{\cU}$. }
\STATE{\hspace{1cm}   $\CD^c_{n-1} := (\bx_{n-1}^{\CD^c_{n-1}},u_{n-1}^{\CD^c_{n-1}})$ of size $M_c$ for estimating $\hat{\mathcal{C}}$. \\}
\STATE{Generate one-step paths: }
\STATE{\hspace{1cm}  $\bx_{n-1}^{i,\CD^a_{n-1}} \mapsto \bx_{n}^{i,\CD^a_{n-1}}$ using $u_{n-1}^{\CD^a_{n-1}}$  for $i=1,\ldots,M_{a}$ }
\STATE{\hspace{1cm}  $\bx_{n-1}^{j,\CD^c_{n-1}} \mapsto \bx_{n}^{j,\CD^c_{n-1}}$ using $u_{n-1}^{\CD^c_{n-1}}$  for $j=1,\ldots,M_{c}$ \\}
\STATE{$y^j_{n-1} \leftarrow  \sum_{k=0}^{K-1} \pi_{(n-1)_k}(\bx_{(n-1)_k}^{j,\CD^c_{n-1}},u_{(n-1)_k}^{j,\CD^c_{n-1}}) \blu{\Delta n_k}  + \underset{{u \in \hat{\cU}_n(\bx_n^{j,\CD^c_{n-1}})} }{\max} \Big\{ \hat{\cC}(n,\bx_{n}^{j,\CD^c_{n-1}},u) \Big\} $ $\forall j$}
\ENDFOR
\RETURN $\{ \hat{\mathcal{C}}_n(\cdot,\cdot), \hat{\mathcal{U}}_n(\cdot) \}_{n=1}^{N-1}$
\end{algorithmic}
\end{algorithm}

\section{Admissible set estimation}
\label{sec:admissibleSetEstimation}

In this section we propose two different approaches to estimate the admissible set of controls  $\cU_{n}$ in equation \eqref{eq:specificConstraint}:
\begin{itemize}
\item \textbf{Probability estimation:} Given a state $\bX_{n} = \bx$ and  $u \in \cW$, we estimate, via simulation, the probability of violating the constraint
\[
\hat{p}_n(\bx, u) \simeq \mathbb{P}\Big( G_{n}(\bx, u)>0\Big).
\]
It follows that $u \in \hat \cU_{n}(\bx)\;\Leftrightarrow\;\hat{p}_n(\bx, u) < p$.
Particularly, to compute $\hat{p}_n(\bx,u)$ we consider Gaussian process smoothing of empirical probabilities, logistic regression and parametric density fitting.

 \item \textbf{Quantile estimation:} We approximate the quantile $q_n(\bx,u)$ of $G_{n}(\bx,u)$ via  empirical ranking, support vector machines and quantile regression methods. 
 The admissible sets $\cU_{n}(\bx)$ and $\cX^a_{n}(u) $ are then defined as:
\[
\hat{\cU}_{n}(\bx) :=  \Big\{ u: \hat{q}_n(\bx,u) \le 0 \Big\}  \text{ and }\quad \hat{\cX}_n^a(u)  :=  \Big\{ \bx: \hat{q}_n(\bx,u) \le 0 \Big\}.
\]
\end{itemize}

To implement all of the above techniques we use Monte Carlo simulation, specifying first the simulation design and then sampling (independently across draws) the $G$'s or $Y$'s to be used as training data.  We  work in a flexible framework where samples of $G_{n}(\bx,u)$ are generated in batches of $M_b$ simulations from each design site $\{\bx^i,u^i\}_{i=1}^{M_a}$. The case of $M_b = 1$ corresponds to a classical regression approach, while large $M_b \gg 1$ can be interpreted as nested Monte Carlo averaging along $M_b$ inner samples.
\remark{In section~\ref{sec:continuation_value}, we parameterized the elements of the approximation space $\cH_n^c$ for estimation of the continuation value function $\hat{\cC}(\cdot,\cdot)$ via vectors $\va$ i.e.~$f_n^c(\bx,u) \equiv f_n^c(\bx,u; \va)$, cf.~\eqref{eq:l2project} and \eqref{eq:global_poly}. To distinguish, in the following sections we use $\vb$ to parameterize the approximators in $\cH_n^a$ for estimating the admissible set: $f_n^a(\bx, u) \equiv f_n^a(\bx,u; \vb)$} in eq.~\eqref{eq:l2project_prob}. The meaning and dimension of $\vb$ will vary from method to method.
\subsection{Probability estimation}
\label{sec:nestedMC}

\subsubsection{Interpolated nested Monte Carlo (INMC)}
\label{sec:inmc}

Recall the NMC method from Section~\ref{sec:rmcintro} where we select $M_a$ design sites of state-action pairs and simulate multiple paths from each site to \emph{locally} assess the probability of $G_n(\bx,u) >0$ (in what follows, we suppress in the notation the dependence on $n$). Specifically, for each design site $(\bx^i, u^i), \ i=1, \ldots, M_a$, we simulate $M_b$ batched samples from the distribution $G(\bx^i,u^i)$ as $\{ g^b(\bx^i,u^i) \}_{b=1}^{M_b}$. The unbiased point estimator of $p(\bx^i,u^i)$ is:
\begin{equation}
\label{eq:phat}
\bar{p}(\bx^i,u^i) := \sum_{b=1}^{M_b} \frac{\mathbbm{1}_{g^b(\bx^i,u^i)>0}}{M_b}.
\end{equation}
Since \eqref{eq:phat} only yields $M_b$ local estimates $\bar{p}(\bx^i, u^i)$, for Algorithm~\ref{algo_Generalized} we have to extend them to an arbitrary state-action $(\bx, u) \mapsto \hat{p}_{\text{INMC}}(\bx, u)$. This is achieved by interpolating $\bar{p}(\bx^i, u^i)$'s, e.g.~linearly. The admissible set with confidence level $\rho$  becomes:$$\hat{ \cU}^{(\rho)}_{\text{INMC}}(\bx) := \left \{ u :\hat{p}_{\text{INMC}}(\bx,u) \leq p -z_{\rho} \sqrt{\frac{\hat{p}_{\text{INMC}}(\bx,u)(1-\hat{p}_{\text{INMC}}(\bx,u))}{M_b}} \right \}.$$
However, especially for $M_b$ small,  interpolation performs poorly because the underlying point estimates $\bar{p}(\bx^i, u^i)$ are \emph{noisy}. Therefore, smoothing should be applied via a statistical \emph{regression} model. Regression  borrows information cross-sectionally to mitigate the estimation noise, reducing the variance of $\bar{p}$.

\subsubsection{Gaussian process regression (GPR)}
\label{sec:gpr}
GPR is a flexible non-parametric regression method that views the map $(\bx,u) \rightarrow p(\bx,u) $ as a realization of a Gaussian random field so that any finite collection of $\{p(\bx,u)\}, (\bx,u) \in \cX \times \cW$ is multivariate Gaussian. For any $n$ design sites $\{(\bx^i,u^i )\}_{i=1}^{n}$, GPR posits that
\begin{align*}
p(\bx^1,u^1),\ldots,p(\bx^n,u^n) \sim \cN(\vv{m}_n , \bK_n)
\end{align*}
with mean vector $\vv{m}_n := [ m(\bx^1,u^1; \vb),\ldots, m(\bx^{n},u^{n}; \vb)]$
and $n\times n$ covariance matrix $\bK_n$ comprised of $\kappa(\bx^i,u^i,\bx^{i'},u^{i'}; \vb), \text{ for } 1 \leq i,i' \leq n $. The vector $\vb$ represents all the hyperparameters for this model.

Given the training dataset $\{(x^i,u^i),\bar{p}^i \}_{i=1}^{M_a}$ (where $\bar{p}^i$ is a shorthand for $\bar{p}(\bx^i,u^i)$), GPR infers the posterior of $p(\cdot, \cdot)$ by assuming an observation model of the form $\bar{p}(\bx,u) = p(\bx,u) + \epsilon$ with a Gaussian noise term $\epsilon \sim \cN(0,\sigma^2_{\epsilon})$. Conditioning equations for multivariate normal vectors imply that the posterior predictive distribution $p(\bx,u)|\{(x^i,u^i),\bar{p}^i \}_{i=1}^{M_a}$ at any arbitrary site $(\bx,u)$ is also Gaussian with the posterior mean $\hat{p}_{\text{GPR}}(\bx,u)$ that is the proposed estimator of $p(\bx,u)$:
\begin{align}
\label{eq:gp_mean}
&\hat{p}_{\text{GPR}}(\bx,u) := m(\bx,u)+ K^T(\bK + \sigma^2 \bI)^{-1}(\vv{p}-\vv{m} ) = \mathbb{E}\Big[p(\bx,u)\big|\vv{x},\vv{u}, \vv{p}\Big] \\ \notag
&\text{where}\quad \vv{x} = [\bx^1,\ldots,\bx^{M_a} ]^T,  \ \  \vv{u} = [u^1,\ldots,u^{M_a} ]^T, \ \  \vv{p} = [\bar{p}^1,\ldots,\bar{p}^{M_a} ]^T , \\ \notag
&\phantom{\text{where}\quad} K^T = [\kappa(\bx,u, \bx^1, u^1; \vb), \ldots, \kappa(\bx,u,\bx^{M_a},u^{M_a};  \vb)],\\
&\phantom{\text{where}\quad} \vv{m}= [ m(\bx^1,u^1; \vb),\ldots,  m(\bx^{M_a},u^{M_a}; \vb)],
\end{align}
and $\bK$ is $M_a \times M_a$ covariance matrix described through the kernel function $\kappa(\cdot,\cdot;\beta)$.

The mean function is often assumed to be constant $m(\bx,u;  \vb) = \beta_0$ or described using a linear model $m(\bx,u;  \vb) = \sum_{k=1}^{K} \beta_k \phi(\bx^i,u^i)$ with $\phi(\cdot,\cdot)$ representing a polynomial basis. A popular choice for the kernel $\kappa(\cdot,\cdot,\cdot,\cdot)$ is squared exponential (see equation~\eqref{eqn:gp_kernel1}) with $\{ \{ \beta_{\mathrm{len},k} \}_{k=1}^d, \beta_{\mathrm{len},u} \}$ termed the lengthscales and $\sigma_p^2$ the process variance of $p(\cdot, \cdot)$:
\begin{equation}
 \kappa(\bx^i,u^i,\bx^{i'},u^{i'}) = \sigma_p^2\exp{\Big(- \sum_{k=1}^d \frac{(x^{i,k} - x^{i',k} )^2  }{\beta_{\mathrm{len},k}}-\frac{(u^i - u^{i'} )^2  }{\beta_{\mathrm{len},u}}\Big)}.\label{eqn:gp_kernel1}
\end{equation}
\blu{The hyperparameters $\vb := ( \{\beta_{k} \}_{k=1}^{K}, \{\beta_{\mathrm{len},k} \}_{k=1}^d, \beta_{\mathrm{len},u},  \sigma^2_{p}, \sigma^2_{\epsilon})$ represent different attributes of $\hat{p}_{\text{GPR}}(\cdot,\cdot)$. The parameters $\{ \beta_k\}_{k=1}^{K}$ determine the trend of $\hat{p}(\cdot,\cdot)$ over the input domain.  The lengthscales 
$\{ \beta_{\mathrm{len},k} \}_{k=1}^d$ and $\beta_{\mathrm{len},u}$ determine spatial smoothness or how quickly the function changes. Small lengthscales make $\hat{p}$ to be ``bumpy'', while large lengthscales make $\hat{p}$ close to linear. The process variance $\sigma^2_{p}$ determines the amplitude of fluctuations in $\hat{p}_{\text{GPR}}$  and $\sigma^2_{\epsilon}$ represents the sampling noise variance. These hyperparameters are estimated by maximizing the log-likelihood function based on the dataset $\{(x^i,u^i),\bar{p}^i \}_{i=1}^{M_a}$. Besides squared exponential kernel (Equation~\ref{eqn:gp_kernel1}), other popular kernels include Mat\'ern-3/2 and Mat\'ern-5/2~\cite{DiceKriging}.}

A conservative estimate $\hat{p}_{\text{GPR}}^{(\rho)}(\bx,u) $ at confidence level $\rho$ is obtained by explicitly incorporating the (estimated) standard error of $\bar{p}(\bx^i,u^i)$ into the GPR smoothing. Namely, we adjust the training dataset to  $\{(x^i,u^i),\bar{p}^i_{\rho} \}_{i=1}^{M_a}$, where $\bar{p}^i_{\rho} := \bar{p}(\bx^i,u^i) + z_{\rho} \sqrt{\frac{\bar{p}(\bx^i,u^i)(1-\bar{p}(\bx^i,u^i))}{M_b}}$. The resulting $\hat{p}_{\text{GPR}}^{(\rho)}(\bx,u) $ is the counterpart of \eqref{eq:gp_mean} using $\{(x^i,u^i),\bar{p}^i_{\rho} \}_{i=1}^{M_a}$.

In Figure~\ref{fig:gpfit} we present the dataset $\{L^i,I^i,0,\bar{p}^i \}_{i=1}^{M_a}$ (background colormap) for the microgrid case study. The thick red line indicates the contour $\{\hat{p}_{GPR} = 5\% \}$, dividing the state space $\cX$ for $u=0$ into admissible $\cX^a(0)$ (left of red line) and inadmissible region $(\cX^a(0))^c$ (right of red line).

\subsubsection{Logistic regression (LR)}
\label{sec:logisticReg}
In the previous section, we created local batches to estimate ${p}(\bx^i, u^i)$ pointwise and then regressed these estimates to build a global approximator. A classical alternative is to learn the probability of $G(x,u) > 0$ using a logistic regression model. This setup uses a single sample $g(\bx^i,u^i)$ from $G(\bx^i,u^i)$ from each design site $(\bx^i,u^i)$  and transforms it to a binary response $y^i = \mathbbm{1}_{g(\bx^i,u^i)>0}$.
 The probability  $\hat{p}(\bx,u)$ is then modeled as a generalized linear model with a logit link function
\begin{align}
\mathbb{P}\Big(Y = 1| \bx,u \Big)  & = \frac{1}{1 + e^{-\vb^T\phi(\bx,u)}} =: \hat{p}_{LR}(\bx, u; \vb).
\label{eq:logistic}
\end{align}
The basis functions $\phi(\bx,u)$ could be polynomials, e.g.~quadratic or cubic in coordinates of $(\bx,u)$. The regression coefficients ${\vb}$ are fitted using the dataset $\{\bx^i,u^i,\by^i\}_{i=1}^{M_s}$, as the solution to
\begin{align}
\argmax_{\vb}  \sum_{i=1}^{M_s} \left\{ y^i \log p_{LR}(\bx^i,u^i;\vb) +  (1-y^i) \log (1-p_{LR}(\bx^i, u^i; \vb)) \right\}.
\label{eq:logistic_likelihood}
\end{align}

We may again create a more conservative estimate $\hat{ \cU}^{(\rho)}_{LR}(\bx)$ of $\hat{ \cU}_{LR}(\bx)$ at confidence level $\rho$ by utilizing the standard error for $\hat{p}_{LR}$ using the Delta method \cite{deltaMethod}:
$$\hat{ \cU}^{(\rho)}_{LR}(\bx) := \left \{ u :\hat{p}_{LR}(\bx,u, \vb) \leq p -z_{\rho} \sqrt{\hat{p}_{LR}(\bx,u)(1-\hat{p}_{LR}(\bx,u))\phi^T\text{Var}(\vb)\phi } \right \}.$$

In Figure~\ref{fig:Logisticfit}, we present the original realizations $y^i \in \{0,1\}$ (in blue) for a design in the input subspace $(L,I,u=0)$ of the microgrid case study. The figure indicates the resulting logistic regression fit $\hat{p}_{LR}(L,I,0)$ at levels $1\%$, $5\%$ and $10\%$ (i.e.~contour lines of $\hat{p}_{LR}(\hat{\vb}) \in \{ 0.01, 0.05, 0.1\}$). The admissibility set for $u=0$, $\cX_n^a(0)$ is the region to the left of the thick red contour.

\begin{remark}
Similar to Section~\ref{sec:nestedMC}, we can simulate batched samples from each design site for the logistic regression, leading to ``binomial'' observation likelihood instead of~\eqref{eq:logistic_likelihood}.
\end{remark}

\begin{remark}
A non-parametric variant of equation~\eqref{eq:logistic} is kernel logistic regression, where the basis functions are $\phi_j(\bx,u) = \kappa(\bx,u,\bx^j,u^j)$ for a kernel function  $\kappa$ centered at $(\bx^j, u^j)$. \blu{One common choice is radial basis functions} (RBF) where $\kappa(\bx,u,\bx^j,u^j) = \exp{(- \gamma_1 \|\bx - \bx^j \|^2_2  - \gamma_2 \|u - u^j \|^2_2)}.$ RBF can be interpreted as the squared-exponential kernel for a logistic Gaussian Process model, with a fixed bandwidth parameter $\gamma_i$. In contrast, in GPR the bandwidths are estimated through MLE.
\end{remark}

\subsubsection{Parametric density fitting (PF)}
\label{sec:parametericDensity}

This approach fits the distribution $G(\bx,u)$, and then analytically infers the probability $\mathbb{P}\Big( G(\bx, u)>0\Big)$ from the corresponding cumulative distribution function. This is done by proposing a parametric family $\{ f(\cdot; \Theta) \}$ of densities, fitting the underlying parameters $\Theta$ based on an empirical sample from $G$ and then evaluating the resulting analytical probability $\bar{p}_{PF}(\bx,u) := \int_{0}^{\infty} f_{G(\bx,u)}(z|\hat{\Theta}(\bx,u)) dz $.
PF yields a ``universal'' solution across a range of constraint levels~$p$.

At a design site $(\bx,u)$, the probability $p(\bx,u)$ is estimated in a two-step procedure: first estimated locally over a design $\CD^a = \{\bx^i, u^i\}$ and then regressed over the full input domain $\cX\times \cW$. For the first step, we apply nested Monte Carlo to generate a collection of realized $\{ g^b(\bx^i,u^i) \}_{b=1}^{M_b}$ that is used to construct a parametric density via the
maximum likelihood estimate:
\begin{equation}
\hat{\Theta}^i := \underset{\Theta}{\argmax} \sum_{b=1}^{M_b}\log f_{G}(g^b(\bx_i, u_i)|\Theta).
\label{eq:pf_theta}
\end{equation}
In the second step, we evaluate $\tilde{p}_{PF}(\bx^i,u^i) := \int_{0}^{\infty} f_{G}(z|\hat{\Theta}(\bx^i,u^i))$ and extend it to the full domain $\cX \times \cW$ based on the computed $\{\bx^i,u^i,\tilde{p}_{PF}(\bx^i,u^i)\}_{i=1}^{M_a}$ using $\cL_2$ projection:
\begin{equation}
\hat{p}_{PF} =  \underset{\hat p \in \mathcal{M}_T}{\argmin} \ \| \hat p(\bx^i,u^i) - \tilde{p}_{PF}(\bx^i, u^i) \|^2,
\label{eq:pf_phat}
\end{equation}
where $\mathcal{M}_T$ is an approximation space chosen for regression.
The admissible set $\cU(\bx)$ is estimated as:
$$ \hat{\cU}_{PF}(\bx) := \left \{ u :\hat{p}_{PF}(\bx,u)  \leq p \right \}.$$

\begin{figure*}[!ht]
        \centering
        \begin{subfigure}[b]{0.32\textwidth}
            \centering
            \includegraphics[width=\textwidth]{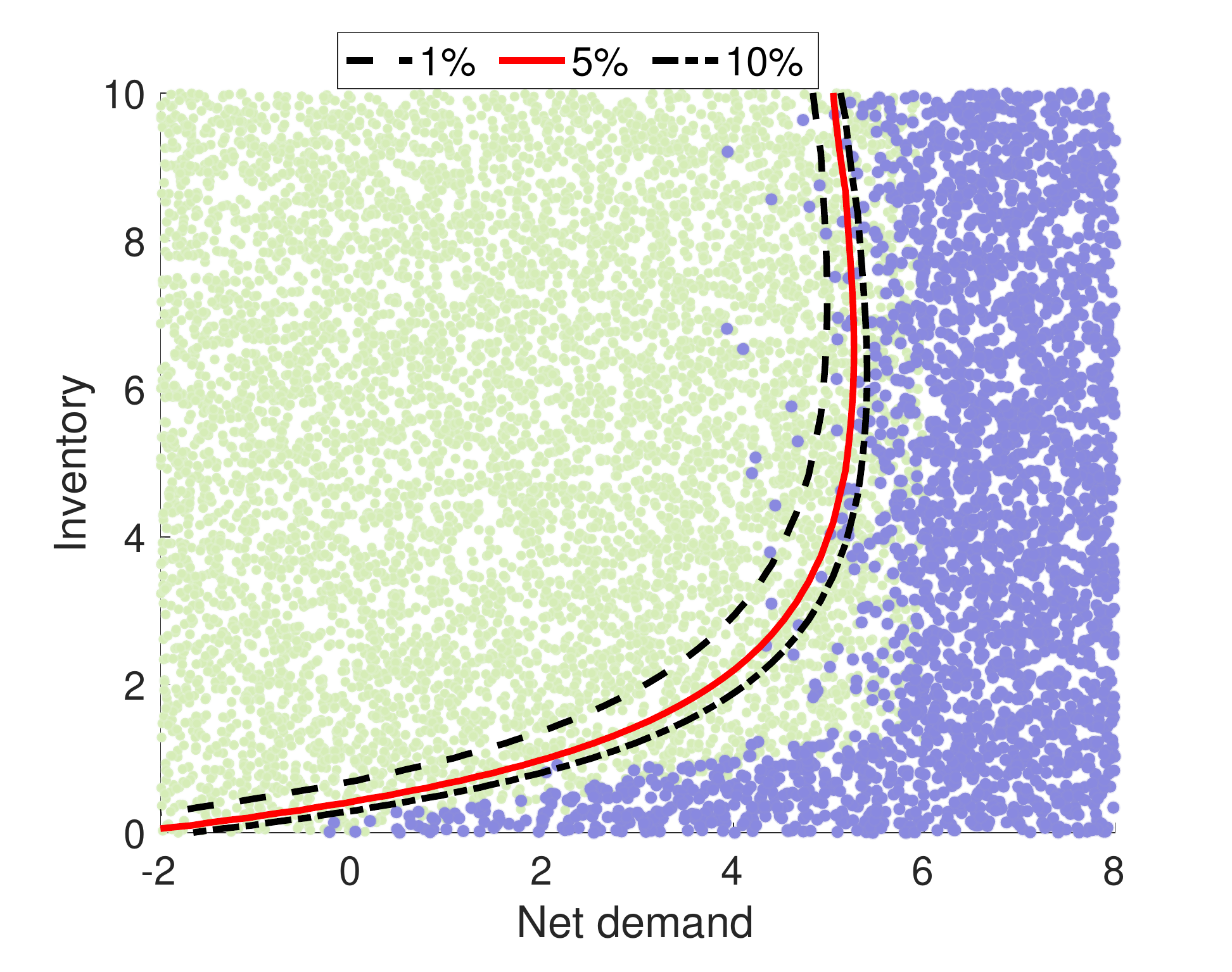}
            \caption[]%
            {{ LR $(M_a = 10000, M_b=1)$}}
            \label{fig:Logisticfit}
        \end{subfigure}
        \begin{subfigure}[b]{0.32\textwidth}
            \centering
            \includegraphics[width=\textwidth]{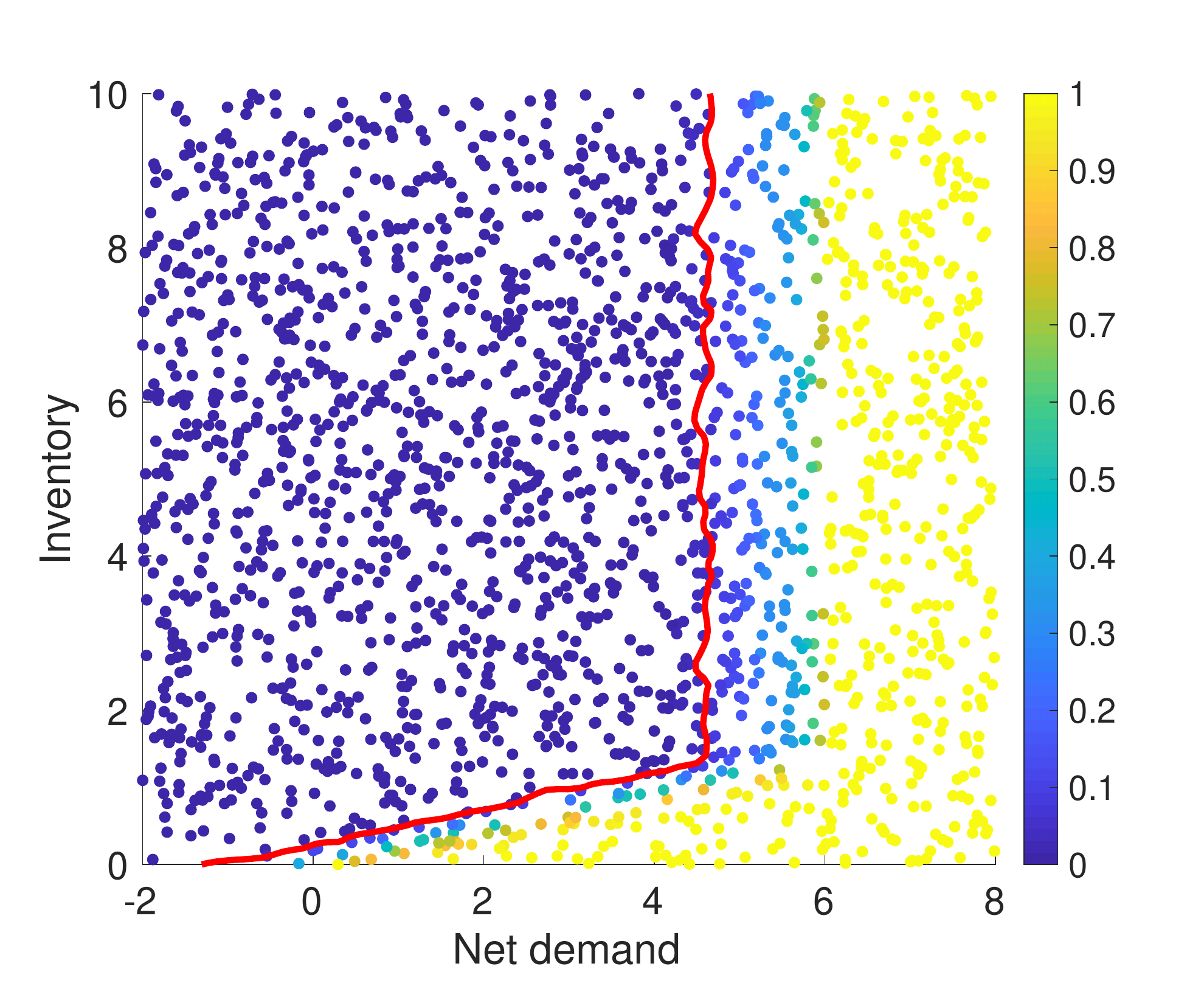}
            \caption[]%
            {{ GPR $(M_a=2000, M_b=50)$}}
            \label{fig:gpfit}
        \end{subfigure}
        \begin{subfigure}[b]{0.32\textwidth}
            \centering
            \includegraphics[width=\textwidth]{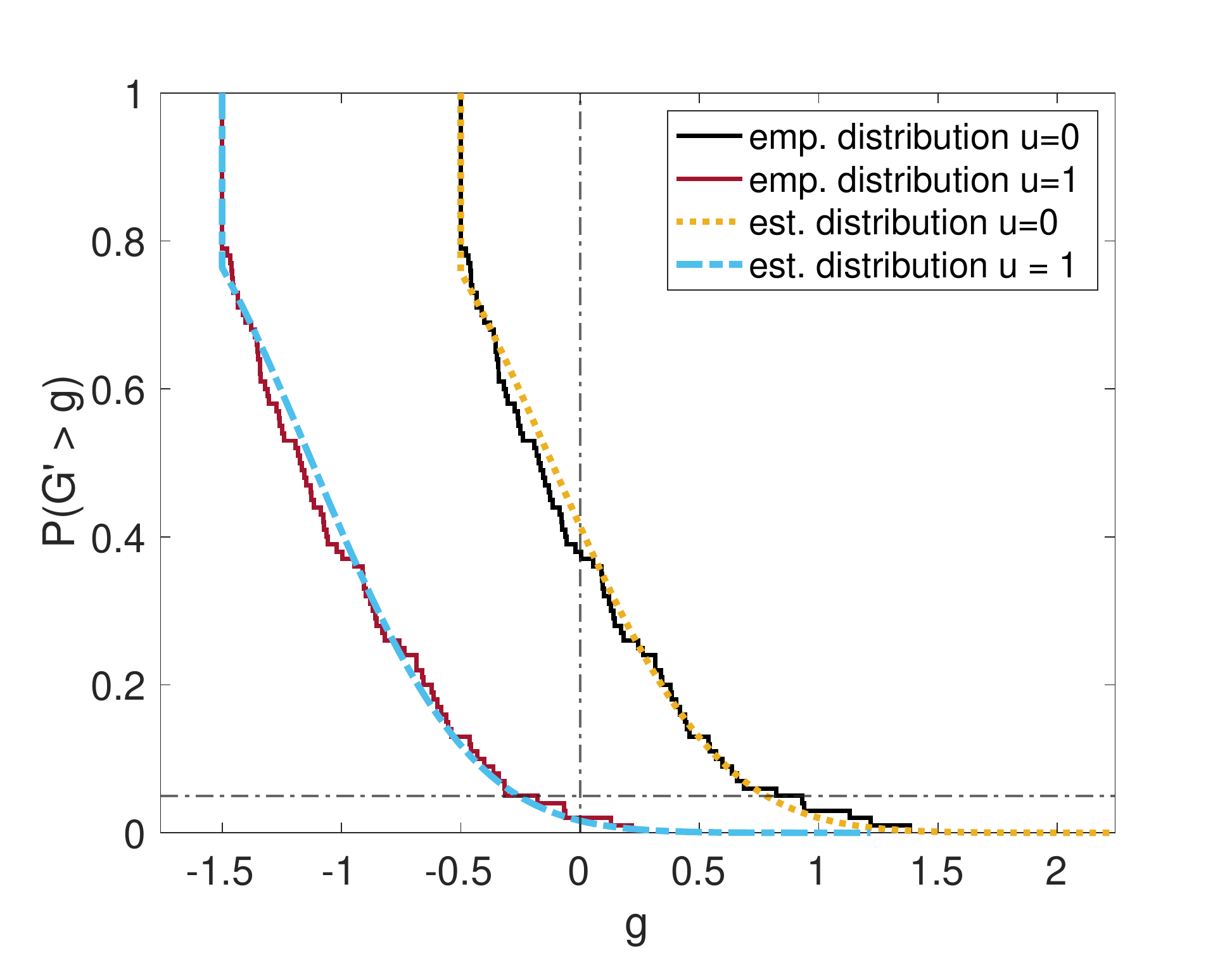}
            \caption[]%
            {{ PF $(M_b=100)$}}
            \label{fig:densityvisual}
        \end{subfigure}
        \begin{subfigure}[b]{0.32\textwidth}
            \centering
            \includegraphics[width=\textwidth]{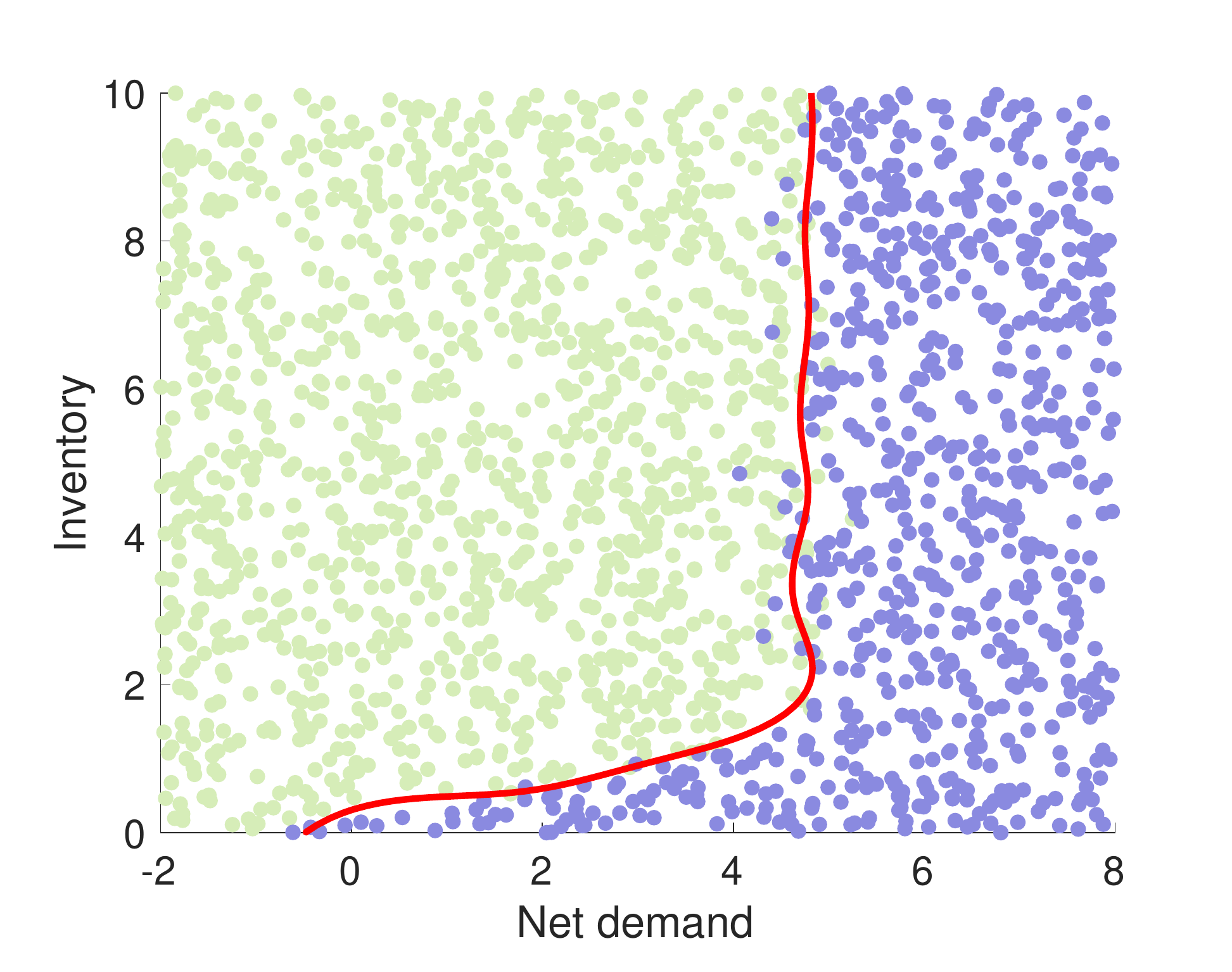}
            \caption[]%
            {{ SVM  $(M_a = 2000, M_b=50)$}}
            \label{fig:svmfit}
        \end{subfigure}
          \begin{subfigure}[b]{0.32\textwidth}
            \centering
            \includegraphics[width=\textwidth]{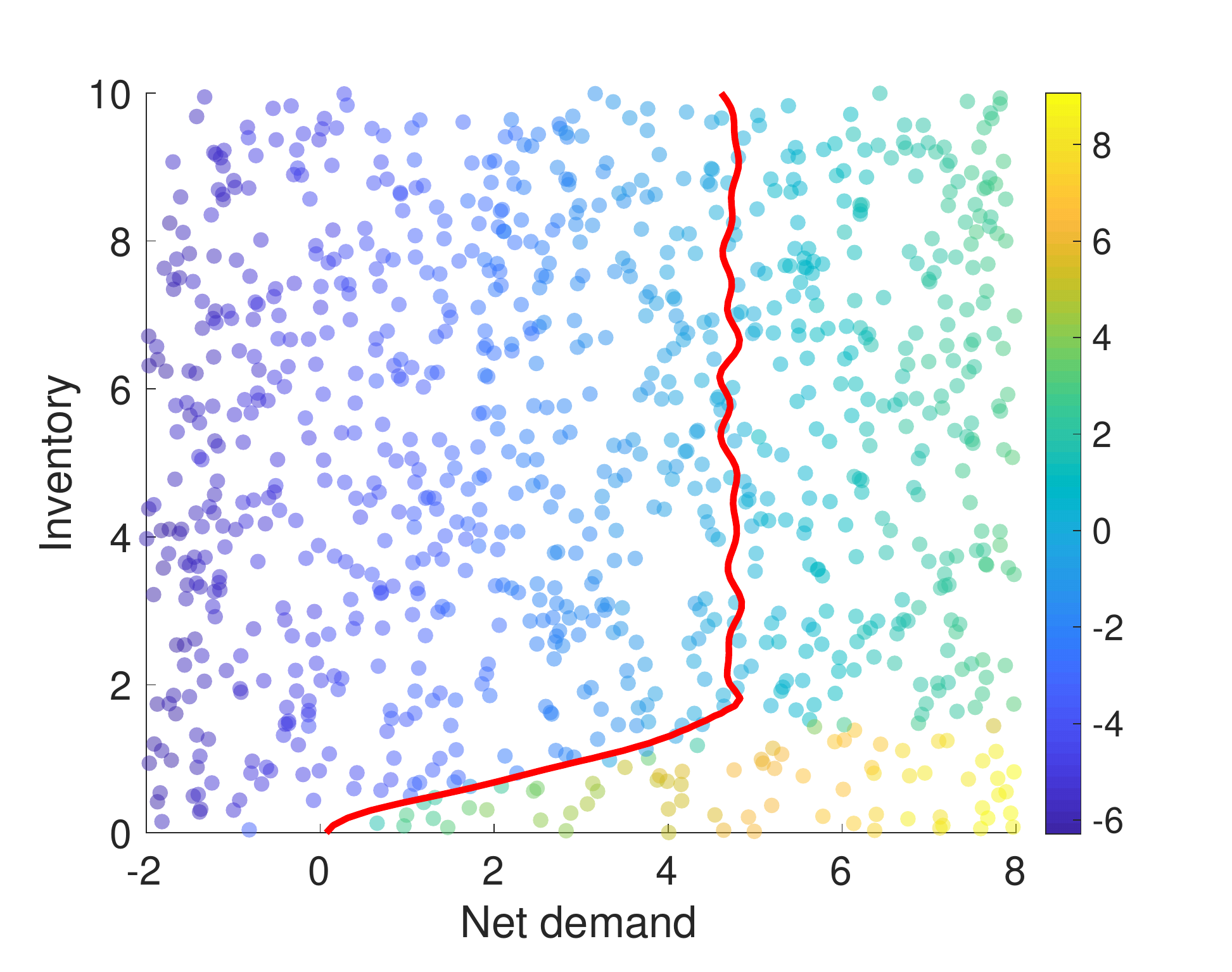}
            \caption[]%
            {{ EP $(M_a = 1000, M_b=100)$}}
            \label{fig:varfit}
        \end{subfigure}
        \begin{subfigure}[b]{0.32\textwidth}
            \centering
            \includegraphics[width=\textwidth]{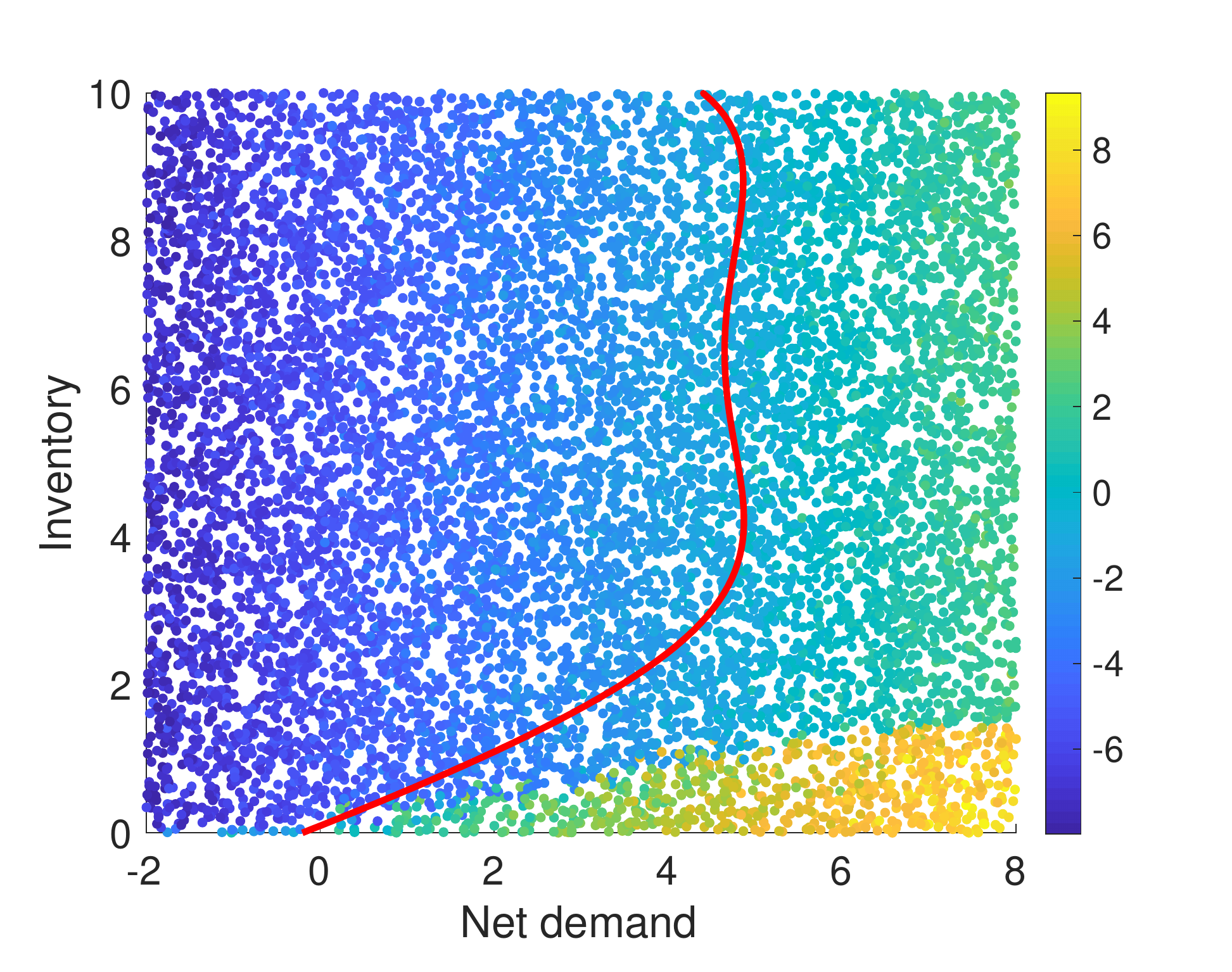}
            \caption[]%
            {{ QR $(M_a=10000, M_b=1)$}}
            \label{fig:qrfit}
        \end{subfigure} \vspace*{-10pt}
        \caption[ ]
         {\footnotesize{ Training data and fitted models for the methods of Section~\ref{sec:admissibleSetEstimation} at $u=0$. Top row: probability estimation schemes, bottom row: quantile estimation schemes. \textit{Top/left panel:} Training set $\{L^i,I^i,y^i \}_{i=1}^{M_a}$ for the LR model, color-coded according to the value of $y^i \in \{0,1\}$, along with the estimated contours for $\hat{p}_{LR}(L,I)$ at levels $\{1\%, 5\%, 10\% \}$. \textit{Top/center:} Training set $\{L^i, I^i, \bar{p}^i \}_{i=1}^{M_a}$ color-coded according to $\bar{p}^i$ for GPR along with the contour $\{\hat{p}_{GPR}(L,I) = 5\%\}$.  \textit{Top/right:} parametric density fitting at $L_0 = 5.5, I_0 = 1.48$ and $u \in \{0,1\}$. We show the empirical and fitted inverse cdf $\mathbb{P}(G' > g)$ based on a truncated Gaussian distribution.  \textit{Bottom/left:} Training set $\{L^i,I^i,y^i \}_{i=1}^{M_a}$ for SVM (color-coded according to $y^i \in \{-1,1\}$) and the  decision boundary in red. \textit{Bottom/center:} Training set $\{L^i,I^i,\bar{q}^i \}_{i=1}^{M_a}$ color-coded according to $\bar{q}^i$ for EP and the contour $\{\hat{q} = 0\}$. \textit{Bottom/right:} Training set $\{L^i,I^i,g^i \}_{i=1}^{M_a}$ color-coded according to $g^i$ for QR along with the contour $\{\hat{q}_{QR}(L,I) = 0\}$. All models share the same ground truth, so the red contours are identical up to model-specific estimation errors. }}
        \label{fig:fitExamples}
    \end{figure*}

A transformation of the distribution $G(\bx,u)$ might be important for above distribution fitting. For example, in the context of microgrid, in Section \ref{sec:motivationMicrogrid}, $G = \mathcal{L} \big(\sup_{s \in [t_n,t_{n+1})} S(s)\big)$ has a point mass at $0$ and thus, any  continuous distribution will lead to poor statistical estimation.
Using a transformation that preserves the probability of the target event,
\begin{equation}
\mathbb{P}\Big(\sup_{s \in [t_n,t_{n+1})} S(s)>0|\mathcal{F}_{n}\Big) = \mathbb{P} \left(\sup_{s \in [t_n,t_{n+1})} [L(s) - u_n - \frac{I(s)}{\delta s} \wedge B_{\max}]>0 \, \big|\mathcal{F}_{n} \right),
\end{equation}
we work with $G'(L_n, I_n,u_n) := \mathcal{L} \big(\sup_{s \in [t_n,t_{n+1})} [L(s) - u_n - \frac{I(s)}{\delta s} \wedge B_{\max}]\big)$.
In Figure~\ref{fig:densityvisual} we present the empirical and estimated probability $z \mapsto \mathbb{P}(G'(L_n, I_n, u_n)>z)$ when $L_n = 5.5, I_n = 1.48$  and $u_n \in \{0, 1\}$ for the microgrid example.  We model the distribution $G'$ using a truncated normal distribution, $\mathbb{P}(G' \leq g) = \Phi(\frac{g - \theta_2}{\theta_3})\mathbbm{1}_{g \geq \theta_1}$, with parameters $\Theta = (\theta_1, \theta_2,\theta_3)$  representing the location of censoring, the mean and the standard deviation respectively. At $L_n = 5.5, I_n = 1.48, u_n = 1.0$ and inner simulation budget $M_b=100$, the estimated parameters $(\hat{\theta}_1, \hat{\theta}_2, \hat{\theta}_3) = (-1.5,  -1.12, 0.53) $ result in probability $\tilde{p}_{PF}(5.5,1.48,1.0) = 0.016 $. The corresponding probability after $\cL_2$ projection (equation~\eqref{eq:pf_phat}) is $\hat{p}_{PF}(5.5,1.48,1.0) = 0.017$. Thus at $p=0.05$, the control $u =1.0 \in \hat{\cU}_n$ is admissible. However, at $u_n =0$, $(\hat{\theta}_1, \hat{\theta}_2, \hat{\theta}_3) = ( -0.5,	-0.12, 0.55)$, $\tilde{p}_{PF}(5.5,1.48,0.0) = 0.414$ and $\hat{p}_{PF}(5.5,1.48,0.0) = 0.429 $, thus the control $u =0 \notin \hat{\cU}_n$ is inadmissible.

\subsection{Quantile estimation}

In this section we consider methods for modeling and estimating $q(\bx^i,u^i)$, the $(1-p)$-th quantile of the distribution $G(\bx^i,u^i)$. Admissibility corresponds to the quantile being negative.

\subsubsection{Empirical percentiles (EP)}
\label{sec:var}

As before, we start by choosing $M_a$ design sites of state-action pairs and generate batched samples $\{ g^b(\bx^i,u^i) \}_{b=1}^{M_b}$ from each design site $(\bx^i,u^i)$.
The empirical estimate of $q(\bx^i,u^i)$ is simply the
 $(1-p)^{th}$ percentile of the realized $\{ g^b \}_{b=1}^{M_b}$ (which requires $M_b > p^{-1}$):
$$
\bar{q}(\bx^i,u^i) = \text{percentile}\left(\{ g^b \}_{b=1}^{M_b}, 100(1-p)\% \right).
$$
Similar to previous methods, we extend to arbitrary $(\bx,u) \mapsto \hat{ q}(\bx, u)$ using regression on the dataset $\{\bx^i,u^i, \bar{q}(\bx^i,u^i)\}_{i=1}^{M_a}$ and an approximation space $\mathcal{M}_q$. The set of admissible controls for $\bx$ is:
$\hat{\cU}_{EP}(\bx) :=  \Big\{ u: \hat{q}(\bx,u) \le 0 \Big\}.$ In Figure~\ref{fig:varfit} we show the estimated  $\hat{q}(\cdot,\cdot,0) $ indicated via the background colormap. The thick red line indicates the zero-contour $\hat{q} = 0$, so that the admissibility set for $u=0$, $\cX_n^a(0)$, is the region to the left of the contour.
\begin{remark}
This approach is similar to the INMC approach discussed in Section~\ref{sec:nestedMC}, however, here we model the quantile rather than the probability of exceeding zero.
Furthermore, we can use the regression standard error of $\hat{q}(\cdot,\cdot)$ to construct  a more conservative estimate of the admissible set $\cU_{EP}(\bx)$.
\end{remark}

A popular alternative to adjusting $\bar{q}$'s via regression standard errors is to replace the empirical percentile with the empirical conditional tail expectation (CTE):
\begin{align*}
    \overline{\text{CTE}}(\bx^i,u^i) & := \frac{\sum_{b=1}^{M_b} g^b\mathbbm{1}_{g^b \geq \bar{q}(\bx^i,u^i)} }{\sum_{b=1}^{M_b} \mathbbm{1}_{g^b \geq \bar{q}(\bx^i,u^i)}},
\end{align*}
and then running a regression on the training set $(\bx^i, u^i, \overline{\text{CTE}}(\bx^i,u^i))$ to obtain the fitted CTE surface
$\widehat{\text{CTE}}(\bx,u)$ and finally $\hat{\cU}_{CTE}(\bx) :=  \big\{ u: \widehat{\text{CTE}}(\bx,u) \le 0 \big\}$. This idea is similar to regularizing the Value-at-Risk estimation with the Conditional VaR.

\subsubsection{Support Vector Machines (SVM)}
For a fixed control $u$, finding the admissible set $\cX^a_{n}(u)$ in \eqref{eq:specificConstraintThroughPartition} can be interpreted as classifying each input $\bx$ as being in $\cX^a_{n}(u)$ or not. Therefore, we consider the use of  classification techniques, specifically support vector machines (SVM). This approach does not estimate the $(1-p)$-quantile $q(\bx,u)$, but rather its $0$-level set with respect to $(\bx,u)$.  The starting point is to use the nested Monte Carlo simulations to compute $\bar{p}(\bx^i,u^i)$ with much smaller batch size $M_b$ compared to Section~\ref{sec:nestedMC}. Next, we construct a binary classification objective with a training dataset $\{\bx^i,u^i,y^i\}_{i=1}^{M_a}$ where the $\pm 1$-labels are
\begin{align}
    y^i :=
\begin{cases}
    1,& \text{if } \bar{p}(\bx^i,u^i)<p;\\
    -1,              & \text{otherwise}.
\end{cases}
\label{eq:svm_noisy}
\end{align}
The boundary separating the two classes is evaluated by solving the optimization problem:
\begin{equation}
\label{eq:SVMsoft}
\min_{\vb\in\mathbb{R}^K}\Big\{\sum_{i=1}^{M_a} \Big(1-y^i [\vb^T\phi(\bx^i,u^i) + \beta_0] \Big)_+ + \frac{C}{2\cdot M_a}||\vb||^2\Big\},
\end{equation}
where $\phi(\bx,u) = \big[\phi_1(\bx,u),\phi_2(\bx,u),\ldots,\phi_K(\bx,u)\big]^T $ are the $K$ basis functions and $C$ is the penalty parameter. We estimate the set of admissible controls corresponding to $\bx$ as:
$$\hat{\cU}_{SVM}(\bx) :=  \Big\{ u: \hat{\vb}^T\phi(\bx,u) + \hat{\beta}_0 \ge 0 \Big\}.$$

Figure~\ref{fig:svmfit} displays the estimated $\hat{\cX}^a_n(u)$ and the corresponding dataset $(L^i, I^i,0, y^i)$ ($u = 0$ is fixed). The region where $u=0$ is admissible is to the left of the (thick red) decision boundary.

\begin{remark}
A conservative estimate $\hat{\cU}^{(\rho)}_{SVM}$ is obtained by biasing the decision boundary to the left by
re-labeling the training points in~\eqref{eq:svm_noisy} via:
 \begin{align}
    y^i=
\begin{cases}
    1,& \text{if } \bar{p}(\bx^i,u^i) + z_{\rho}\sqrt{\frac{ \bar{p}(\bx^i,u^i) (1- \bar{p}(\bx^i,u^i) )}{M_b}} < p\\
    -1,              & \text{otherwise}.
\end{cases}
\label{eq:svm_noisy_relabel}
\end{align}
\end{remark}

\subsubsection{Quantile Regression (QR)}

QR directly constructs a parametric model for $q(\bx,u)$:
\[
\hat{q}(\bx,u; \vb) := \sum_k \beta_k \phi_k(\bx,u).
\]
To estimate the coefficients $\vb \in \mathbb{R}^K$, we use the dataset $\{\bx^i,u^i,g^i \}_{i=1}^{M_a}$ (where $g^i $ is a sample from the distribution $G(\bx^i,u^i)$) to maximize the negative log likelihood:
\begin{align*}
\hat{\vb} & =\arg\min_{\vb\in\mathbb{R}^K}\bigg\{ \sum_{i=1}^{M_a} \mathcal{L}^{(p)}\Big(\textit{g}^i -\sum_{k=1}^K \beta_k\phi_k(\bx^i,u^i)\Big)\bigg\}, \\
\end{align*}
where $\mathcal{L}^{(p)}(y) =y(p-1_{\{y<0\}})=p y_+ + (1-p)y_-$.
As for the parametric density fitting, a transformation of $G(\bx,u)$ might be beneficial when applying quantile regression. Figure~\ref{fig:qrfit} presents the dataset $\{L^i,I^i,0,g^i\}_{i=1}^{M_a}$ (background colormap) and the estimated admissible set $\hat{\cX}^a(0)$ which is the region to the left of the
contour $\{\hat{q}_{QR}(L,I) = 0 \}$ (thick red line).

Relying on the Delta method again to compute the variance of the estimated quantile $\hat{q}(\bx,u; \hat{\vb})$ as $ \phi(\bx,u)'Var(\hat{\vb})\phi(\bx,u)$, the admissible set at $\bx$ at confidence level $\rho$ is:
$$\widehat{\cU}^{(\rho)}_{QR}(\bx) :=  \Big\{ u: \hat{q}(\bx,u; \hat{\vb}) + z_{\rho} \sqrt{\phi(\bx,u)^T Var(\hat{\vb})\phi(\bx,u) }\le0 \Big\}.$$

\section{Case Studies}
\label{sec:numerics}

Recall the problem introduced in Section~\ref{sec:motivationMicrogrid} where we control the operation of a diesel generator in order to match demand at  minimal cost while maintaining the probability of blackout between each decision epoch below a given threshold $p$. In this section, we discuss two variants of such microgrid control. In the first example, we assume a time-homogeneous net-demand process which reduces the problem of estimating admissible set to a pre-processing step. In the second example, we use time-dependent net demand process calibrated to data obtained from a microgrid in Huatacondo, Chile. Time-inhomogeneity requires to estimate the admissible set at every step. The microgrid features a perfectly efficient battery, so that the respective power output at $t_{n_k}$ (recall $t_{n_k}$ is a generic time instance on the finely discretized time grid) is given by:
 $$B_{n_k}=\max \Big(\min\Big(L_{n_k}-u_{n}, \frac{I_{n_k}}{\Delta n_k}  \Big), - \frac{I_{\max}-I_{n_k}}{\Delta n_k} \Big).$$
Table \ref{table:microgridParameters} lists other microgrid parameters, i.e.~capacity of the battery $I_{\max}$, maximum charging rate $B_{\min}$, maximum discharging rate $B_{\max}$ and diesel switching cost $\cK$.

\begin{table}[tbhp]
\caption{\footnotesize Parameters for the Microgrid example. }
\label{table:microgridParameters}
\centering
\begin{tabular}{c} 
\hline
$I_{\max} = 10$ (kWh), $B_{\min} = -6, B_{\max} = 6$ (kW), $\cK = 5$ \\ \hline
$T=48$ (hours), $\Delta t= 0.25$ (hours) \\ \hline
\end{tabular}
\end{table}

\subsection{Implementation details}

\textbf{Numerical Gold Standard:} In the absence of analytic benchmark, we use a high-budget gold standard to compare the output from the models discussed in Section~\ref{sec:admissibleSetEstimation}. For each fixed time-step $t_n$ we discretize the domain $\cX = (L,I)$ into $10,000$ design sites over a grid of $100\times 100$. For each design site $(L^i,I^j), \ i,j \in \{1,\ldots,100\}$ and $u^k \in  0 \cup \{1 = u_1,\ldots, u_{101}=10 \}$, we evaluate $\hat{p}(L^i,I^j,u^k)$ using \eqref{eq:phat} with batch size $M_b=10,000$. Thus, the total simulation budget is $100\times100\times102\times10000 \approx 10^{10}$. We then evaluate the local minimal admissible control $$u^{\min}_n(L^i,I^j) =  \min \left\{ u: \hat{p}(L^i,I^j,u)<p\right\} .$$ To evaluate $u^{\min}_n(L,I)$ at new sites we employ linear interpolation on the dataset $\{L^i,I^j,u^{\min}_n(L^i,I^j) \}_{i,j=1}^{100}$.

To estimate the continuation function, we use the piecewise continuous approximation of Section~\ref{sec:piecewise} with $M_I=15,M_u=15$ and $M_L=1500$ sites in $L$. The design $\CD^c$  is constructed as the Cartesian product $\{L^1, L^2, \ldots, L^{M_L}\}\times \{I^0, I^1, \ldots, I^{M_I}\} \times \{u^0, u^1, \ldots, u^{M_u}\} $, where $L^i, i = 1, \ldots, M_L$ are sampled uniformly from the interval $[L_{\min}, L_{\max}] = [-8, 8]$. The inventory $\{I^0, I^1, \ldots, I^{M_I}\}$ and control $\{u^1, \ldots, u^{M_u}\}$  discretizations  are equispaced in $[0, I_{\max}]$ and $[\underline{u}, \overline{u}]$ respectively.

\textbf{Low budget policies:}
We approximate the continuation value function $\cC$ using a piecewise continuous approximation with degree-3 in $L$ combined with interpolation in other dimensions. The design $\CD^c$ is the same as for the numerical gold standard with discretization levels $M_L = 1000, M_I=10,M_u=10$.

We approximate the admissible set $\cU$ using the methods described in Section~\ref{sec:admissibleSetEstimation} and compare the performance of each method by using a fixed set of $M'=20,000$ out-of-sample simulations. To address the discontinuity in $\cW = 0\cup [\underline{u},\overline{u}]$, we implement two separate statistical models to learn $\cU_n(\cdot)$. As an example, with logistic regression of Section~\ref{sec:logisticReg} we estimate two sets of parameters in equation~\eqref{eq:logistic}: the first one uses one-step paths generated with $u=0$ and a two-dimensional regression of $y^{i,(1)}$ against $(L^i,I^i)$. The second one uses design sites in the three-dimensional space $(L,I,u)$ where $u\in [1,10]$ and a 3-D regression of $y^{i,(2)}$ against $(L^i,I^i,u^i)$. For both regressions we choose a Sobol sequence space-filling experimental design: $\CD^a = (L^i,I^i) \in [-2,8]\times[0,10]$ when $u = 0$, and $\CD^a = (L^i,I^i,u^i) \in [-2,8]\times[0,10]\times [1,10]$ otherwise. The control space $[1,10]$ is discretized into 51 levels.

Additional parameters used for each method are specified in Table~\ref{table:design admissible set}. We found that Matern-3/2 kernels work better than \eqref{eqn:gp_kernel1} for smoothing $\bar{p}(L,I,u)$ (GPR) and $\tilde{p}(L,I,u)$ (PF) because the respective input-output maps feature steep transitions as a function of $(L,i,u)$. It is known that ``rougher'' kernels are better suited for such learning tasks compared to the $C^\infty$-smooth squared exponential kernel \eqref{eqn:gp_kernel1} by allowing the fitted $\hat{p}$ to have more ``wiggle room''. On the other hand, in the context of EP and CTE the input observations of $\hat{q}(L,I,u)$ and $\overline{\text{CTE}}(L,I,u)$ are quite smooth in $(L,i,u)$ and both GP kernel families perform equally well. The algorithms are implemented in \texttt{python 2.7}. We used ``\texttt{GaussianProcessRegressor}'' and ``\texttt{SVM.SVC}'' functions from \texttt{sklearn} library for GPR and SVM respectively. For LR and QR we used ``\texttt{Logit}'' and ``\texttt{quantile\_regression}'' functions from \texttt{statsmodels} library.

\begin{table}[tbhp]
\caption{\footnotesize{ Parameters for the estimation of the admissible sets for each method. We use total simulation budget of $10^5$ for all models except the Gold Standard.}}\label{table:design admissible set}
{\small
\begin{center}
\begin{tabular}{@{}llr@{}}
\toprule
Method &  Budget ($M_a \times M_b$) &   \multicolumn{1}{c}{Further parameters} \\ \midrule
Gaussian Process (GPR) &  $2000 \times 50 $ & Matern-3/2 kernel \\
Logistic Regression (LR) &  $10^5 \times 1$ &   Degree-2 polynomials\\
Parametric Density Fitting (PF) &  $2000 \times 50$ &  Truncated Gaussian, Matern-3/2 kernel \\
Empirical Percentile (EP) &   $1000 \times 100 $ &  Squared exponential kernel \\
Conditional Tail Expectation (CTE) &  $1000\times 100 $ &  Squared exponential kernel \\
Quantile Regression (QR) &   $10^5 \times 1$ &   Degree-4 polynomials \\
Support Vector Machine (SVM) & $2000 \times 50 $&   C =1, RBF kernel \\
Gold Standard (GS) & $10^6\times10^4 $ &  budget = $10^{10}$ \\ \bottomrule
\end{tabular}
\end{center}
}
\end{table}

\subsection{Example 1: Microgrid with Stationary Net Demand}
\label{sec:example1}

In this subsection, we assume time-homogeneous Ornstein-Uhlenbeck dynamics of the net demand process
\begin{equation}
    dL(t)  = -\lambda L(t) dt + \sigma \mathrm{d} B(t)   \ \implies L(t)  = L(0) e^{-\lambda t } + \sigma \int_{0}^t e^{-\lambda (t-s ) }\mathrm{d}B(s),
    \label{eq:stationary_netdemand}
\end{equation}
where $(B(t))$ is a standard Brownian motion. This scenario reduces the complexity of learning the probability constraints since we need to estimate the admissible set $\cU_0(\cdot)$ only once as a pre-processing step before starting the approximate dynamic programming scheme for the continuation values. The simplified setting offers a good testbed
to evaluate the performance of different admissible set estimation methods of Section~\ref{sec:admissibleSetEstimation}; we show that the relative  performance remains similar as we extend to more realistic dynamics in Section~\ref{sec:example2}. For this example, we assume the mean reversion parameter $\lambda = 0.5$ and volatility $\sigma = 2$.

\begin{figure*}[htbp]
        \centering
        \begin{subfigure}[b]{0.32\textwidth}
            \centering
            \includegraphics[width=\textwidth]{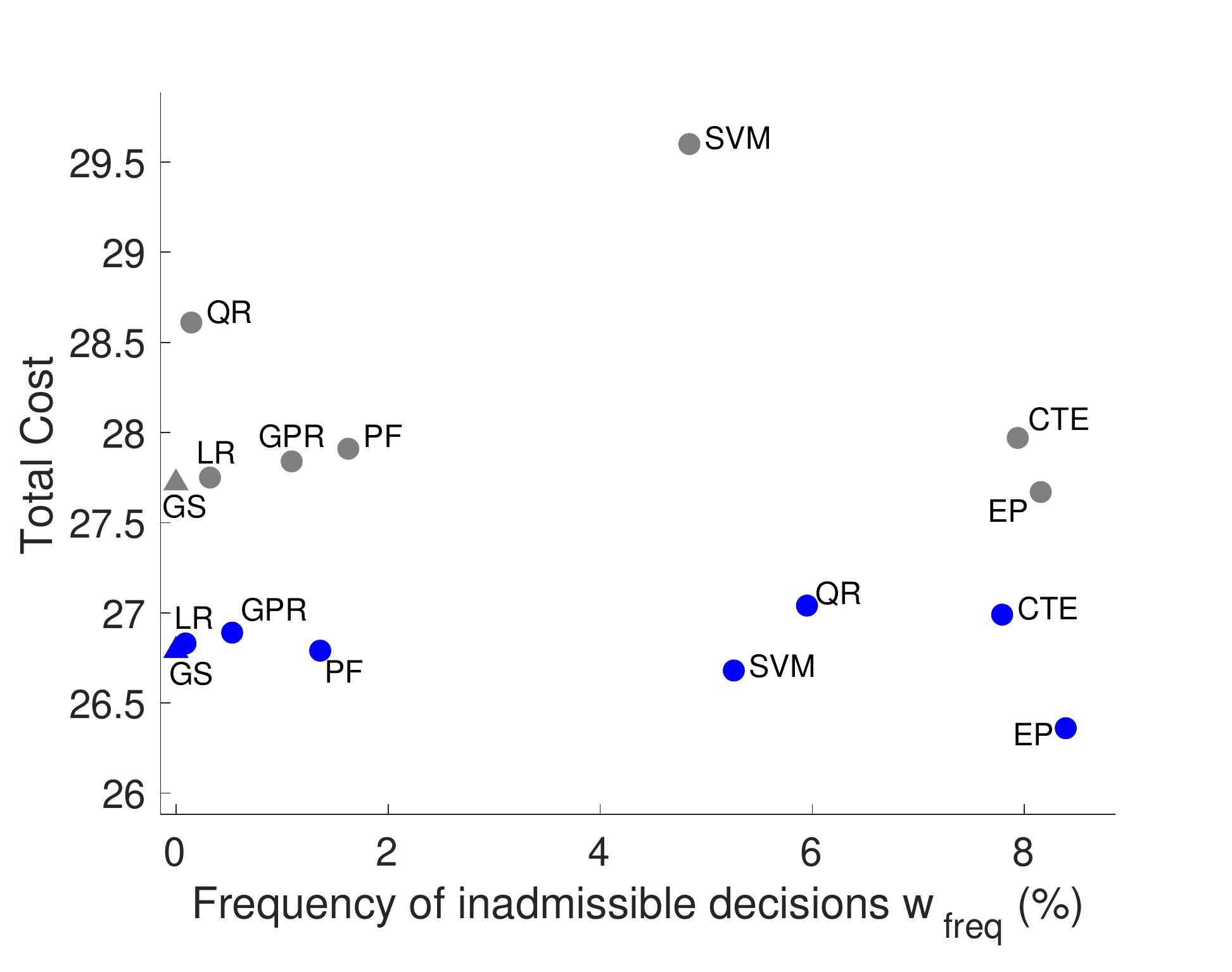}
            \caption[]{}
            \label{fig:cost_violations1}
        \end{subfigure}
        \begin{subfigure}[b]{0.32\textwidth}
            \centering
            \includegraphics[width=\textwidth]{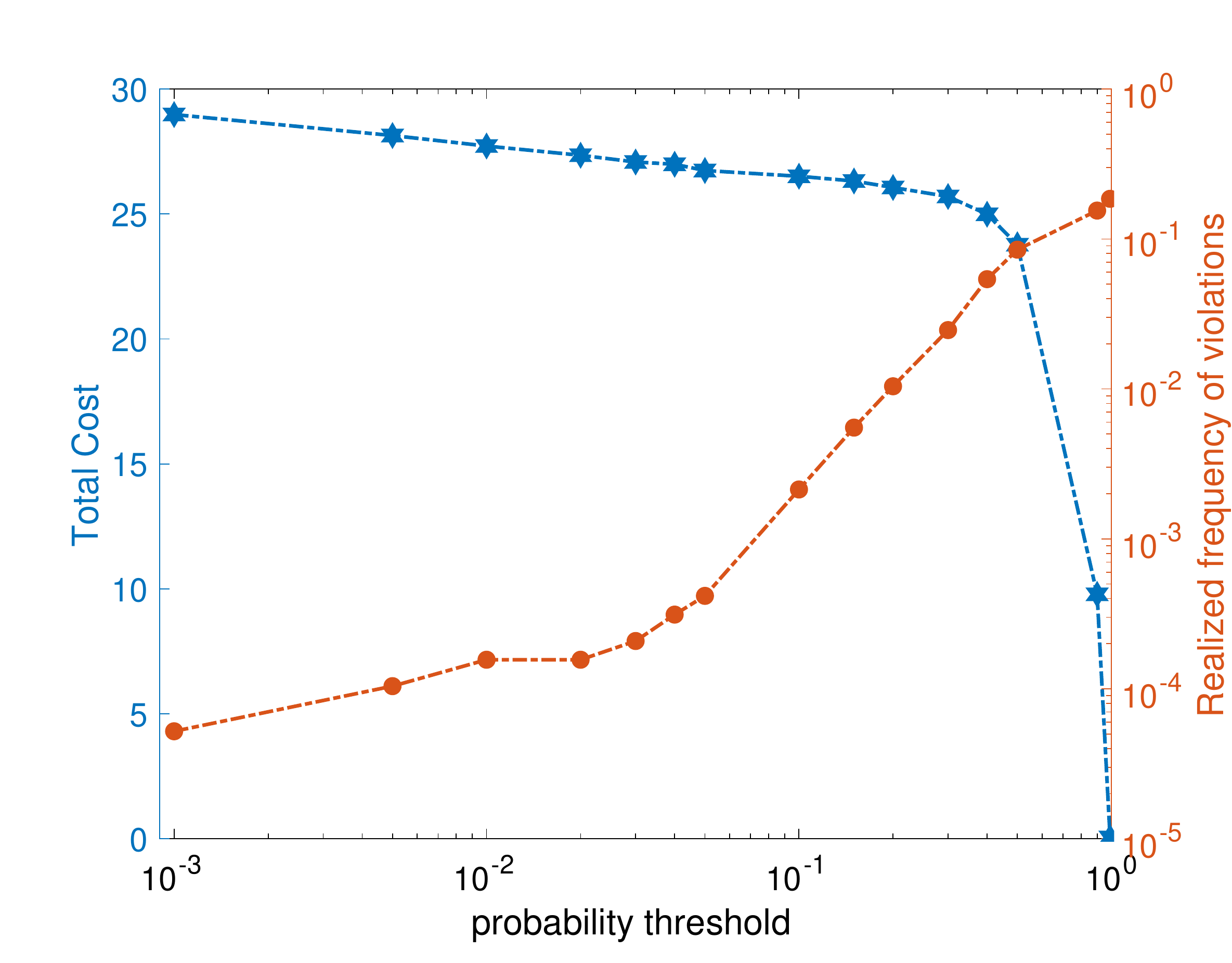}
            \caption[]{}
            \label{fig:Sensitivity}
        \end{subfigure}
        \begin{subfigure}[b]{0.32\textwidth}
            \centering
            \includegraphics[width=\textwidth]{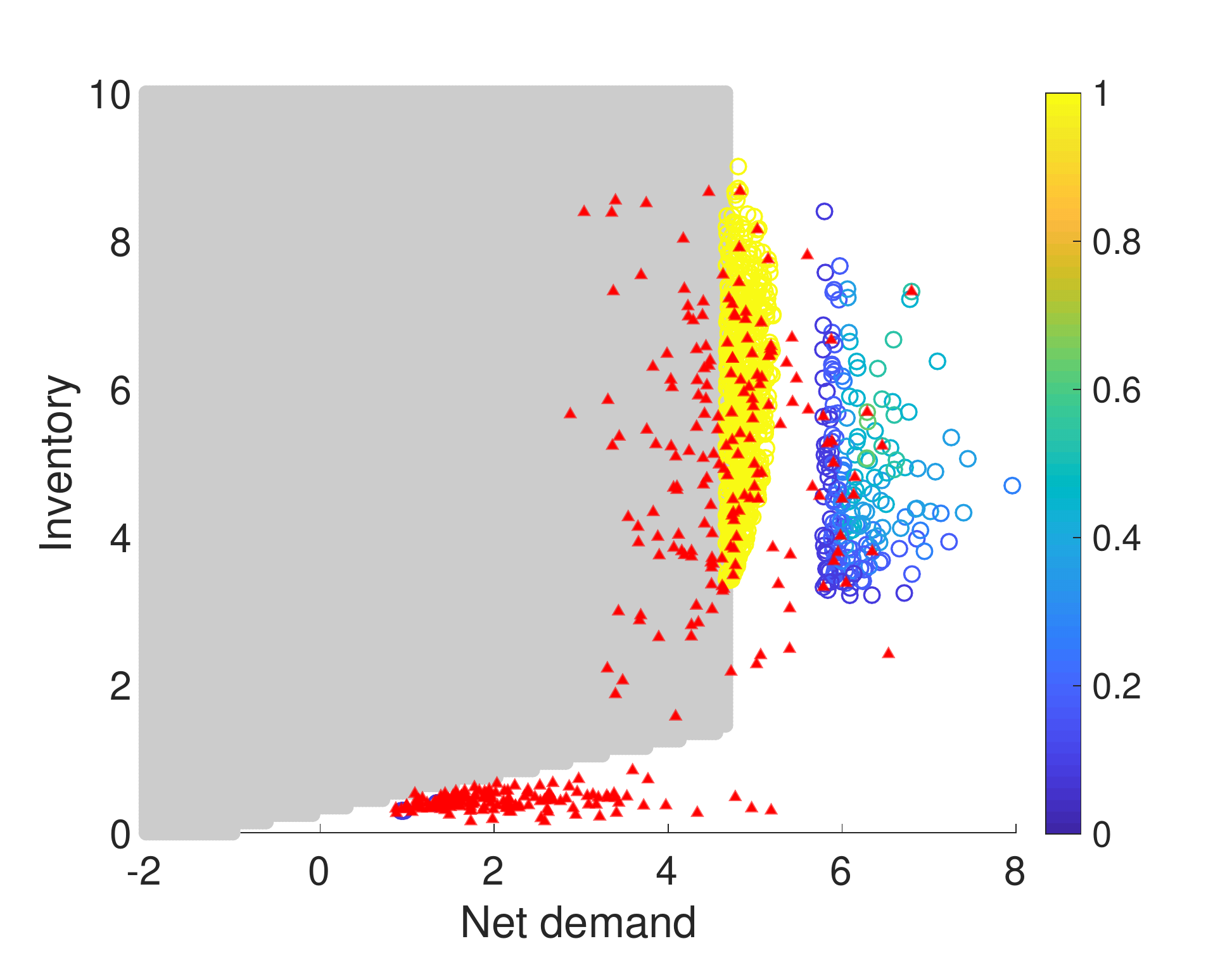}
            \caption[]{}
            \label{fig:realized_blackouts}
        \end{subfigure} \vspace*{-10pt}
        \caption[ ]
        {{ \footnotesize{ \emph{Left panel}: Trade-off between cost $\hat{V}_0(0,5)$ and frequency of inadmissible decisions $w_{freq}$ for the stationary model. Dark blue points correspond to $p=5\%$ probability constraint threshold and light grey ones to $p=1\%$. \emph{Center:} Total cost  $\hat{V}_0(0,5)$ (left axis, blue stars) and realized frequency of violations
        $w_{rlzd}$ (right axis, red circles) as functions of $p$ employing the LR model. \emph{Right:} Locations $(L,I)$ of realized violations
        ${\sup_{s \in [t_n,t_{n+1}) }S^{m'}(s)>0}$ (red triangles), inadmissible decisions
        ${ \hat{u}(n,\bx_n^{\hat{u},m'}) - u^{\min}_n(\bx_n^{\hat{u},m'}) <0 }$ (circles with color representing the inadmissibility margin) on 5000 out-of-sample simulations using LR model. The constraint is binding in the white region and is not binding in the grey region. }}}
        \label{fig:cost_violations}
    \end{figure*}

Figure~\ref{fig:cost_violations1} plots the resulting costs $\hat{V}_0(0,5)$ versus the frequency of inadmissible decisions $w_{freq}$ for different methods of Section~\ref{sec:admissibleSetEstimation}. We show the results both for  $p=0.05$ (dark blue), and $p=0.01$ (light grey) and benchmark  both cases against the numerical gold standard. Since the probabilistic constraints form the crux of the problem, we require schemes to maintain $\hat{u} \in \cU_n$ as much as possible, i.e., $w_{freq} \approx 0$. At $p=5\%$, we observe 0.09\% , 0.54\% and 1.36\% frequency of inadmissible decisions with logistic regression (LR), Gaussian process regression (GPR) and parametric density fitting (PF), respectively. Such accuracy might be deemed acceptable. However $w_{freq}$ is much worse (as high as 8.4\% with EP) for the other methods. While all the methods are a priori consistent, admissible set estimation via probability-based methods clearly seems to outperform quantile-based ones. Our experiments suggest that at low simulation budget, estimators of ${p}(\bx,u)$ have significantly lower bias compared to estimators of ${q}(\bx,u)$, thus partially explaining the difference. For a more stringent threshold $p=1\%$, we find the cost of all the methods to increase, without significant difference in the frequency of inadmissible decisions $w_{freq}$. Indeed, Figure~\ref{fig:cost_violations} illustrates the trade-off between lower costs and lower $w_{freq}$ (i.e.~more conservative estimate of the constraints).

Table~\ref{tab:example1} expands Figure~\ref{fig:cost_violations} by also reporting the corresponding $\tilde{\cT}$ statistic, the average inadmissibility margin $w_{avm}$ and realized frequency of violations (i.e.~blackouts) $w_{rlzd}$ defined as:
\begin{align}
w_{avm} &:= \blu{\frac{\sum_{n,m'}|\hat{u}_n(\bx_n^{\hat{u},m'}) - u^{\min}
_n(\bx_n^{\hat{u},m'}) | \mathbbm{1}_{ \hat{u}_n(\bx_n^{\hat{u},m'}) - u^{\min}_n(\bx_n^{\hat{u},m'}) <0 }}{\sum_{n,m'} \mathbbm{1}_{ \hat{u}_n(\bx_n^{\hat{u},m'}) - u^{\min}_n(\bx_n^{\hat{u},m'}) <0 };}} \label{eq:wav} \\
w_{rlzd}&:=\frac{1}{N\cdot M'}\sum_{n,m'}\mathbbm{1}_{\sup_{s \in [t_n,t_{n+1}) }S^{m'}(s)>0}. \label{eq:wblackouts}
\end{align}
We  find the realized frequency of violations $w_{rlzd}$ to be lowest for LR, GPR and PF. The  average inadmissibility margin $w_{avm}$ is also lowest for GPR and PF (the large value of  $w_{avm}$  for LR is attained in very small region as evident from $w_{freq}\approx 0$). The $\tilde{\cT}$ statistic is negative  for LR, GPR and PF and positive for the rest, meaning that all other methods fail to statistically respect the probability constraints when binding. Due to small frequency of inadmissible decisions $w_{freq}$, cost $\hat{V}_0(0,5)$ similar to the numerical gold standard and negative test statistic $\tilde{\cT}$, we recommend LR, GP and PF methods for the problem at hand.

Next, we test the sensitivity of the cost in terms of the probability threshold $p$  (employing logistic regression $\hat{\cU}_{LR}$) in Figure~\ref{fig:Sensitivity}. Increasing $p$ decreases $V$ as the set of admissible controls $\cU$ monotonically increases in $p$. For example, any admissible control at $p=1\%$ threshold is also feasible for $p > 1\%$, thus the respective cost at $1\%$ will be at least as much as at, say, $10\%$ threshold.

As previously discussed, the constraint is binding for only approximately $10\%$ of time-steps. In fact, that probability varies across the methods since the estimate of $\hat{\cU}$ affects the choice of $\hat{u}_n$ and ultimately the \emph{distribution} of $\hat{X}_n$. Intuitively, the realized system states are driven by the estimates of the probabilistic constraints. Typically, more conservative estimates of $\cU$ will push $\hat{X}_{0:N}$ away from the ``risky'' regions. This is also confirmed in Figure~\ref{fig:cost_violations} where as $p \rightarrow 1$, $w_{rlzd} \rightarrow 20\% = w_{bind}$ while in Table~\ref{tab:example1} $w_{bind} \simeq 10\%$.

The metrics $w_{freq}$~\eqref{eq:w_freq},  $w_{bind}$~\eqref{eq:test2}, $w_{rlzd}$~\eqref{eq:wblackouts} are closely linked. As the inadmissible decisions can occur only when the constraint is binding, $u^{\min} > 0$, we expect $w_{freq} \leq w_{bind}$ and $w_{freq} \approx w_{bind}$ for a method with a bias in overestimating the admissible set (e.g.~$\cX^{a,EP}(u) \supset \cX^{a,GS}(u) \ \forall u \in\cW$). The realized violations (blackouts)
$w_{rlzd}$ can be represented as a sum of three:
\[
w_{rlzd} = p_1 w_{freq}+ p_2 (w_{bind} - w_{freq}) + p_3 (1 - w_{bind}), \qquad p_1+p_2+p_3 = 1,
\]
where the weights $p_1, p_2, p_3$ depend on the distribution of the controlled trajectories.
The first term represents the instances when the  constraint is binding but the controller chooses an inadmissible control (i.e.~mis-estimates $\hat{\cU}$). The second term represents instances when the constraint is binding and correctly estimated, but due to random shocks violations take place (with a conditional frequency below the specified $p=0.05$). The last term represents instances when the constraint is not binding but some violations still occur with the intrinsic conditional frequency strictly less than $p$. Note that due to $w_{bind} \ll 1$, most of the violations are of the latter type, i.e.~take place when $u^* = 0$ and the conditional violation probability is below $p$.  We illustrate these scenarios in Figure~\ref{fig:realized_blackouts} using the LR model. Thus, the first term counts the instances when violations occur at the same time as controller makes an inadmissible decision (circle encircling triangle), the second term counts  the triangles when $I \approx 0$, and the third term the triangles in the grey region where the constraint is not binding (violations
when $u^{\min} = 0$).

Although we observed poor performance of quantile based methods, asymptotically (with respect to the simulation budget) we expect them to perform similar to the probability based methods. As an example, in Appendix~\ref{sec:appendix_simulation_budget} Table~\ref{tab:svm_budget}, we present the performance of SVM for thresholds $p=5\%$ and $p=1\%$ with increasing budget. For $p=5\%$ and by increasing the simulation budget from $10^5$ to $10^8$, we find the frequency of inadmissible decisions $w_{freq}$ to drop from $5.93\%$ to $1.5\%$, average inadmissibility amount $w_{avm}$ from $0.78$ kW to $0.27$ kW,  frequency of realized blackouts $w_{rlzd}$ from $2.80\%$ to $0.30\%$ and the test statistic  which rejected the method at $10^5$ simulation budget ($\cT \gg 0$) suggests to accept it ($\cT \ll 0$) at $10^8$ simulation budget. We observe similar behavior at $p=1\%.$ 
\blu{The main challenge with quantile-based methods is the underlying bias in learning $q(\mathbf{x}, u)$. This bias is known to converge to zero slowly, necessitating a relatively large $M_b$. A glimpse of this can be observed in Appendix~A where SVM does not perform adequately all the way up to total budget of $10^8$. Increasing $M_b$ with $M_a$ fixed offers only a limited improvement that tapers off quickly, as was observed previously in\cite{aditya2018}. Keeping a constant budget, a high $M_b$ forces a low $M_a$ which offsets these performance gains as the regression is not  able to properly explore the space.}

\blu{There are many techniques to improve $\bar{q}$, for example variance reduction tools. We have experimented with variance reduction for the SVM method by using antithetic variables. Although the results were marginally better, they were not meaningfully different from those reported in Appendix~A. A completely different way to improve EP and SVM could be to judiciously choose the simulation design $\CD^a$. Paper~\cite{aditya2018} explores various approaches to build simulation designs for estimating the continuation value function and documents their strong impact on performance. We anticipate that a similar strategy for $\CD^a$ may improve convergence.}

\begin{table}[tbhp]
{
\caption{\footnotesize Cost of running the microgrid $\hat{V}_0(0,5)$, frequency of inadmissible decisions
$w_{freq}$, average inadmissibility margin $w_{avm}$, realized frequency of violations (i.e. blackouts) $w_{rlzd}$, test statistic $\tilde{\cT}$ and frequency of binding constraint $w_{bind}$ for the example in Section~\ref{sec:example1}.}\label{tab:example1}
\begin{center}
\begin{tabular}{@{}lcccccc@{}}
\toprule
Method & $\hat{V}_0(0,5)$ (\$) & $w_{freq}$ (\%) & $w_{avm}$ (kW) & $w_{rlzd}$ (\%) & $\tilde{\cT}$ & $w_{bind}$ (\%) \\ \midrule
GS & 26.79 & 0.00 & 0.00 & 0.37 & - & - \\
LR & 26.83 & 0.09 & 0.82 & 0.03 & -125 & 8.69 \\
GPR & 26.89 & 0.53 & 0.16 & 0.11 & -98 & 8.10 \\
PF & 26.79 & 1.36 & 0.27 & 0.21 & -69 & 8.51 \\
SVM & 26.68 & 5.26 & 0.55 & 1.83 & 388 & 9.67 \\
QR & 27.04 & 5.95 & 0.33 & 0.98 & 145 & 9.49 \\
CTE & 26.99 & 7.79 & 0.43 & 1.63 & 320 & 9.93 \\
EP & 26.36 & 8.39 & 0.49 & 1.98 & 403 & 10.45 \\ \bottomrule
\end{tabular}
\end{center}
}
\end{table}

\textbf{Conservative estimators for $\cU$.} Algorithms for SCPC are expected to respect the probabilistic constraints, so that it is critical to minimize the occurrence of inadmissible decisions. As discussed in Section~\ref{sec:admissibleSetEstimation}, one way to raise the statistical guarantee for admissibility of  $\hat{\cU}$ is by adding a margin of error $\xi(\bx,u)$. This yields a more conservative (i.e.~smaller) $\hat{\cU}$ and hence lowers $w_{freq}$. In Table~\ref{tab:example1_conservative} we examine three scenarios for $\xi(\bx,u)$ sorted from least to most conservative (in all cases we maintain $p=5\%$):
\begin{itemize}
\item Scenario 1: unadjusted $\xi = 0\%$ (same as Table~\ref{tab:example1});
\item Scenario 2: $\xi^{(\rho)}(\bx,u)$ at 95\% confidence level, $z_{\rho} = 1.96$;
\item Scenario 3: fixed $\xi = 4\%$, which is equivalent to lowering the violation threshold to $p - \xi = 1\%$.
\end{itemize}

\begin{figure*}[htb]
        \centering
        \begin{subfigure}[b]{0.44\textwidth}
            \centering
            \includegraphics[width=\textwidth]{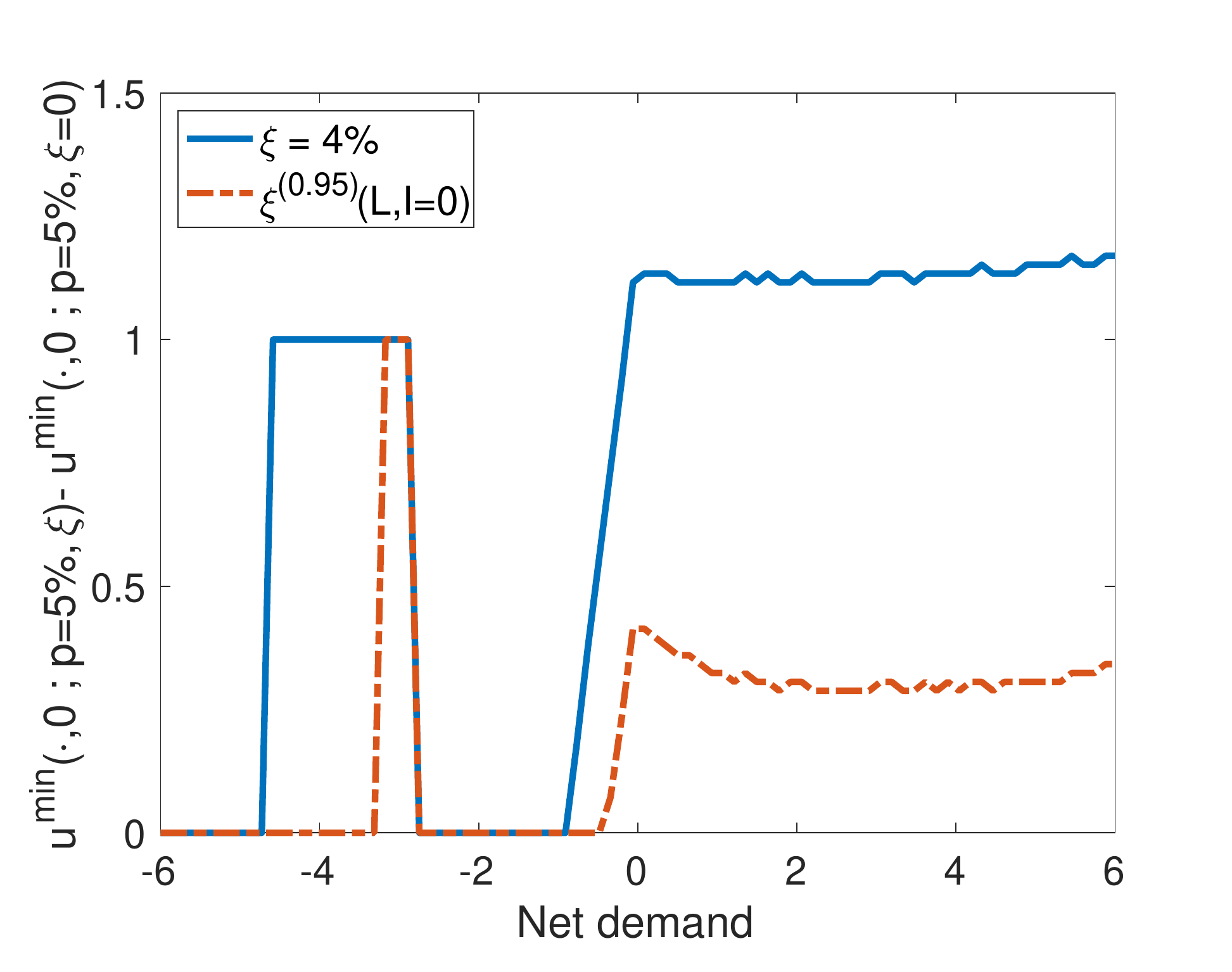}
            \caption[]%
            {{  LR}}
            \label{fig:Logistic_scenarios}
        \end{subfigure}
        \begin{subfigure}[b]{0.44\textwidth}
            \centering
            \includegraphics[width=\textwidth]{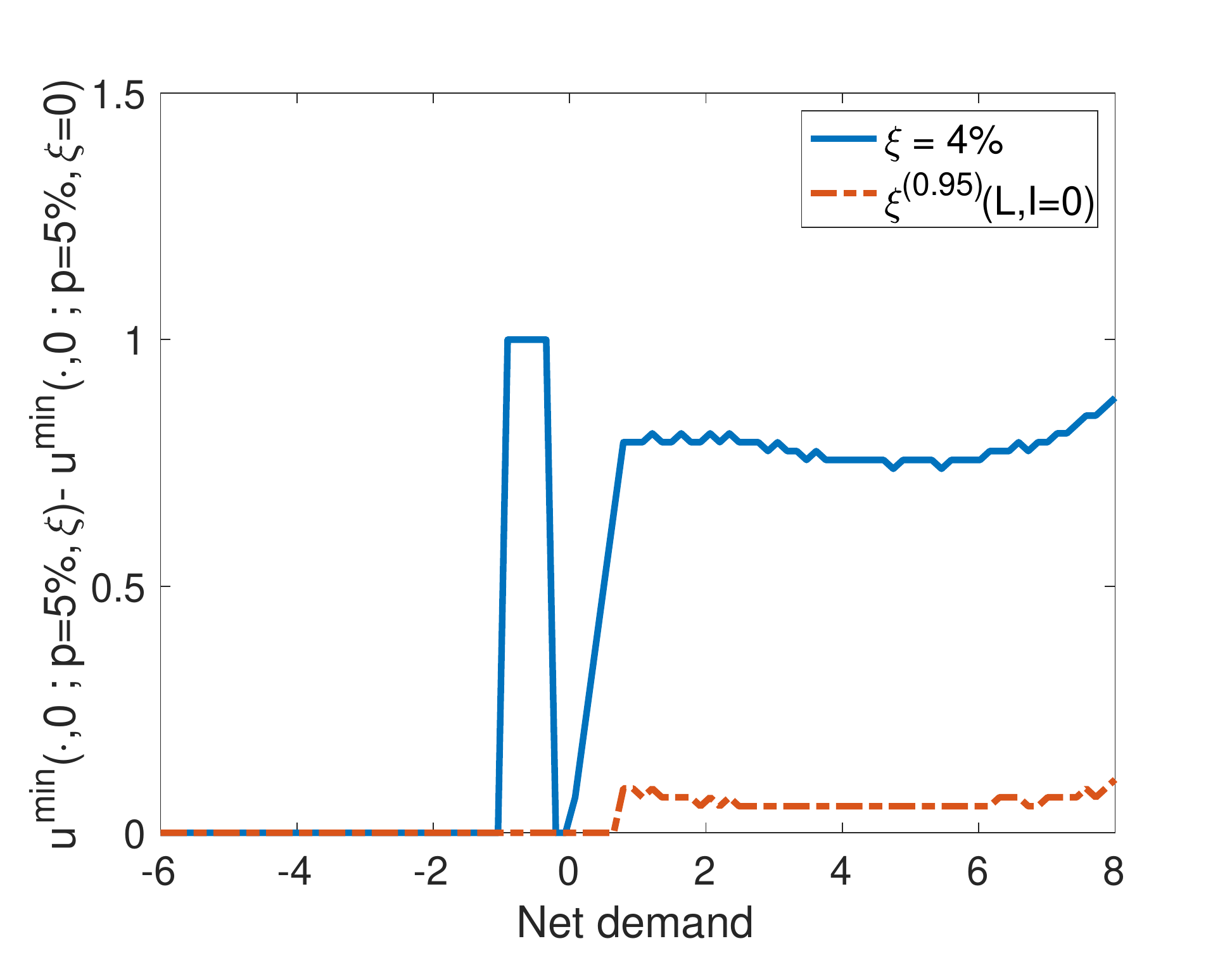}
            \caption[]%
            {{ QR}}
            \label{fig:Quantile_scenarios}
        \end{subfigure}  \vspace*{-10pt}
        \caption[ ]
        {{\small Impact of the margin of error $\xi(\cdot, \cdot)$ on minimum admissible control $\hat{u}^{\min}$. We plot the difference between minimum admissible control for scenario 2 ($\hat{u}^{\min}(\cdot,I; \xi^{(0.95)}(L,I) )$) and scenario 3 ($\hat{u}^{\min}(\cdot,I; \xi = 4\% )$) with respect to scenario 1 ($\hat{u}^{\min}(\cdot,I; \xi = 0\% )$) using LR (left panel) and QR (right panel) models. }}
        \label{fig:adaptive_non_Adaptive}
    \end{figure*}

Table~\ref{tab:example1_conservative} confirms the intuition that the frequency of inadmissible decisions $w_{freq}$ should be decreasing from scenario 1 to 3. This is illustrated in Figure~\ref{fig:adaptive_non_Adaptive} that shows how  the minimum admissible control is affected by $\xi(\bx,u)$. Although adding a margin of error does lower $w_{freq}$, this mechanism does not really alter the relative performance of the different methods. Thus,  for all three scenarios, we find
SVM, CTE and EP to give unreliable estimates of $\cU$ (since $\tilde{\cT} \gg 0$). An exception is QR which yields high $w_{freq}$ for $\xi=0$ but does become acceptable ($\tilde{\cT} < 0$) in scenario 3. In contrast, LR, GPR and PF perform well throughout. Table~\ref{tab:example1_conservative1} in Appendix~\ref{sec:appendix_rho} provides additional results as we vary the confidence level to $\rho = 90\%, 99\% \text{ and } 99.95\%$, with the same general conclusions. (Observe that a $w_{freq}$ very close to zero likely implies that $\hat{\cU} \subset \cU$ is strictly smaller and the controller is so conservative as to rule out some admissible actions.)

We generally expect the ultimate cost $\hat{V}_0(0,5)$ to increase as $\hat{\cU}$ becomes more conservative, see the estimated $\hat{V}$'s across each row of Table~\ref{tab:example1_conservative}. The increase in costs arises  due to two factors: when the diesel generator is started sooner (due to $u = 0$ becoming inadmissible as $\xi$ is raised), and the higher level of $\hat{u}$ once the diesel is ON. This can be seen in Figure~\ref{fig:adaptive_non_Adaptive} where in Scenarios 2 and 3 the controller switches the generator at a lower net demand and once the diesel is running picks a higher power output ($\hat{u}^{\min}(\cdot,I; p=5\%, \xi ) - \hat{u}^{\min}(\cdot,I; p=5\%, \xi=0 )> 0$). We stress that the link between $\hat{\cU}$ and $\hat{V}$ is complicated by the fact that as $\hat{\cU}$ changes, so does the distribution of the controlled paths. So for example in Table~\ref{tab:example1_conservative} the cost for QR falls in Scenario 2, although it remains within two Monte Carlo standard errors.

\begin{table}[tbhp]
\caption{\footnotesize Impact of margin of error $\xi$ on the estimated cost of running the microgrid $\hat{V}_0(0,5)$, frequency of inadmissible decisions
$w_{freq}$, and test statistic $\tilde{\cT}$ from \eqref{eq:test1}. The probabilistic constraint is $p=5\%$. }
\centering
{
\begin{tabular}{@{}lccc|ccc|ccc@{}}
\toprule
 & \multicolumn{3}{c|}{$\xi = 0\% $} & \multicolumn{3}{c|}{$\xi^{(0.95)}(\bx,u)$}   & \multicolumn{3}{c}{$\xi = 4\%$} \\ \midrule
Method & $\hat{V}_0(0,5)$ & $w_{freq}$  & $\tilde{\cT}$ & $\hat{V}_0(0,5)$ & $w_{freq}$ & $\tilde{\cT}$ & $\hat{V}_0(0,5)$ & $w_{freq}$  & $\tilde{\cT}$ \\ \midrule
GS & 26.79 & 0.00 & - & - & - & - & - & - & - \\
LR & 26.83 & 0.09 & -125 & 26.95 & 0.08 & -124 & 27.86 & 0.04 & -112 \\
GPR & 26.89 & 0.53 & -98 & 28.00 & 0.01 & -110 & 28.12 & 0.00 & -107 \\
PF & 26.79 & 1.36 & -69 & - & - & - & 27.91 & 0.44 & -96 \\
SVM & 26.68 & 5.26 & 388 & 29.65 & 3.41 & 225 & 29.60 & 3.41 & 225 \\
QR & 27.04 & 5.95 & 145 & 26.89 & 5.17 & 72 & 28.61 & 0.00 & -117 \\
CTE & 26.99 & 7.79 & 320 & 27.36 & 7.52 & 274 & 28.44 & 6.83 & 248 \\
ER & 26.36 & 8.39 & 403 & 26.97 & 7.78 & 225 & 28.13 & 7.08 & 283 \\ \bottomrule
\end{tabular}
 } 
\label{tab:example1_conservative}
\end{table}

\textbf{Take-aways.} Our experiments demonstrate the following: (i) To accurately estimate
admissible sets of the form~\eqref{eq:specificConstraint} we recommend to use LR, GPR or PF which all model the underlying probability of violations $p(\bx,u)$. Although asymptotically equivalent, the approach of quantile estimation leads to poor estimates $\hat{\cU}_n$ for practical budgets. (ii) Frequency of inadmissible decisions can be partly controlled by using conservative estimates $\hat{\cU}^\xi_n$ at the expense of higher costs. However, even a conservative $\hat{\cU}^\xi$ fails to make quantile-based methods acceptable, except for QR. (iii) For a new application, our suggested approach is to first evaluate the test statistic $\tilde{\cT}$ at $\xi=0\%$ using one of the recommended methods. Depending on how close is $\tilde{\cT}$ to zero, one can then adjust $\hat{\cU}$' via $\xi$ to improve the statistical guarantees on the frequency of inadmissible decisions $w_{freq}$.

\subsection{Example 2: Microgrid with seasonal demand}
\label{sec:example2}

\begin{figure*}[htb]
        \centering
        \begin{subfigure}[b]{0.44\textwidth}
            \centering
            \includegraphics[width=\textwidth]{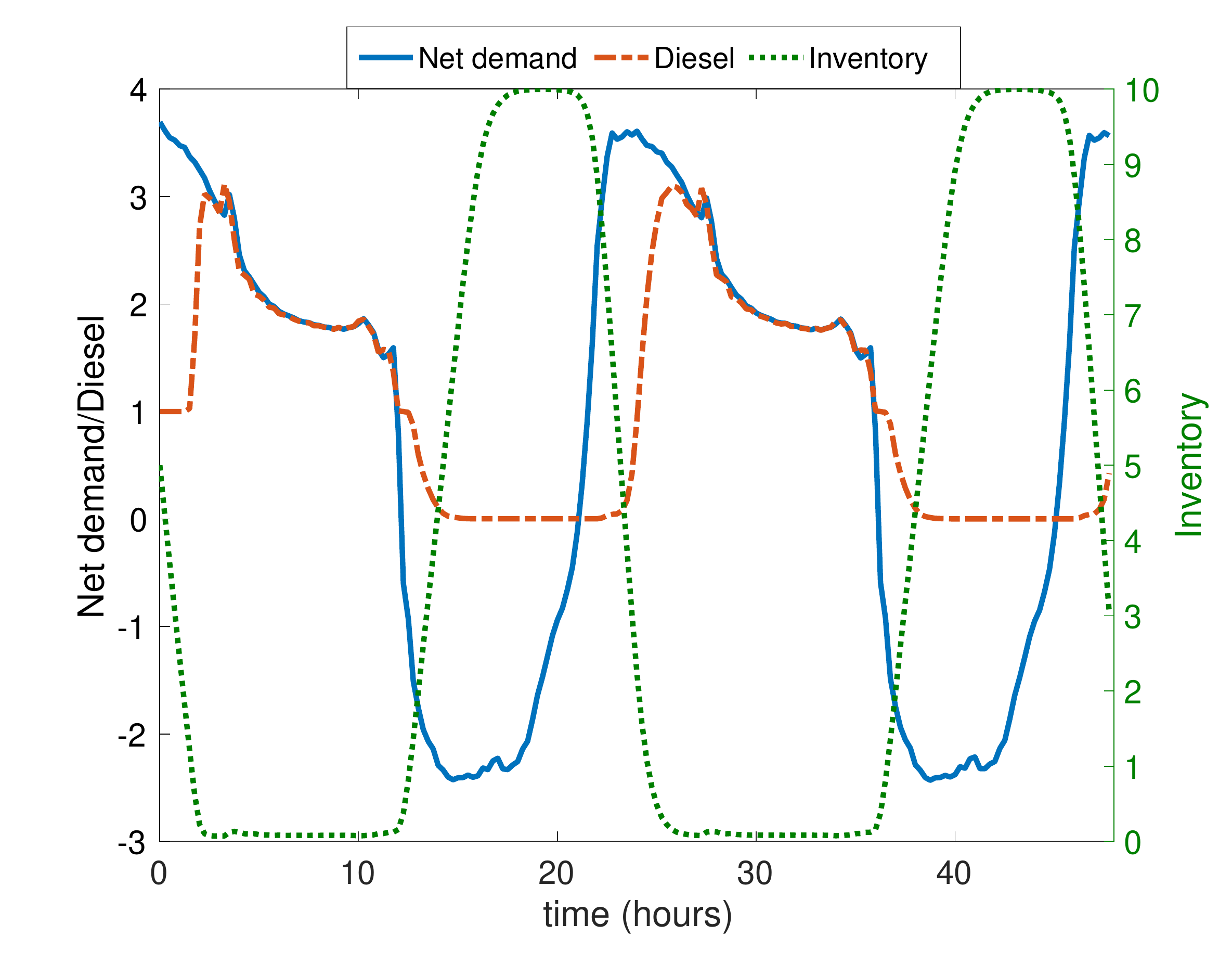}
        \end{subfigure}
        \quad
        \begin{subfigure}[b]{0.44\textwidth}
            \centering
            \includegraphics[width=\textwidth]{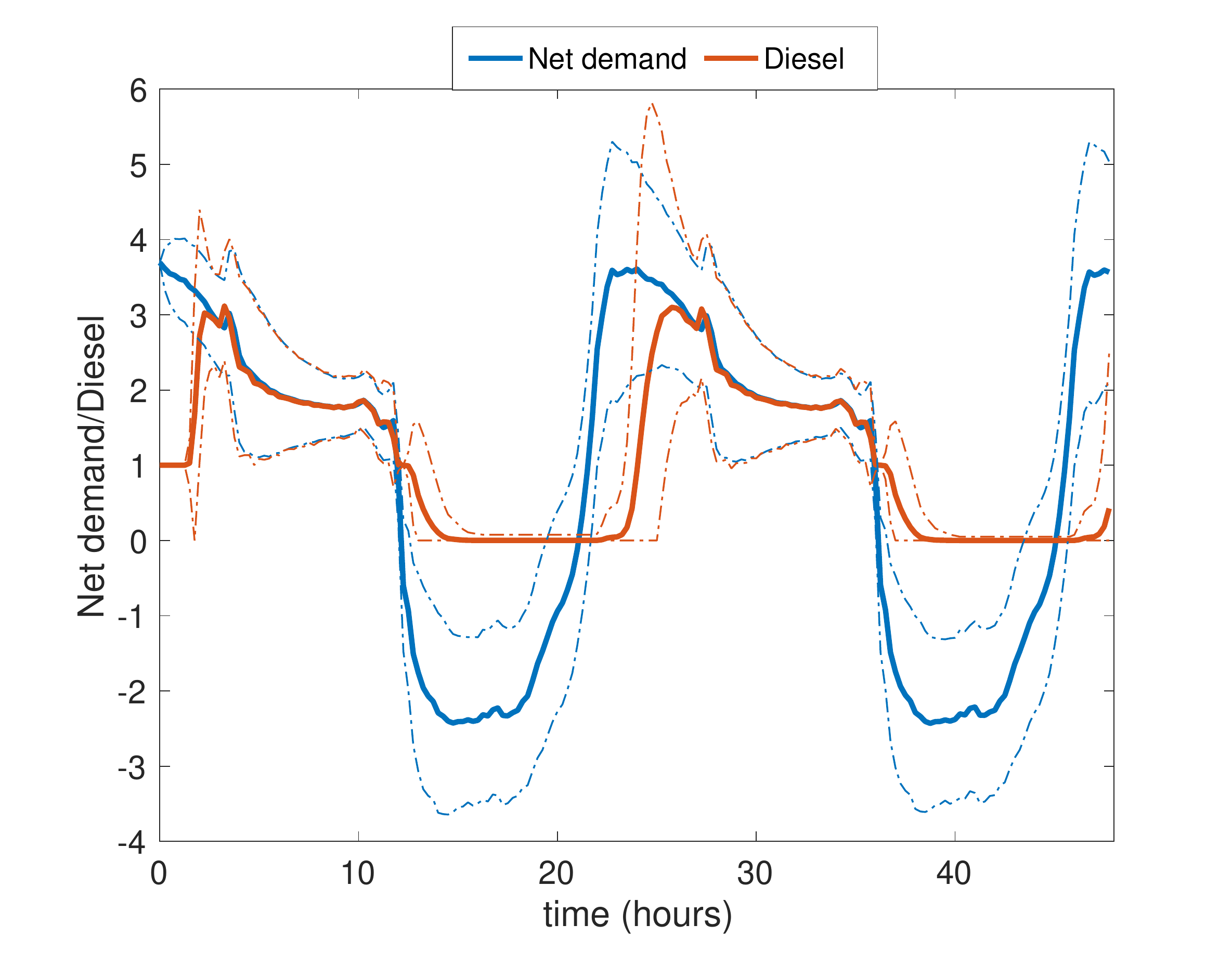}
        \end{subfigure}
        \caption[ ]
        { \footnotesize Model parameters, average trajectory of the state variables, control and their variance. Left panel: Average values of net demand $\frac{1}{M'}\sum_{m'=1}^{M'}L_n^{\hat{u},m'}$, inventory $\frac{1}{M'}\sum_{m'=1}^{M'}I_n^{\hat{u},m'}$ and optimal control (diesel) $\frac{1}{M'}\sum_{m'=1}^{M'}\hat{u}_n^{m'}$ processes using the gold standard strategy. Right panel: 95\% confidence bands for net demand $L_n^{\hat{u}}$ and  realized optimal diesel control  $\hat{u}_n$. \blu{Net demand and diesel output is measured in kW and Inventory in kWh.}}
        \label{fig:mean_path}
    \end{figure*}

Unlike the previous example, where we assumed time-homogeneous net demand, in practice there is seasonality: during the day renewable generation is high and net demand is often negative; during morning/evening demand exceeds supply making $L(t) > 0$. To incorporate this seasonality we use time-dependent Ornstein Uhlenbeck process (see \cite{heymann16} for a similar microgrid control problem):
\begin{align}\label{eq:seasonal-ou}
dL(t) = \left[\frac{\partial \mu}{\partial t}(t) +  \lambda\big(\mu(t) - L(t)\big)\right]dt + \sigma(t) \mathrm{d} B(t).
\end{align}
Here, $\lambda$ represents the speed of mean reversion towards the seasonal mean $\mu(t)$, while $\sigma(t)$ represents the time-varying volatility. Using It\^o's lemma and integration by parts one can prove that
\[
L(t)=  \mu(t) + e^{-\lambda t}\big(L(0) - \mu(0)\big) + \int_{0}^t e^{-\lambda(t-s)}\sigma(s) \mathrm{d} B(s).
\]
 Thus,
 \[
 \mathbb{E}[L(t)] =  \mu(t) + e^{-\lambda t}(L(0) - \mu(0)).
 \]
We calibrate $\mu(t)$ and $\sigma(t)$ in \eqref{eq:seasonal-ou} using iterative methodology described in \cite{heymann16} and the data from a solar-powered microgrid in Huatacondo, Chile\footnote{\blu{https://microgrid-symposiums.org/microgrid-examples-and-demonstrations/huatacondo-microgrid/}}. Specifically, we compute the mean and variance of the residual demand over 24 hours at 15-minute intervals using data from Spring 2014, i.e.~compute  $\{\mu_1,\mu_2,\ldots,\mu_{96} \}$ and $\{\sigma_1,\sigma_2,\ldots,\sigma_{96} \}$. The estimated $\mu(t)$ can be seen in the left panel of Figure~\ref{fig:mean_path} that plots the empirical average of $L(t)$. As expected, during the day, i.e., $t \in [12, 20]$ (noon-8:00 pm), the expected net-demand is negative ($\mu(t) < 0$) while it is positive ($\mu(t) > 0$) in the morning and during the night. The volatility $\sigma(t)$ is higher during the day due to the intermittent and unpredictable nature of solar irradiance. The mean reversion parameter was estimated to be $\lambda = 0.3416$.

To visualize the interplay of the net demand, inventory and optimal control, the left panel of Figure~\ref{fig:mean_path} presents the average trajectories of the three processes over 48 hours. During the morning hours when the demand $L(t)$ is high and the battery is empty, the controller uses the diesel generator. During the day when the renewable output is high and $L(t)$ is negative, the controller switches off the diesel and the battery charges itself. However, the non-trivial region is when the average net-demand changes sign, either from positive to negative around noon or negative to positive in the evening. During the former time interval, the optimal control process is in $\{0,1\}$ (recall that minimum diesel output is $1$). Similarly, during the evening when the net demand becomes positive (as the renewable output declines), the controller quickly ramps up the diesel to match $L(t) \gg 0$. The right panel of Figure~\ref{fig:mean_path} repeats the average control and inventory curves, also showing their 2-standard deviation bands (in terms of the out-of-sample trajectories of $\hat{L}^{\hat{u}}_{0:T}$). As expected, the time periods around ramp-up or ramp-down of the diesel generator is when $\hat{u}_n$ experiences the greatest path-dependency and dispersion and differs most from the demand curve.

Comparing Table~\ref{tab:example2}, which lists the estimated cost $\hat{V}_0(\mu(0),5)$ along with related statistics, with Figure~\ref{fig:cost_violations} indicates that incorporating seasonal net-demand process does not change the relative order of performance between the methods. The cost goes up as the diesel generator has to be used throughout the mornings and the evenings to match demand.

\begin{table}[tbhp]
\caption{\small Cost of running the microgrid $\hat{V}_0(\mu(0),5)$, frequency of inadmissible decisions
$w_{freq}$, average inadmissibility margin $w_{avm}$, realized frequency of  violations $w_{rlzd}$ and frequency of the constraint being binding $w_{bind}$ for the case study in Section~\ref{sec:example2}.}
\label{tab:example2}
\centering
\begin{tabular}{@{}lcccccc@{}}
\toprule
Method & $\hat{V}_0(\mu(0),5)$ & $w_{freq}$ (\%) & $w_{avm}$ (kW) & $w_{rlzd}$ (\%) & $\tilde{\cT}$ & $w_{bind}$ (\%) \\ \midrule
GS & 53.38 & 0 & 0 & 0.30 & - & - \\
LR & 53.78 & 0.03 & 0.79 & 0.01 & -301 & 45.2 \\
GPR & 54.04 & 1.17 & 0.14 & 0.19 & -220 & 31.0 \\
PF & 54.55 & 0.02 & 0.26 & 0.01 & -226 & 25.7 \\
SVM & 40.52 & 43.37 & 0.91 & 43.37 & 5,306 & 46.4 \\
QR & 52.56 & 42.87 & 0.28 & 38.41 & 4,772 & 46.3 \\
CTE & 53.02 & 21.62 & 0.21 & 10.43 & 1,079 & 46.0 \\
EP & 52.82 & 21.91 & 0.23 & 11.57 & 1,227 & 46.1 \\ \bottomrule
\end{tabular}
\end{table}

As in the previous example, the performance of LR, GPR  and PF almost matches the gold standard despite significantly lower simulation budget. In this setting the constraint is binding approximately 45\% of the time (except for GPR and PF where it is 30\% and 25\% of the time). Frequency of inadmissible decisions $w_{freq}$ is 0.03\% for LR, 1.17\%  for GPR, and 0.02\% for PF. In contrast $w_{freq}$ is 43\% for QR, 22\% for EP, 43\% for SVM and 22\% for CTE, implying that all these schemes are highly unreliable for learning $\hat{\cU}$. The average inadmissibility margin $w_{avm}$ is also significantly lower for GPR ($0.14$ kW) and PF  ($0.26$ kW)  compared to the rest of the methods. Here again we observe larger inadmissibility margin and very low frequency of inadmissible decisions for logistic regression. Similar behavior is also evident for the test statistic $\tilde{\cT}$ and realized frequency of violations $w_{rlzd}$.

\begin{figure*}[tb]
        \centering
            \includegraphics[width=0.32\textwidth]{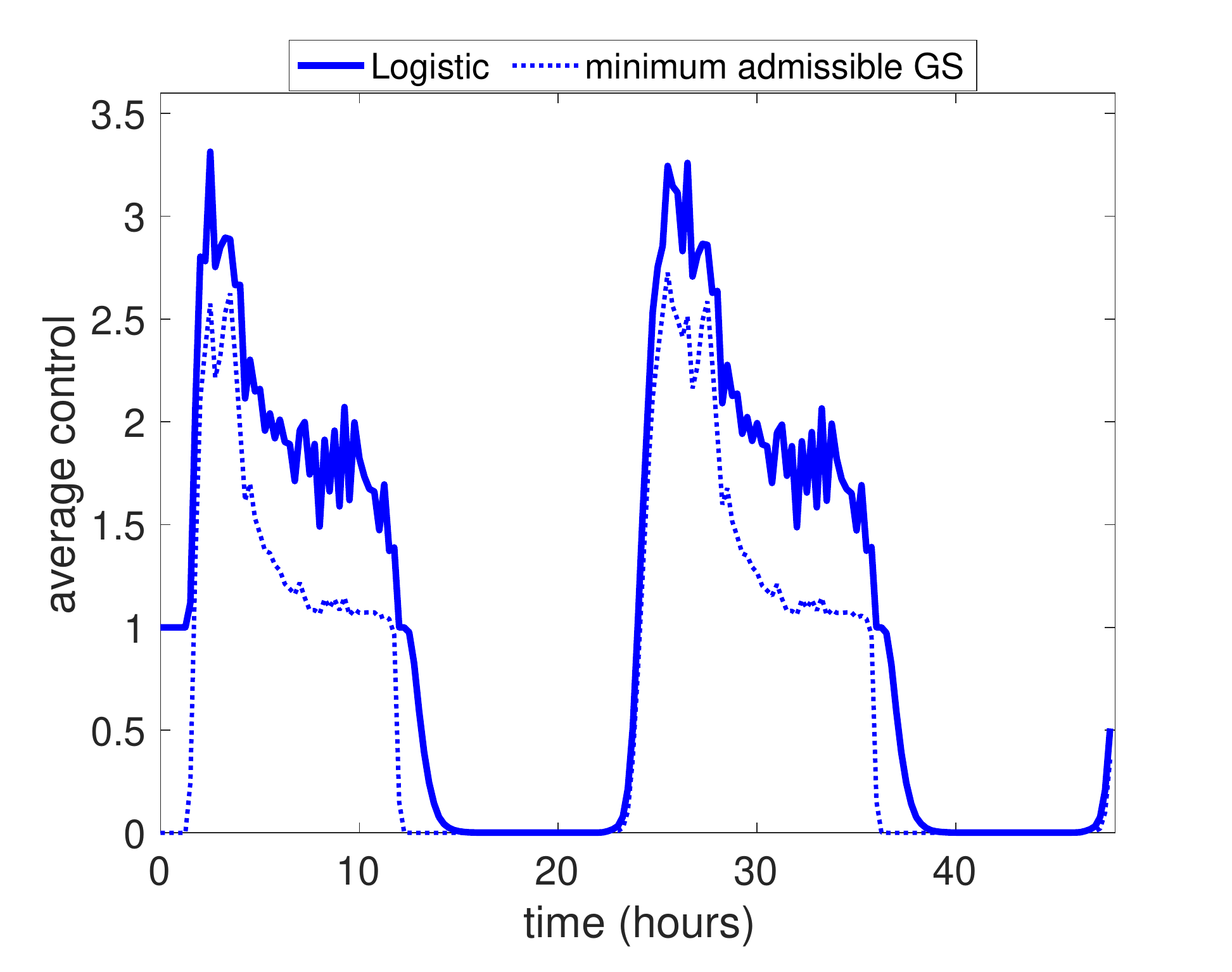}
            \includegraphics[width=0.32\textwidth]{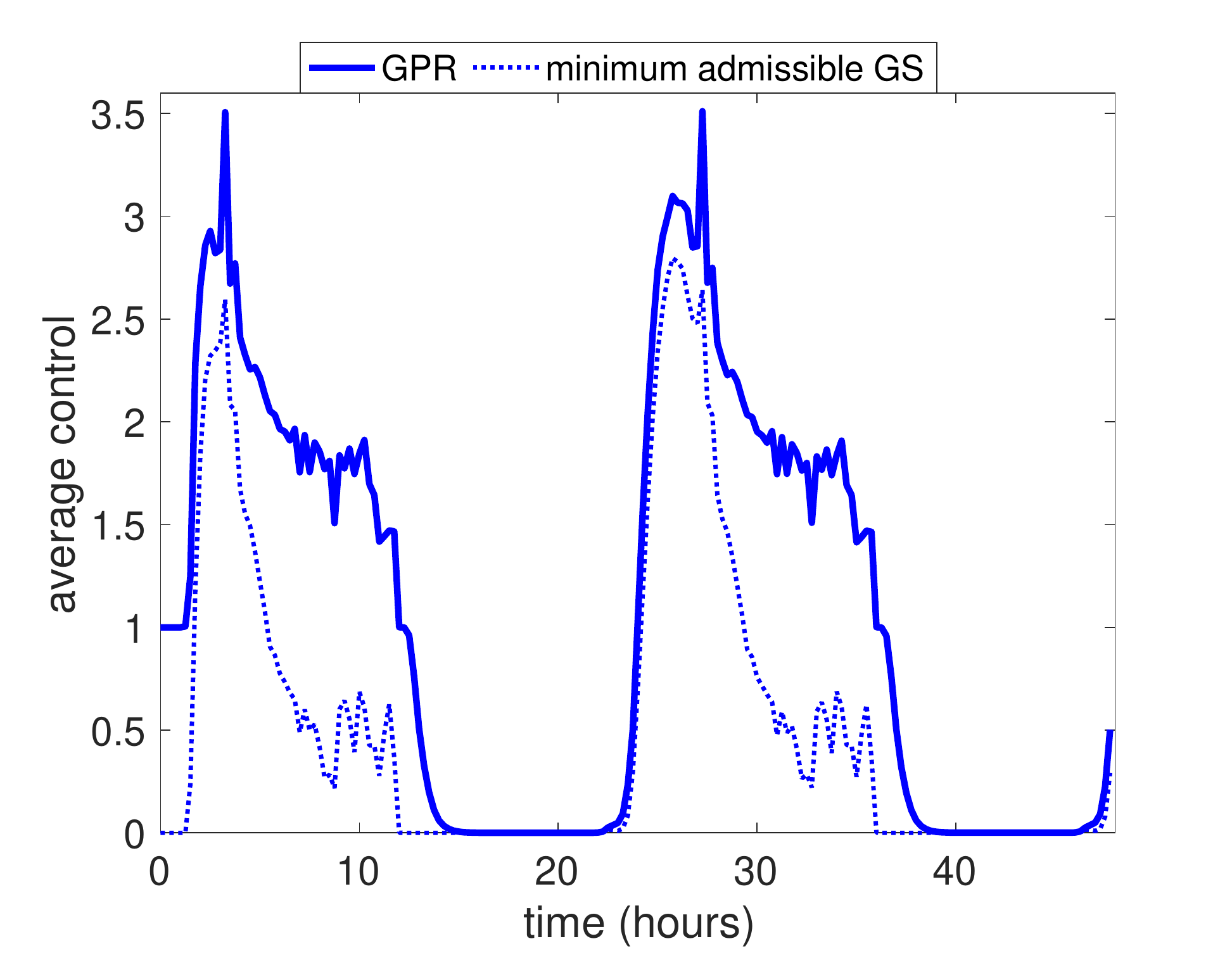}
            \includegraphics[width=0.32\textwidth]{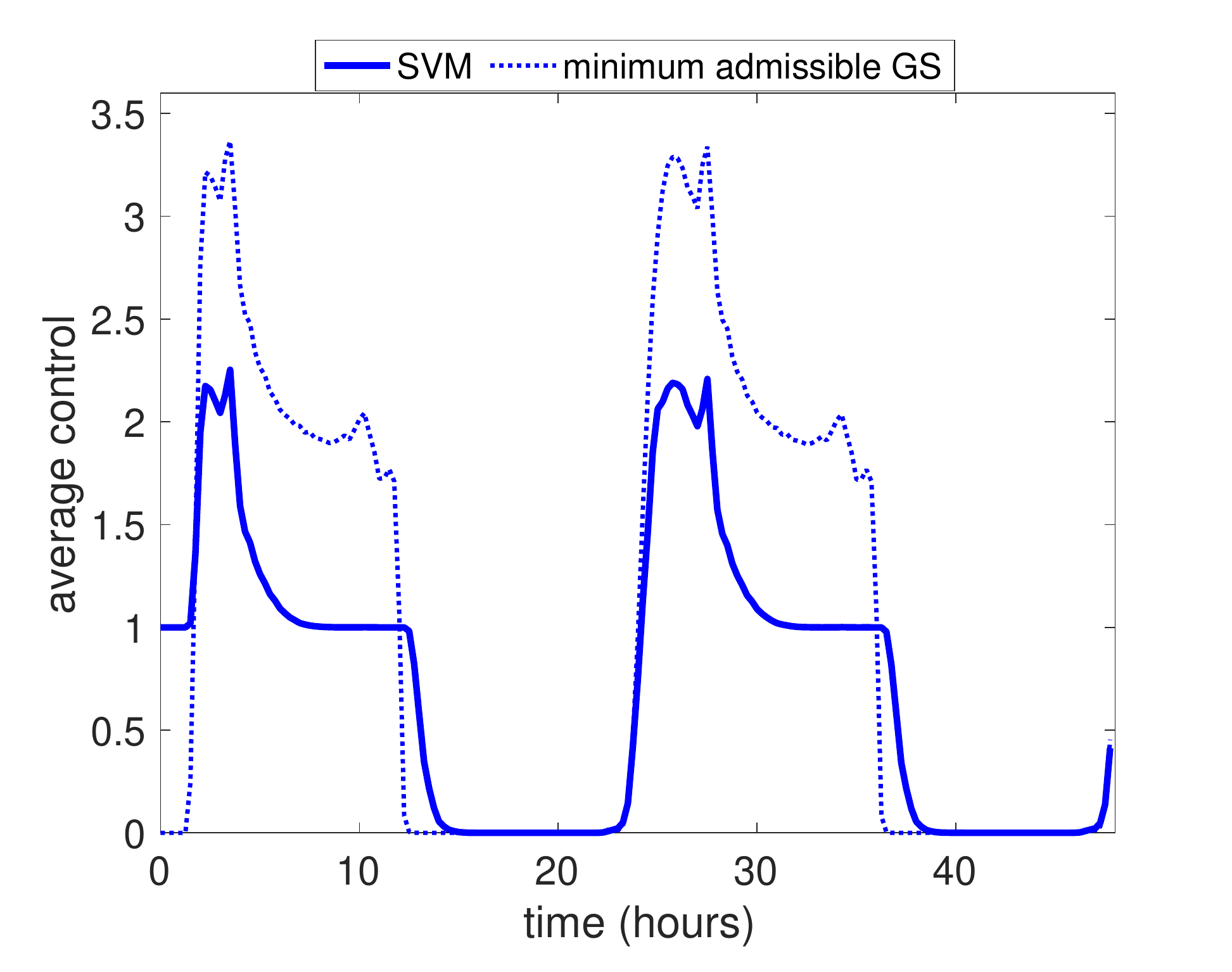}
        \caption[ ]
        {\small{ Average control $Ave(\hat{u}_n)$ for LR, GPR and SVM and the average minimum admissible control $Ave(u^{\min}_n)$ using Gold Standard across forward controlled trajectories.  }}
        \label{fig:average_inventory}
    \end{figure*}

To illustrate the typical behavior over a trajectory, Figure~\ref{fig:average_inventory} plots the average control $Ave({\hat{u}_n}):= \frac{1}{M'}\sum_{m=1}^{M'}\hat{u}_n(\bx_n^{\hat{u},m'})$ corresponding to different methods and the average minimum admissible control $Ave({u^{\min}_n}) := \frac{1}{M'}\sum_{m'=1}^{M'}u^{\min}_n(\bx_n^{\hat{u},m'})$ computed using the gold standard. Notice that the latter is dependent upon the controlled trajectories $\bx_n^{\hat{u}}$, resulting in different  $Ave({u^{\min}_n})$ across methods. We expect $Ave({\hat{u}_n})$ above $Ave({u^{\min}_n})$ if a given method does not violate the constraint most of the time. This is true for  LR and GPR, but SVM quite obviously fails, as the dashed line in the leftmost panel of Figure~\ref{fig:average_inventory} is significantly higher than the solid line at numerous time steps. Furthermore, the conservative nature of GPR is reflected in the large difference between the average minimum admissible control and the average optimal control. This is also evident through $w_{bind} \approx 30\%$ for GPR compared to approximately 45\% for the rest of the methods.

\section{Conclusion}

We developed a statistical learning framework to solve stochastic optimal control problems with local probabilistic constraints.  The key objective of our algorithm is to efficiently estimate the set of admissible controls $\cU(\cdot)$ and the continuation value function $\cC(\cdot,\cdot)$ covering a general formulation of the state process dynamics and rewards. Since SCPC problems require estimating the admissible set repeatedly during the backward induction, we use regression based functional representation of $\bx \mapsto \cU(\bx)$. This perspective also provides a natural way of uncertainty quantification for admissibility, in particular offering conservative estimates that bring statistical guarantees regarding $\hat{\cU}$. At the same time, our dynamic emulation algorithm allows parallel computation of $\cU$ and $\cC$ for additional computational efficiency.

Thanks to the plug-and-play functionality of the dynamic emulation algorithm, it was straightforward to test a large variety of schemes for learning $\cU$. Our numerical results suggest that estimating probabilistic constraints via logistic regression, Gaussian process smoothing and parametric density fitting is more accurate than estimating the corresponding quantile (empirical ranking, SVM or quantile regression). A future line of research would be to additionally parametrize (e.g.~using another GP model) the optimal control map $\bx \mapsto \hat{u}_n(\bx)$ \cite{gpdp} which would speed-up the algorithm in the context of continuous action spaces. \blu{Another direction would be to extend the one-dimensional control framework described in this paper to a multi-dimensional setting. In reference to the microgrid example, multi-dimensional control will allow us to control the diesel output and the demand response or to control output from multiple dispatchable generators.} 

\vspace{10pt}


\appendix
\section{Effect of Simulation Budget $(M_b \times M_a)$}
\label{sec:appendix_simulation_budget}
\begin{table}[H]
\small
\caption{ \footnotesize Impact of simulation budget on performance of SVM for the case study in Section~\ref{sec:example1} and probability thresholds $p=5\%$ and $p=1\%$. The reported values are averages over 10 runs of each scheme. The total simulation budget is divided into batch size $M_b$ and number of design sites $M_a$. For total budget $10^5$: $(M_b,M_a)=(100,1000)$; for $10^6$: $(M_b,M_a)=(500,2000)$; for $10^7$: $(M_b,M_a)=(2000,5000)$; for $10^8$: $(M_b,M_a)=(10000,10000)$.}
\label{tab:svm_budget}
\centering
\begin{tabular}{@{}lcccccrr@{}}
\toprule
$p$ & Budget & $\hat{V}_0(0,5)$ (\$) & $w_{freq}$ (\%) & $w_{avm}$ (kW) & $w_{rlzd}$ (\%) & $\tilde{\cT}$ & $w_{bind}$ (\%) \\ \midrule
\multirow{4}{*}{5\%} & $10^5$ & 26.38 & 5.93 & 0.78 & 2.80 & 665 & 9.73 \\
 & $10^6$ & 26.55 & 5.28 & 0.55 & 1.84 & 386 & 9.77 \\
 & $10^7$ & 26.68 & 4.96 & 0.53 & 1.64 & 330 & 9.75 \\
 & $10^8$ & 26.79 & 1.50 & 0.27 & 0.30 & -51 & 9.22 \\ \midrule
\multirow{4}{*}{1\%} & $10^5$ & 28.32 & 6.63 & 0.93 & 2.43 & 1,460 & 9.87 \\
 & $10^6$ & 28.26 & 5.17 & 0.66 & 1.09 & 631 & 9.56 \\
 & $10^7$ & 28.52 & 0.55 & 0.24 & 0.03 & -39 & 8.78 \\
 & $10^8$ & 28.41 & 0.15 & 0.22 & 0.01 & -51 & 8.82 \\ \bottomrule
\end{tabular}
\end{table}

\section{Effect of Adaptive Margin of Error Level $\rho$}
\label{sec:appendix_rho}
\begin{table}[H]
\small
\centering
\caption{\footnotesize Impact of conservative $\cU^{(\rho)}$ estimators for the case study in Section~\ref{sec:example1}. The probabilistic constraint is set at $p=5\%$. }
\label{tab:example1_conservative1}
{%
\begin{tabular}{@{}lccc|ccc|ccc@{}}
\toprule
 & \multicolumn{3}{c|}{$\rho = 90\%$} & \multicolumn{3}{c|}{$\rho = 99\%$}   & \multicolumn{3}{c}{$\rho = 99.95\%$} \\ \midrule
Method & $\hat{V}_0(0,5)$ & $w_{freq}$  & $w_{rlzd}$ & $\hat{V}_0(0,5)$ & $w_{freq}$ &  $w_{rlzd}$ & $\hat{V}_0(0,5)$ & $w_{freq}$  &  $w_{rlzd}$ \\ \midrule
LR & 26.74 & 0.090 & 0.034 & 26.87 & 0.085 & 0.032 & 27.04 & 0.085 & 0.026 \\
GPR & 27.34 & 0.012 & 0.055 & 28.06 & 0.007 & 0.037 & 28.06 & 0.005 & 0.029 \\
SVM & 27.35 & 4.975 & 1.732 & 29.20 & 3.481 & 1.117 & 29.72 & 3.395 & 1.088 \\
QR & 27.20 & 5.373 & 0.793 & 27.18 & 4.880 & 0.676 & 27.04 & 4.409 & 0.591 \\
CTE & 27.93 & 7.158 & 1.153 & 28.31 & 6.766 & 0.888 & 28.61 & 6.163 & 0.714 \\
EP & 26.78 & 7.990 & 1.497 & 27.17 & 7.629 & 1.183 & 27.96 & 7.102 & 0.956 \\ \bottomrule
\end{tabular}
}
\end{table}

\clearpage
\bibliographystyle{siamplain}
\bibliography{references}

\begin{thebibliography}{10}

\bibitem{ShapiroReview}
{\sc S.~Ahmed and A.~Shapiro}, {\em Solving chance-constrained stochastic
  programs via sampling and integer programming}, in INFORMS TutORials in
  Operations Research, Z.-L. Chen and S.~Raghavan, eds., INFORMS, 2014,
  pp.~261--269.

\bibitem{delara17}
{\sc J.-C. Alais, P.~Carpentier, and M.~{De Lara}}, {\em Multi-usage hydropower
  single dam management: chance-constrained optimization and stochastic
  viability}, Energy Systems, 8 (2017), pp.~7--30.

\bibitem{Andrieu2010}
{\sc L.~Andrieu, R.~Henrion, and W.~R{\"o}misch}, {\em A model for dynamic
  chance constraints in hydro power reservoir management}, European Journal of
  Operational Research, 207 (2010), pp.~579--589.

\bibitem{balata17}
{\sc A.~{Balata} and J.~{Palczewski}}, {\em {Regress-Later {M}onte {C}arlo for
  Optimal Inventory Control with applications in energy}}, arXiv:1703.06461,
  (2017).

\bibitem{balata18}
{\sc A.~{Balata} and J.~{Palczewski}}, {\em {Regress-Later Monte Carlo for
  optimal control of Markov processes}}, arXiv:1703.09705,  (2018).

\bibitem{Blackmore2010}
{\sc L.~{Blackmore}, M.~{Ono}, A.~{Bektassov}, and B.~C. {Williams}}, {\em A
  probabilistic particle-control approximation of chance-constrained stochastic
  predictive control}, IEEE Transactions on Robotics, 26 (2010), pp.~502--517.

\bibitem{Blackmore2011}
{\sc L.~{Blackmore}, M.~{Ono}, and B.~C. {Williams}}, {\em Chance-constrained
  optimal path planning with obstacles}, IEEE Transactions on Robotics, 27
  (2011), pp.~1080--1094.

\bibitem{boogert08}
{\sc A.~Boogert and C.~de~Jong}, {\em Gas storage valuation using a {M}onte
  {C}arlo method}, The Journal of Derivatives, 15 (2008), pp.~81--98.

\bibitem{boogert12}
{\sc A.~Boogert and C.~de~Jong}, {\em Gas storage valuation using a
  multi-factor price process}, Journal of Energy Markets, 4 (2011), pp.~29--52.

\bibitem{warin12}
{\sc B.~Bouchard and X.~Warin}, {\em Monte {C}arlo valuation of {A}merican
  options: Facts and new algorithms to improve existing methods}, in Numerical
  Methods in Finance: Bordeaux, June 2010, R.~A. Carmona, P.~Del~Moral, P.~Hu,
  and N.~Oudjane, eds., Springer Berlin Heidelberg, Berlin, Heidelberg, 2012,
  pp.~215--255.

\bibitem{Calafiore2006}
{\sc G.~C. {Calafiore} and M.~C. {Campi}}, {\em The scenario approach to robust
  control design}, IEEE Transactions on Automatic Control, 51 (2006),
  pp.~742--753.

\bibitem{Campi2009}
{\sc M.~C. {Campi} and G.~C. {Calafiore}}, {\em Notes on the scenario design
  approach}, IEEE Transactions on Automatic Control, 54 (2009), pp.~382--385.

\bibitem{ludkovski10}
{\sc R.~Carmona and M.~Ludkovski}, {\em Valuation of energy storage: an optimal
  switching approach}, Quantitative Finance, 10 (2010), pp.~359--374.

\bibitem{gpdp}
{\sc M.~P. Deisenroth, C.~E. Rasmussen, and J.~Peters}, {\em Gaussian process
  dynamic programming}, Neurocomputing, 72 (2009), pp.~1508 -- 1524.

\bibitem{delara2010}
{\sc L.~Doyen and M.~D. Lara}, {\em Stochastic viability and dynamic
  programming}, Systems \& Control Letters, 59 (2010), pp.~629 -- 634.

\bibitem{Geletu2013}
{\sc A.~Geletu, M.~Klöppel, H.~Zhang, and P.~Li}, {\em Advances and
  applications of chance-constrained approaches to systems optimisation under
  uncertainty}, International Journal of Systems Science, 44 (2013),
  pp.~1209--1232.

\bibitem{heymann16}
{\sc B.~Heymann, J.~F. Bonnans, F.~Silva, and G.~Jimenez}, {\em A stochastic
  continuous time model for microgrid energy management}, in 2016 European
  Control Conference (ECC), June 2016, pp.~2084--2089.

\bibitem{Janson2018}
{\sc L.~Janson, E.~Schmerling, and M.~Pavone}, {\em Monte {C}arlo motion
  planning for robot trajectory optimization under uncertainty}, in Robotics
  Research: Volume 2, A.~Bicchi and W.~Burgard, eds., Springer International
  Publishing, Cham, 2017, pp.~343--361.

\bibitem{Tankov2017}
{\sc Y.~Jiao, O.~Klopfenstein, and P.~Tankov}, {\em Hedging under multiple risk
  constraints}, Finance and Stochastics, 21 (2017), pp.~361--396.

\bibitem{Karatzas1998}
{\sc I.~Karatzas and S.~E. Shreve}, {\em Brownian Motion}, Springer New York,
  New York, NY, 1998.

\bibitem{langrene15}
{\sc N.~Langren{\'e}, T.~Tarnopolskaya, W.~Chen, Z.~Zhu, and M.~Cooksey}, {\em
  New regression {M}onte {C}arlo methods for high-dimensional real options
  problems in minerals industry}, in 21st International Congress on Modelling
  and Simulation, 2015.

\bibitem{langrene17}
{\sc N.~{Langren{\'e}} and X.~{Warin}}, {\em {Fast and stable multivariate
  kernel density estimation by fast sum updating}}, Journal of Computational
  and Graphical Statistics, 28 (2019), pp.~596--608.

\bibitem{Liu2017}
{\sc C.~Liu, X.~Wang, Y.~Zou, H.~Zhang, and W.~Zhang}, {\em A probabilistic
  chance-constrained day-ahead scheduling model for grid-connected microgrid},
  in 2017 North American Power Symposium (NAPS), 2017, pp.~1--6.

\bibitem{ls2001}
{\sc F.~A. Longstaff and E.~S. Schwartz}, {\em Valuing {A}merican options by
  simulation: A simple least-squares approach}, The Review of Financial
  Studies, 14 (2001), pp.~113--147.

\bibitem{ludkovski15}
{\sc M.~{Ludkovski}}, {\em Kriging metamodels and experimental design for
  {B}ermudan option pricing}, Journal of Computational Finance, 22 (2018),
  pp.~37--77.

\bibitem{aditya2018}
{\sc M.~{Ludkovski} and A.~{Maheshwari}}, {\em Simulation methods for
  stochastic storage problems: A statistical learning perspective}, Energy
  Systems, 11 (2020), pp.~377--415.

\bibitem{Luedtke}
{\sc J.~Luedtke and S.~Ahmed}, {\em A sample approximation approach for
  optimization with probabilistic constraints}, SIAM Journal on Optimization,
  19 (2008), pp.~674--699.

\bibitem{Nemirovski2006}
{\sc A.~Nemirovski and A.~Shapiro}, {\em Scenario approximations of chance
  constraints}, in Probabilistic and Randomized Methods for Design under
  Uncertainty, G.~Calafiore and F.~Dabbene, eds., Springer London, London,
  2006, pp.~3--47.

\bibitem{Nemirovski2007}
{\sc A.~Nemirovski and A.~Shapiro}, {\em Convex approximations of chance
  constrained programs}, SIAM Journal on Optimization, 17 (2007), pp.~969--996.

\bibitem{Ono2015}
{\sc M.~Ono, M.~Pavone, Y.~Kuwata, and J.~Balaram}, {\em Chance-constrained
  dynamic programming with application to risk-aware robotic space
  exploration}, Autonomous Robots, 39 (2015), pp.~555--571.

\bibitem{Alejandra2019}
{\sc A.~Peña-Ordieres, J.~R. Luedtke, and A.~Wächter}, {\em {Solving
  chance-constrained problems via a smooth sample-based nonlinear
  approximation}}, arXiv:1905.07377,  (2019).

\bibitem{Prekopa1978}
{\sc A.~Pr{\'e}kopa and T.~Sz{\'a}ntai}, {\em Flood control reservoir system
  design using stochastic programming}, in Mathematical Programming in Use,
  M.~L. Balinski and C.~Lemarechal, eds., Springer Berlin Heidelberg, Berlin,
  Heidelberg, 1978, pp.~138--151.

\bibitem{Ludkovski2016}
{\sc S.~A.~P. Quintero, M.~Ludkovski, and J.~P. Hespanha}, {\em Stochastic
  optimal coordination of small {UAV}s for target tracking using
  regression-based dynamic programming}, Journal of Intelligent {\&} Robotic
  Systems, 82 (2016), pp.~135--162.

\bibitem{adhithya18}
{\sc A.~Ravichandran, S.~Sirouspour, P.~Malysz, and A.~Emadi}, {\em A
  chance-constraints-based control strategy for microgrids with energy storage
  and integrated electric vehicles}, IEEE Transactions on Smart Grid, 9 (2018),
  pp.~346--359.

\bibitem{DiceKriging}
{\sc O.~Roustant, D.~Ginsbourger, and Y.~Deville}, {\em {DiceKriging,
  DiceOptim}: Two {R} packages for the analysis of computer experiments by
  kriging-based metamodeling and optimization}, Journal of Statistical
  Software, 51 (2012), pp.~1--55.

\bibitem{tvr}
{\sc J.~N. Tsitsiklis and B.~van Roy}, {\em Regression methods for pricing
  complex {A}merican-style options}, IEEE Transactions on Neural Networks, 12
  (2001), pp.~694--703.

\bibitem{vanAckooij2014}
{\sc W.~van Ackooij, R.~Henrion, A.~M{\"o}ller, and R.~Zorgati}, {\em Joint
  chance constrained programming for hydro reservoir management}, Optimization
  and Engineering, 15 (2014), pp.~509--531.

\bibitem{warin17}
{\sc W.~{Van-Ackooij} and X.~{Warin}}, {\em {On conditional cuts for Stochastic
  Dual Dynamic Programming}}, ArXiv e-prints,  (2017).

\bibitem{warin}
{\sc X.~Warin}, {\em Gas storage hedging}, in Numerical Methods in Finance:
  Bordeaux, June 2010, R.~A. Carmona, P.~Del~Moral, P.~Hu, and N.~Oudjane,
  eds., Springer, Berlin, Heidelberg, 2012, pp.~421--445.

\bibitem{deltaMethod}
{\sc J.~Xu and J.~S. Long}, {\em Confidence intervals for predicted outcomes in
  regression models for categorical outcomes}, Stata Journal, 5 (2005),
  pp.~537--559.

\end{thebibliography}

\end{document}